\documentclass[11pt,twoside,reqno]{amsart}
\usepackage{import}
\usepackage{newclude}
\usepackage{anysize}
\usepackage{amsmath}
\usepackage{amsthm}
\usepackage{amsmath,amscd}
\usepackage[utf8]{inputenc}
\usepackage{amssymb}
\usepackage{scalerel}
\usepackage{stmaryrd}
\usepackage{wasysym}
\usepackage{mathrsfs}
\usepackage{booktabs}
\usepackage{siunitx} 
\usepackage[pagebackref=true,hypertexnames=false]{hyperref}
\renewcommand*{\backref}[1]{}
\renewcommand*{\backrefalt}[4]{({\tiny%
   \ifcase #1 Not cited.%
         \or Cited on page~#2.%
         \else Cited on pages #2.%
   \fi%
   })}
\usepackage{cleveref}
\usepackage{enumitem}
\usepackage{dsfont}
\usepackage{soul}
\hypersetup{colorlinks=true,linkcolor=black,citecolor=black}
\usepackage{color}
\usepackage[all]{xy}
\usepackage{float}

\usepackage{bm}
\usepackage{mathtools}
\usepackage{makecell}
\usepackage{thmtools}
\usepackage{setspace}
\usepackage{comment}
\usepackage{tikz}
\usepackage{tikz-cd}
\usepackage{mathdots}
\usepackage{graphicx}
\usepackage[myheadings]{fullpage}
\usepackage{csquotes}

\DeclareFontFamily{U}{mathb}{}

\DeclareFontShape{U}{mathb}{m}{n}{
   <-5.5>  mathb5
   <5.5-6.5> mathb6
   <6.5-7.5> mathb7
   <7.5-8.5> mathb8
   <8.5-9.5> mathb9
   <9.5-11>  mathb10
   <11->     mathb12
}{}

\DeclareSymbolFont{mathb}{U}{mathb}{m}{n}

\DeclareMathSymbol{\lefttorightarrow}{3}{mathb}{"FC}
\DeclareMathSymbol{\righttoleftarrow}{3}{mathb}{"FD}

\makeatletter
\newcommand{\acts}{%
  \mathrel{\mathpalette\@acts\righttoleftarrow}}
\newcommand{\actedby}{%
  \mathrel{\mathpalette\@acts\lefttorightarrow}}
\newcommand{\@acts}[2]{\reflectbox{$\m@th#1#2$}}
\makeatother

\numberwithin{equation}{section}
\newcommand\mtop{.95in}
\newcommand\mbottom{.95in}
\newcommand\mleft{1in}
\newcommand\mright{1in}
\usepackage[top = \mtop, bottom = \mbottom, left = \mleft, right=\mright]{geometry}
\newcommand{\val}{{v_p}}
\DeclareMathOperator{\Mat}{Mat}

\newtheorem{thm}{Theorem}[section]
\newtheorem{theorem}[thm]{Theorem}

\newtheorem{prop}[thm]{Proposition}

\newtheorem{lemma}[thm]{Lemma}
\newtheorem{conj}{Conjecture}

\newtheorem{consequence}{Consequence}

\crefname{consequence}{Consequence}{Consequences}
\newtheorem{cor}[thm]{Corollary}
\theoremstyle{definition}
\newtheorem{defi}[thm]{Definition}

\newtheorem{rmk}[thm]{Remark}

\newtheorem{example}[thm]{Example}

\newcommand\reallywidehat[1]{%
\savestack{\tmpbox}{\stretchto{%
  \scaleto{%
    \scalerel*[\widthof{\ensuremath{#1}}]{\kern-.6pt\bigwedge\kern-.6pt}%
    {\rule[-\textheight/2]{1ex}{\textheight}}
  }{\textheight}%
}{0.5ex}}%
\stackon[1pt]{#1}{\tmpbox}%
}
\DeclareSymbolFont{bbold}{U}{bbold}{m}{n}
\DeclareSymbolFontAlphabet{\mathbbold}{bbold}

\makeatletter
\def\@tocline#1#2#3#4#5#6#7{\relax
  \ifnum #1>\c@tocdepth 
  \else
    \par \addpenalty\@secpenalty\addvspace{#2}%
    \begingroup \hyphenpenalty\@M
    \@ifempty{#4}{%
      \@tempdima\csname r@tocindent\number#1\endcsname\relax
    }{%
      \@tempdima#4\relax
    }%
    \parindent\z@ \leftskip#3\relax \advance\leftskip\@tempdima\relax
    \rightskip\@pnumwidth plus4em \parfillskip-\@pnumwidth
    #5\leavevmode\hskip-\@tempdima
      \ifcase #1
       \or\or \hskip 1em \or \hskip 2em \else \hskip 3em \fi%
      #6\nobreak\relax
    \hfill\hbox to\@pnumwidth{\@tocpagenum{#7}}\par
    \nobreak
    \endgroup
  \fi}
\makeatother

\newcommand{\R}{\mathbb{R}}
\newcommand{\Z}{\mathbb{Z}}
\newcommand{\Q}{\mathbb{Q}}

\newcommand{\C}{\mathbb{C}}
\newcommand{\F}{\mathbb{F}}

\newcommand{\E}{\mathbb{E}}

\newcommand{\mf}{\mathfrak}

\newcommand{\mc}{\mathcal}

\newcommand{\bbone}{\mathbbold{1}}

\newcommand{\la}{\lambda}

\renewcommand{\Re}[1]{\text{Re}(#1)}

\newcommand{\inj}{\hookrightarrow}

\renewcommand{\L}{\Lambda}

\renewcommand{\Pr}{\mathbb{P}}

\DeclareMathOperator{\new}{new}
\DeclareMathOperator{\Disc}{Disc}
\DeclareMathOperator{\Nm}{Nm}
\DeclareMathOperator{\Mod}{Mod}
\DeclareMathOperator{\Var}{Var}

\DeclareMathOperator{\Gal}{Gal}

\DeclareMathOperator{\GL}{GL}

\DeclareMathOperator{\SO}{\mathcal{SO}}

\DeclareMathOperator{\GSp}{GSp}

\DeclareMathOperator{\Cov}{Cov}

\DeclareMathOperator{\diag}{diag}

\let\Im\relax
\DeclareMathOperator{\Im}{Im}

\newcommand{\Den}{V}
\newcommand{\bQ}{\bar{\mathbb{Q}}}

\title{Refined random matrix conjectures for zeros of $p$-adic $L$-functions}

\author{Heiko Knospe}
\author{Jiahe Shen}
\author{Roger Van Peski}

\date{\today}
\begin{document}

\begin{abstract}     
    We make and test the conjecture that for $p$-adic $L$-functions associated to quadratic Dirichlet characters and their Teichm\"uller twists, all local statistics of zeros in the unit disk match those of $p$-adic random matrix eigenvalues. This generalizes the case of the number of zeros in the unit disk previously conjectured by Ellenberg-Jain-Venkatesh, which was well-supported by numerics, as well as the case of the discriminant, where their numerics were more ambiguous.
    
    Using recent random matrix results on local statistics of eigenvalues, we make the random matrix prediction explicit for many statistics, analogous to local densities and pair-correlations in the complex setting. To test these, we compute data on $p$-adic $L$-functions for characters of conductor up to $10^7$, and for certain statistics in windows above $10^{10}$. We find good agreement with the new random matrix predictions for these statistics, as well as for the discriminant. 
\end{abstract}

\maketitle

\tableofcontents

\section{Introduction}\label{sec: Introduction}

This paper concerns a conjecture, stated as \Cref{conj:distinguished_polynomial_intro} below. It informally says that given an odd prime $p$, the zeros of a $p$-adic $L$-function, chosen at random from certain families, behave exactly like the eigenvalues of a large random matrix over $\Z_p$.

The conjecture has not been written formally before, but the idea is not new. A very original work of Ellenberg-Jain-Venkatesh \cite{ellenberg2011modeling} formulated the special case where the statistic in question is the number of zeros in the unit disk, with no information about their positions. Here the numerical data which they computed on the $p$-adic $L$-function side agreed quite well with the random matrix prediction. 

Stronger conjectures, while not stated explicitly, are discussed informally in the same work:

\begin{displayquote}
\emph{``More ambitiously, one can speculate that the distribution of more refined features of the $p$-adic $L$-function can also be predicted by random matrix models; we discuss one example, concerning the spacing of zeroes. Here the numerical evidence is less compelling but not wholly discouraging... }

\emph{The observed data in this range is not clearly in line with the predictions, and we do not know whether slow convergence or the prediction itself is to blame.''} 

\end{displayquote}

The contributions of the present paper are twofold:
\begin{enumerate}
    \item We use recent work of the second and third authors \cite{shen2026eigenvalues} on the random matrix side to make the random matrix prediction very explicit for certain statistics of zeros of $p$-adic $L$-functions, analogous to the local densities and pair-correlation functions in the classical complex case, and
    \item We present numerical tests of the resulting predictions, and the `less compelling' one of \cite{ellenberg2011modeling}, by computing data on much larger sets of $p$-adic $L$-functions than had been done previously. 
\end{enumerate} 

For some statistics the convergence of the $p$-adic $L$-function data to the prediction is quite slow, to the point that we did not believe the conjecture for substantial portions of the time spent on the project. Some particularly stubborn statistics necessitated additional computations optimized to that particular statistic in order to reach $p$-adic $L$-functions of conductor several orders of magnitude higher than our main method allowed. However, the last substantial disagreements disappeared in the end, giving us enough confidence to formally state the general conjecture alluded to above. 

\subsection{The conjecture}

Let $p$ be an odd prime. We consider $p$-adic $L$-functions $L_p(s,\chi)$, where $\chi$ is a fixed Dirichlet character of the form $\theta \omega^i$ with $\theta$ a quadratic Dirichlet character and $\omega$ the Teichm\"uller character. These are $p$-adic analytic functions of a variable $s$ lying in the disk $\{s \in \C_p: |s| < p^{1-\frac{1}{p-1}}\}$; the reader unfamiliar with the notation may consult \Cref{sec:preliminaries} for definitions, or (perhaps better) simply accept that this is a reasonable family of $p$-adic analytic functions of arithmetic relevance and continue.

To the $p$-adic $L$-function $L_p(s,\chi)$ is associated an \emph{Iwasawa power series} $F_\chi \in \Z_p[[T]]$, uniquely determined by the condition 
\begin{equation}
    F_\chi((1+dp)^s-1)=L_p(s,\chi) \text{  for all $s \in \C_p$ with $|s| < p^{1-\frac{1}{p-1}}$}
\end{equation}
where $d$ is the primitive conductor of $\theta$. Following \cite{ellenberg2011modeling}, we treat zeros of $F_\chi$ rather than $L_p(s,\chi)$ as the objects of primary interest, though the relations between the two are discussed in \Cref{sec:preliminaries}. Specifically, we treat the zeros of $F_\chi$ inside the open unit disk $B_1(0) \subset \C_p$. The number of such zeros is always finite, and we let $P_\chi(T) \in \Z_p[T]$ be the monic polynomial with these zeros. Because the zeros of $P_\chi$ lie in the open unit disk, $P_\chi$ always lies in one of the sets
$$
\mathscr{P}_d:=\left\{
x^d+a_{d-1}x^{d-1}+\cdots+a_0:
a_0,\ldots,a_{d-1}\in p\Z_p
\right\}
$$
of \emph{distinguished polynomials} of degree $d$ for some $d \geq 0$, where we define $\mathscr{P}_0:=\{1\}$. Rather than considering the roots of $F_\chi$ as a random point configuration in the space $B_1(0)$, which has quite complicated structure, we will state results in the space
$$
\mathscr{P}:=\bigsqcup_{d\geq0}\mathscr{P}_d
$$
in which $P_\chi$ lies. Here we endow each $\mathscr{P}_d$ with the $p$-adic topology from the identification with $(p\Z_p)^d$ and endow $\mathscr{P}$ with the disjoint-union topology and corresponding Borel sigma-algebra. 

Let $A_n$ be a random matrix in $\Mat_n(\Z_p)$ where each entry is independent and is distributed by the additive Haar probability measure on $\Z_p$. Let $P_{A_n}$ be the distinguished part of its characteristic polynomial. Then $P_{A_n}$ converges in distribution as $n \to \infty$ to a random element of $\mathscr{P}$ which we denote $P_{A_\infty}$, meaning in particular that 
\begin{equation}\label{eq:An_weak_cvg_intro}
    \lim_{n \to \infty} \E[\Phi(P_{A_n})] = \E[\Phi(P_{A_\infty})]
\end{equation}
for any continuous and bounded $\Phi: \mathscr{P} \to \R$. This was essentially shown in \cite[Theorem 4.1]{ellenberg2011modeling}.
For the additive Haar-random matrix model used here, the corresponding
weak convergence is proved in \cite[Lemma 4.4]{shen2026resultant}.

On the $L$-function side, we consider $p-1$ different families of $p$-adic $L$-functions. Namely, for each twist $i\in\{0,1,\ldots,p-2\}$, we let $X_i$ denote the family of Dirichlet characters $\chi = \theta \omega^i$ where 
\begin{enumerate}
    \item $\theta$ is a non-trivial quadratic character with primitive conductor coprime to $p$,
    \item $\chi$ is even as a function on $\Z$, 
\end{enumerate}
and if $i=1$ we additionally require that $\theta(p) \neq 1$. We refer to $i$ as the \emph{twist} of the family. Furthermore, for any positive $D$, we let $X_i(D)$ be the finite subfamily of such characters for which $\theta$ has primitive conductor coprime to $p$ and at most $D$.

\begin{conj}
\label{conj:distinguished_polynomial_intro}
Fix a twist $i\in\{0,\ldots,p-2\}$. Then if $\chi_D$ is chosen uniformly at random from $X_i(D)$ and $P_{\chi_D}$ is the corresponding distinguished polynomial, then $P_{\chi_D}$ converges weakly to $P_{A_\infty}$ as $D \to \infty$. In other words, for any continuous bounded function $\Phi: \mathscr{P} \to \R$, 
\begin{equation}
    \label{eq:Lfn_rmt_intro}
    \lim_{D \to \infty} \E[\Phi(P_{\chi_D})] 
=\mathbb{E}\left[\Phi(P_{A_\infty})\right].
\end{equation}
\end{conj}

Here $\E$ denotes the expected value. On the right-hand side of \eqref{eq:Lfn_rmt_intro} it is over the $\mathscr{P}$-valued random polynomial $P_{A_\infty}$. On the left-hand side it is over the uniformly random character $\chi_D$, and this side is thus  equal by definition to the average
$$\lim_{D\to\infty}
\frac{1}{|X_i(D)|}
\sum_{\chi\in X_i(D)}
\Phi\left(P_{\chi}\right).$$
We will similarly use $\E$ for the expectation over $\chi_D$ below, though such expectations may of course be written as sums like the above.

\subsection{Consequences of the conjecture} We emphasize that \Cref{conj:distinguished_polynomial_intro} only gives an explicit, testable prediction for special choices of $\Phi$ where the right-hand side of \eqref{eq:Lfn_rmt_intro} can be computed. For any of these $\Phi$, one may test the prediction of \Cref{conj:distinguished_polynomial_intro} by computing the left-hand side of \eqref{eq:Lfn_rmt_intro} at a fixed large $D$, as we explain in the numerical portion.

Below, in \Cref{thm:lambda_intro} through \Cref{cor:ramified_shell_intro}, we show these predictions for various natural choices of $\Phi$. The proofs of these consequences are carried out in Sections~\ref{sec: Random matrix computations} and~\ref{sec:distributional-consequences}.

Ellenberg-Jain-Venkatesh \cite{ellenberg2011modeling} computed the right-hand side in the special case where $\Phi(P) = \bbone[P \in \mathscr{P}_d]$. Here $\Phi(P_{A_n})$ depends only on $A_n \pmod{p}$, and computations over finite fields produced the following. 

\begin{consequence}[$\lambda$-invariant] \label{thm:lambda_intro} 
Assume \Cref{conj:distinguished_polynomial_intro} and fix a twist $i\in\{0,\ldots,p-2\}$, and let $\chi_D$ be chosen uniformly at random from $X_i(D)$ as there. Then, for every $d\geq 0$,
$$\lim_{D\to\infty} \Pr(\deg(P_{\chi_D})=d) = p^{-d}\prod_{j=d+1}^{\infty}(1-p^{-j}). $$
\end{consequence}

In particular, the limiting distribution of the $\lambda$-invariant is conjecturally the same for every twist $i$. An aggregated version of the cases $i=1$ and $i=1+\frac{p-1}{2}$ of the above is the main conjecture of \cite{ellenberg2011modeling}. In addition to \cite{ellenberg2011modeling}, numerical tests of \Cref{thm:lambda_intro} and analogous conjectures for related families of $p$-adic $L$-functions were carried out by Delbourgo-Chao \cite{delbourgo2015invariants}, Delbourgo-Knospe \cite{delbourgo2023iwasawa}, Knospe \cite{knospe2024special}, Kundu-Washington \cite{kundu2024first,kundu2024heuristics}, and Santato \cite{santato2021zeroes}. The following is also a special case of \Cref{conj:distinguished_polynomial_intro} where the right-hand side was computed in \cite{ellenberg2011modeling}.

\begin{consequence}[Discriminant conditional on $\lambda=2$] \label{thm:disc_ejv_intro} 
Assume \Cref{conj:distinguished_polynomial_intro} and fix a twist $i\in\{0,\ldots,p-2\}$, and let $\chi_D$ be chosen uniformly at random from $X_i(D)$ as there. Then, for every $s\geq 1$, 
\begin{align*}
    \lim_{D\to\infty} \Pr(\val(\Disc(P_{\chi_D}))=s|\deg(P_{\chi_D})=2) &:=\lim_{D\to\infty} \frac{\#\{\chi\in X_i(D): \deg(P_\chi)=2,\ \val(\Disc(P_\chi))=s\}}{\#\{\chi\in X_i(D):\deg(P_\chi)=2\}} \\ 
    &= \begin{cases} p^{1-s}(1-p^{-2})\left(1-p^{-(s+1)/2}\right), & s\text{ odd}\\ p^{1-s}(1-p^{-2})\left(1-\dfrac{p^{-s/2-1}}{1+p^{-1}}\right), & s\text{ even} \end{cases}.
\end{align*}
\end{consequence}

As mentioned above, the numerical agreement with the data in \cite{ellenberg2011modeling} was less compelling. We test \Cref{thm:disc_ejv_intro} for larger conductor ranges $D$ and find that this disagreement between the data and the conjecture largely disappears, see \Cref{tab:lambda2-discriminant-low-high-comparison}.

However, our main interest is not in testing these previously stated special cases, but rather on new statistics capturing the following questions:
\begin{enumerate}
    \item How do the roots of $F_\chi$ typically sort into different extensions $K/\Q_p$?
    \item Given $K$, how are roots distributed within the unit ball in $K$?
\end{enumerate}

Random matrix results from \cite{shen2026eigenvalues}, along with some computations given in \Cref{sec: Random matrix computations}, allow computation of the right-hand side of \eqref{eq:Lfn_rmt_intro} in many new cases bearing on these questions. 

One technical point is that because \Cref{conj:distinguished_polynomial_intro} is made only for continuous bounded $\Phi$, it does not imply averaged statements such as the value---or even existence---of the limit
\begin{equation}
    \lim_{D \to \infty} \E[\deg(P_{\chi_D})].
\end{equation}
This is because by restricting to bounded functions $\Phi$, it neglects the effects of rare events where $P_{\chi_D}$ has very high degree, which can influence such averages. For instance, if somehow for every $n$, the $2^n$th character $\chi$ in $X_i$ (ordered by conductor of the associated $\theta$) had the property that $\deg(P_\chi)=2^n$, this will not affect the limiting frequencies of $\deg(P_\chi)$, but it will have a very large effect on the limiting mean and variance of $\deg(P_\chi)$. Empirically, it seems highly unlikely that a situation like this occurs, and we make an additional conjecture to this effect:

\begin{conj}\label{conj:bad_sets_don't_effect_averages_intro}
    Fix a twist $i \in \{0,\ldots,p-2\}$, and let $B \subset X_i$ be any subset with asymptotically negligible density, 
    \begin{equation*}
        \lim_{D \to \infty} \frac{|B \cap X_i(D)|}{|X_i(D)|} = 0.
    \end{equation*}
    Then, for any fixed $k \geq 0$,
    \begin{equation*}
        \lim_{D \to \infty} \frac{1}{|X_i(D)|} \sum_{\chi \in X_i(D) \cap B}\deg(P_{\chi})^k = 0.
    \end{equation*}
\end{conj}

Together, these conjectures imply that asymptotically negligible families $B \subset X_i$ do not affect limiting statistics. For instance, if $\chi_D$ and $\tilde{\chi}_D$ are characters chosen uniformly at random from $X_i(D)$ and $X_i(D) \cap B^c$ respectively, they imply (see \Cref{thm:moment-convergence-from-conjectures}) the existence and equality of limiting moments
\begin{equation}\label{eq:moments_match_chi_tilde}
     \lim_{D \to \infty} \E[\deg(P_{\chi_D})^k] = \lim_{D \to \infty} \E[ \deg(P_{\tilde{\chi}_D})^k] = \E[\deg(P_{A_\infty})^k].
\end{equation}

In fact, under \Cref{conj:distinguished_polynomial_intro}, the equality of the left- and right-hand sides of \eqref{eq:moments_match_chi_tilde} is \emph{equivalent} to \Cref{conj:bad_sets_don't_effect_averages_intro}, see \Cref{rmk: moment tightness and negligible density sets}. To reiterate, \Cref{conj:distinguished_polynomial_intro} is the main conjecture and gives very fine information, while \Cref{conj:bad_sets_don't_effect_averages_intro} is a much cruder tail-bound type statement which allows one to control averages.

We may now return to the promised predictions for moments of zero counts. For $U \subset \bar{\Z}_p$, we will use shorthand
$$Z_{U,\chi} = \#\{z \in U: F_{\chi}(z)=0\} = \#\{z \in U: P_{\chi}(z)=0\}$$
for the number of zeros lying in $U$. First we consider the simplest case of roots lying in the unit disk $p\Z_p$ of the trivial field extension $\Q_p$. Here we use standard notation
$$\Cov(X,Y)=\E[(X-\E[X])(Y-\E[Y])].$$

\begin{consequence}[Density and pair-correlations of zeros in $\Z_p$]\label{cor:Zp_zeros_intro}
    Assume \Cref{conj:distinguished_polynomial_intro} and \Cref{conj:bad_sets_don't_effect_averages_intro}, fix a twist $i\in\{0,\ldots,p-2\}$, and let $\chi_D$ be chosen uniformly at random from $X_i(D)$. Then for any disk $a+p^r\Z_p \subset p\Z_p$ with $a \in p\Z_p$ and $r \geq 1$, the number of roots $Z_{a+p^r\Z_p,\chi_D}=\#\{z \in a+p^r\Z_p: P_{\chi_D}(z)=0\}$ has limiting mean
    \begin{equation}
        \label{eq:zp_density_intro}
        \lim_{D \to \infty} \E[Z_{a+p^r\Z_p,\chi_D}] = p^{-r}.
    \end{equation}
    For any two disjoint such disks $a_1+p^r\Z_p$ and $a_2+p^r\Z_p$ with $a_1 \not \equiv a_2 \pmod{p^r}$, the numbers of roots have limiting covariance
\begin{equation}\label{eq:covariances_intro}
    \lim_{D \to \infty} \Cov(Z_{a_1+p^r\Z_p,\chi_D},Z_{a_2+p^r\Z_p,\chi_D}) = -p^{-2r}\theta_3\left(-\sqrt{p};\frac{|a_1-a_2|^2}{p}\right),
\end{equation}
    where
    $$
        \theta_3(z;t) = \prod_{i \geq 1} (1-t^i) (1+t^{i-1/2}z) (1+t^{i-1/2}/z) = \sum_{k \in \Z} t^{\frac{k^2}{2}}z^k
    $$
    is a Jacobi theta function.
\end{consequence}

The $p^{-r}$ in \Cref{cor:Zp_zeros_intro} is the Haar measure of the disk $a+p^r\Z_p$. \Cref{cor:Zp_zeros_intro} says informally that when $D$ is large, the zeros of $F_{\chi_D}$ in $p\Z_p$ form a random collection of points of uniform average density $1$ with respect to the Haar measure, and that they repel each other, since the right-hand side of \eqref{eq:covariances_intro} is negative. The latter fact is immediate from the product form for $\theta_3$. 

Another concrete witness to this repulsion is that from \Cref{cor:Zp_zeros_intro} and the summation formula for $\theta_3$ we also have 
\begin{equation}
    \lim_{D \to \infty} \frac{\E[ Z_{a_1+p^r\Z_p,\chi_D} \cdot Z_{a_2+p^r\Z_p,\chi_D}]}{\E[ Z_{a_1+p^r\Z_p,\chi_D}] \cdot \E[ Z_{a_2+p^r\Z_p,\chi_D}]} = (1+p^{-1})|a_1-a_2| + O(|a_1-a_2|^4).
\end{equation}
In other words, when the numbers of zeros in two small disks are normalized to have mean $1$, the expected number of pairs of zeros with one in each disk is linear in the distance, in particular very small at small distance.

\begin{consequence}[The variance in $p^r\Z_p$]
\label{cor:Zp_higher_moments_intro}
With the same assumptions and setting as \Cref{cor:Zp_zeros_intro}, we have
\begin{multline}
\label{eq:zp_variance_intro}
\lim_{D\to\infty}\E\left[
\left(Z_{a+p^r\Z_p,\chi_D}
-\E[Z_{a+p^r\Z_p,\chi_D}]
\right)^2\right]\\
=p^{-r}-p^{-2r}+(1-p^{-1})\sum_{k\ge0}(-1)^k(1+p^{-k-1})\frac{p^{-\frac{k^2+k}{2}}p^{-r(k^2+2k+3)}}{1-p^{-(k^2+2k+2)}}.
\end{multline}
\end{consequence}

We hence have predictions for the first and second moments, and \Cref{cor:Zp_third_moment} below gives very accurate bounds on the random matrix prediction for the third moment.

Note that one should expect that the variance should be the same as the mean to leading order as the disk becomes small, since it is highly unlikely that more than one zero lies in such a disk and so one does not see the repulsion, only Poisson-type statistics. The leading order $p^{-r}$ in the variance above is exactly the mean, as one would expect.

In quadratic extensions we have similar formulas for the average number of zeros in a small disk, but the number of zeros in the whole extension is already small enough that we simply test the number in the whole extension below. Since an extension $K/\Q_p$ of degree $2$ has prime degree, every element of $K\setminus\Q_p$ generates $K$ over $\Q_p$. In these cases, we therefore write $K^{\new}:=K\setminus\Q_p$.

\begin{consequence}[Zeros in quadratic extensions]
\label{cor:quadratic_intro}
With the same assumptions and setting as \Cref{cor:Zp_zeros_intro}, let $K/\Q_p$ be a quadratic extension, and let
$$
B_K:=\{z\in K^{\new}:|z|<1\}.
$$
If $K/\Q_p$ is unramified, then
\begin{align}
\begin{split}
\lim_{D\to\infty}\E[Z_{B_K,\chi_D}]&=\sum_{k\geq 0}(1-p^{-k-1})p^{-\frac{k^2+k}{2}}\\
&\times \frac{(1+p^{-1})p^{-(k^2+2k+3)}(1-p^{-(k^2+2k+3)})-2p^{-1}p^{-(k^2+2k+4)}(1-p^{-(k^2+2k+2)})}
{(1-p^{-(k^2+2k+2)})(1-p^{-(k^2+2k+3)})} \\ 
&=p^{-3}-2p^{-5}+2p^{-6}+p^{-7}+O(p^{-8}).
\end{split}
\end{align}
If $K/\Q_p$ is ramified and $p$ is odd, then 
\begin{align}
\begin{split}
\lim_{D\to\infty}\E[Z_{B_K,\chi_D}]&=p^{-2}(1-p^{-1})\sum_{k\geq 0}\frac{(1-p^{-2k-2})p^{-(k^2+k)}}{(1-p^{-(k^2+2k+2)})(1-p^{-(k^2+2k+3)})}\\
&=p^{-2}-p^{-3}+p^{-4}-p^{-6}+O(p^{-8}).
\end{split}
\end{align}
\end{consequence}
For the centered third moment of the number of zeros in $p^r\Z_p$ and for zeros in cubic extensions, the corresponding random matrix computations do not simplify to exact formulas comparable to those above. Nevertheless, we obtain explicit upper and lower bounds which are sufficiently tight for meaningful comparison with our numerical data. We state and prove these results in \Cref{cor:Zp_third_moment} and \Cref{cor:cubic_extensions} in \Cref{sec:distributional-consequences}.

For extensions of still higher degree, we do not in general have closed forms for the random matrix prediction. There is, however, one particularly simple statistic for totally ramified extensions of arbitrary degree. If $K/\Q_p$ is totally ramified and $\pi$ is a uniformizer of $K$, the expected number of roots in the annulus $\pi\mc{O}_K\setminus\pi^2\mc{O}_K$
admits an explicit prediction. It is natural here to aggregate over all totally ramified extensions $K/\Q_p$ of a fixed degree, which also matches the statistic accessible in our numerical computations.

\begin{consequence}
    \label{cor:ramified_shell_intro}
    Assume \Cref{conj:distinguished_polynomial_intro} and \Cref{conj:bad_sets_don't_effect_averages_intro}, fix an integer $m \geq 1$ and a twist $i\in\{0,\ldots,p-2\}$, and let $\chi_D$ be chosen uniformly at random from $X_i(D)$. Then
    \begin{equation}
        \lim_{D \to \infty} \sum_{K} \E[Z_{\pi\mc{O}_K \setminus \pi^2 \mc{O}_K,\chi_D}] = \frac{m(p-1)}{p^{m+1}},
    \end{equation}
    where the sum ranges over all totally ramified extensions $K/\Q_p$ of degree $m$ with conjugate subfields counted separately, and $\pi$ is a corresponding uniformizer.
\end{consequence}

\begin{rmk}\label{rmk:cohen-lenstra}
The oldest special case of \Cref{conj:distinguished_polynomial_intro} is technically not the conjecture of \cite{ellenberg2011modeling}, but rather the Cohen-Lenstra prediction \cite{cohen-lenstra} for the $p$-part of the class number of quadratic imaginary number fields, though this was arrived at in a different manner. In this case, the constant term of the corresponding $p$-adic $L$-function is (up to $p$-adic unit) the $p$-part of the class number, since $X_i$ excludes the case where it has a trivial zero. See for instance \Cref{thm:base_layer_K_0}. On the random matrix side, the constant term is just the determinant, the limiting distribution of which follows by the computation of Friedman-Washington \cite{friedman-washington} that the size of the cokernel has Cohen-Lenstra distribution. 
\end{rmk}

\subsection{Numerical tests} \label{subsec:numerics_intro}

For each $p \in \{3,5,7\}$, each twist $i \in \{0,\ldots,p-2\}$, and each character $\chi \in X_i(10^7)$, we compute the polynomial $P_\chi$ as described in \Cref{sec:computational_methods}. We may then explicitly compute the left-hand side in any of the above consequences at fixed $D=10^7$. In order to obtain several `independent trials,' as it were, we will instead split $[0,10^7]$ into four equal intervals, yielding four subsets of $X_i(D)$ corresponding to conductors in each interval, and average over each subset. In each case, the random matrix prediction is shown on the plot, and the plot's caption indicates which of the above consequences it corresponds to. Some plots appear here in the Introduction, while additional ones referenced here appear later in \Cref{sec:plots}.

\begin{figure}%
        \centering
        \vspace{5ex}
        \includegraphics[scale=.95]{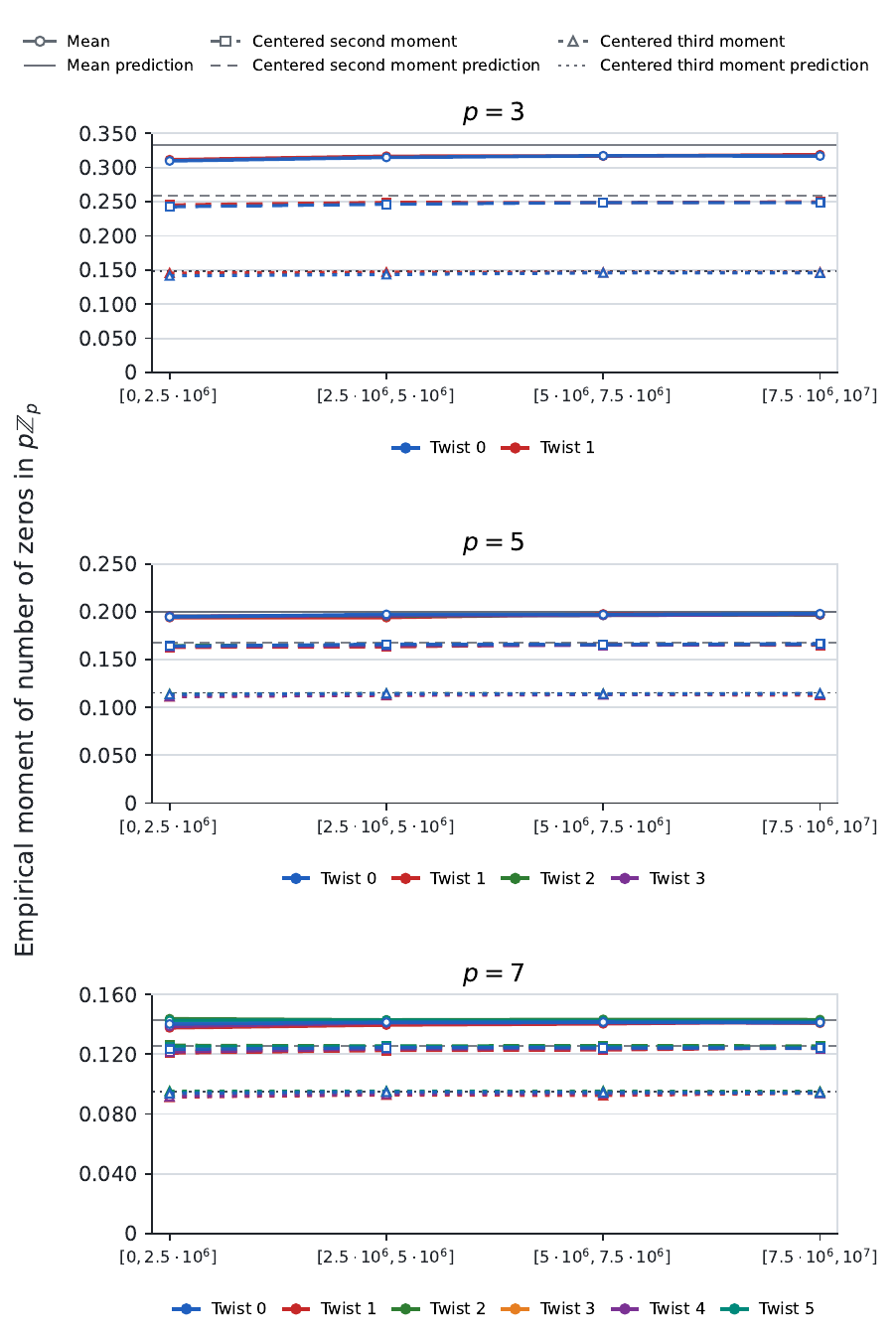}
        \caption{The first three empirical moments of the number of zeros in $p\Z_p$, for four different conductor ranges given on the $x$-axis, the primes $p=3,5,7$, and all twists $0,\ldots,p-2$. Random matrix predictions are from \Cref{cor:Zp_zeros_intro} for the mean, \Cref{cor:Zp_higher_moments_intro} for the second moments, and \Cref{cor:Zp_third_moment} for the third. While we only have bounds for the third moment prediction, their difference is within the width of the prediction line.}%
        \label{fig:Zp_moments_all_primes}
\end{figure}

First consider only the zeros in $\Z_p$, which are the subjects of Consequences~\ref{cor:Zp_zeros_intro} and~\ref{cor:Zp_higher_moments_intro}. These in particular predict the first three moments of the number of zeros in the unit disk $p\Z_p$, shown in \Cref{fig:Zp_moments_all_primes}. The prediction for the third moment was not stated earlier since the formulas are quite long, but is given later in \Cref{cor:Zp_third_moment} in \Cref{sec:distributional-consequences}.

Note that for $p=3$ the empirical averages visually undershoot the prediction. In the Cohen-Lenstra literature it is well-known that convergence is particularly slow for $p=3$, requiring conductor ranges several orders of magnitude higher than what we compute here. It is unsurprising that if the constant term of $P_\chi$ is not too close to the random matrix prediction for conductors less than $10^7$, then other statistics are similarly skewed. We will also see that in the  computations we carry out for conductors $\sim 10^{10}$ and $p=3$, the statistics we compute there match the random matrix prediction much better.

\begin{figure}
        \centering
        \includegraphics[scale=.85]{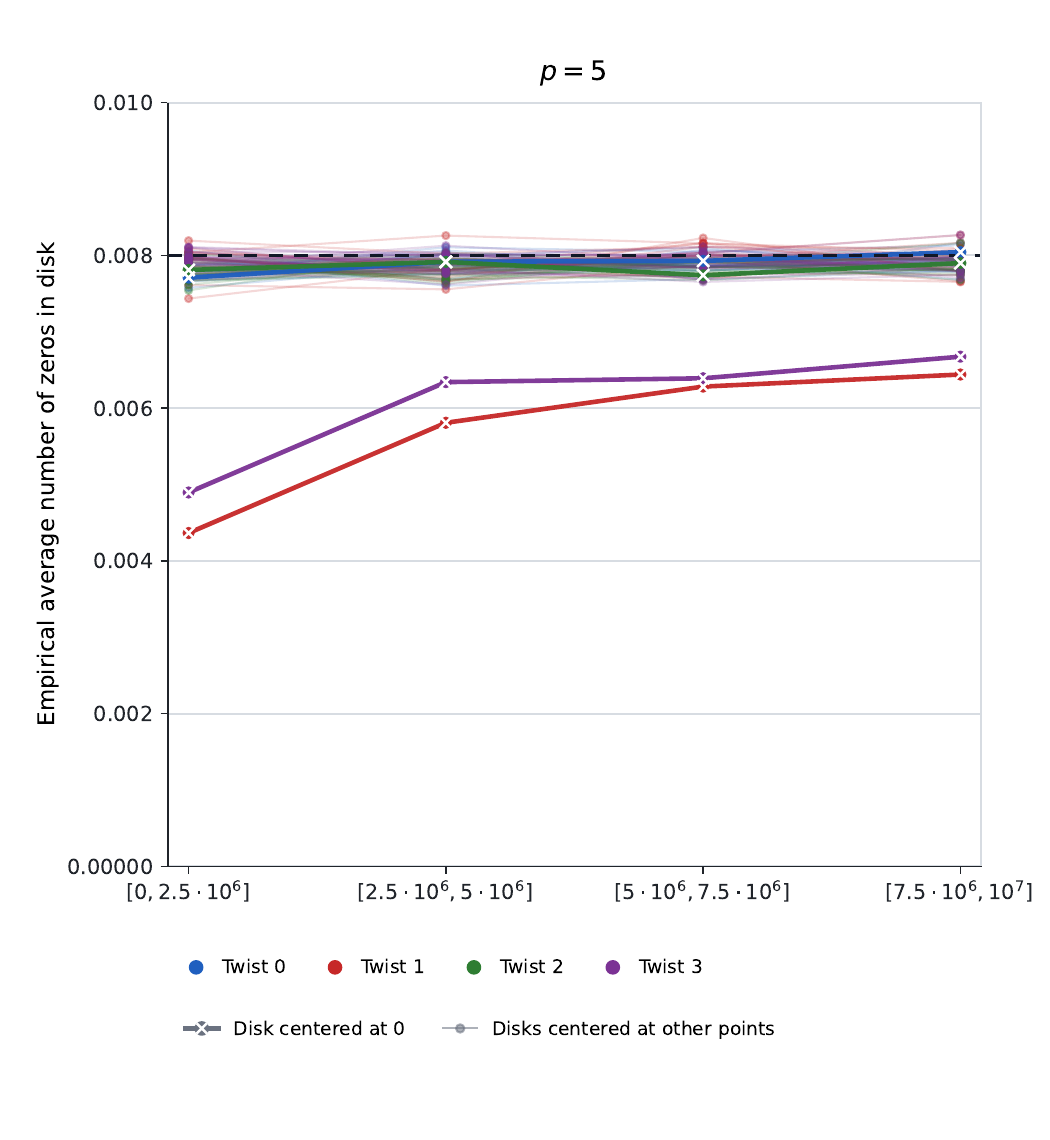}
    \caption{Empirical averages of the number of zeros in small disks in $p\Z_p$ for $p=5$. For each twist, averages are taken over four conductor ranges, shown on the $x$-axis, and the height of a point above a given conductor range is the empirical average number of zeros, while the color of the point corresponds to the twist. Horizontal dotted line shows the random matrix prediction $\frac{1}{p^3}$. Lines are drawn between points corresponding to the same twist and disk at different conductor ranges. The points corresponding to averages over the central disk $p^3\Z_p$ are indicated with crosses and the lines between them are drawn in bold. The other points correspond to averages over shifted disks $ap+p^3\Z_p$, $a \in \{1,\ldots,p^2-1\}$. Note that empirical averages for twists $0,2$ are very close to the random matrix prediction even for the central disk $p^3\Z_p$, while for twists $1$ and $1+\frac{p-1}{2}=3$ the central disk has notably fewer zeros. Analogous plots for $p=3$ and $p=7$ are given in \Cref{fig:p3_small_disk_means} and \Cref{fig:p7_small_disk_means} respectively.}
        \label{fig:p5_small_disk_means}
\end{figure}

\begin{figure}
        \centering
        \vspace{5ex}
        \includegraphics[scale=.95]{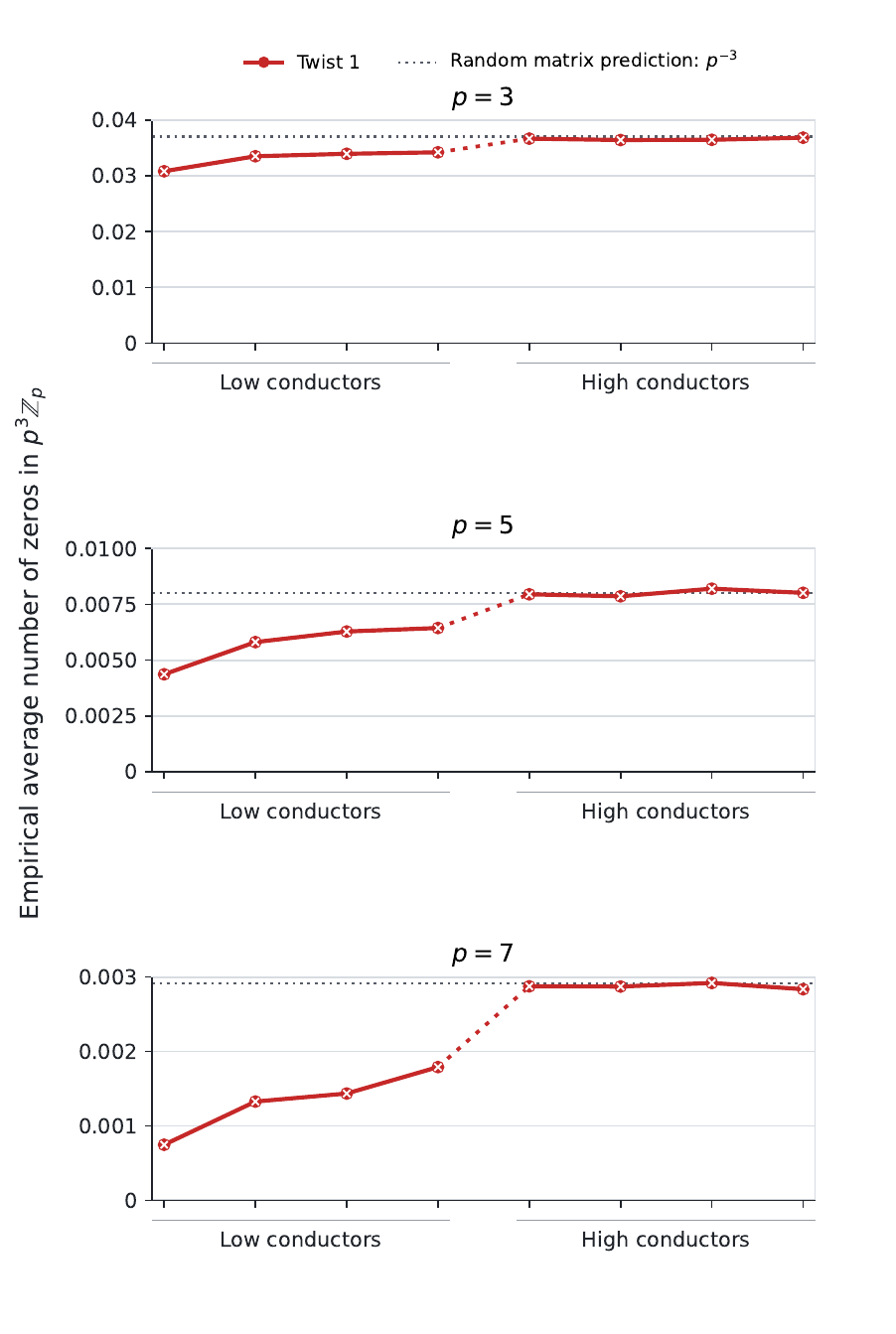}
        \caption{The empirical mean number of zeros in $p^3\Z_p$ for the four conductor ranges $(0,.25\cdot 10^7]$, $(.25\cdot 10^7,.5 \cdot 10^7]$, $(.5\cdot 10^7,.75 \cdot 10^7]$, $(.75 \cdot 10^7,10^7]$ pictured in earlier figures (left), and the empirical means in the translates $(10^{10},10^{10}+.25\cdot 10^7],\ldots$ of these intervals by $10^{10}$ (right).}%
        \label{fig:twist1_smallball_high_low}
\end{figure}

These statistics concern zeros in the entire unit disk, but it is natural to test more local ones as well. \Cref{cor:Zp_zeros_intro} gives predictions for the average numbers of zeros in small disks and the covariances between them.

\begin{figure}[H]
        \centering
        \includegraphics[scale=.85]{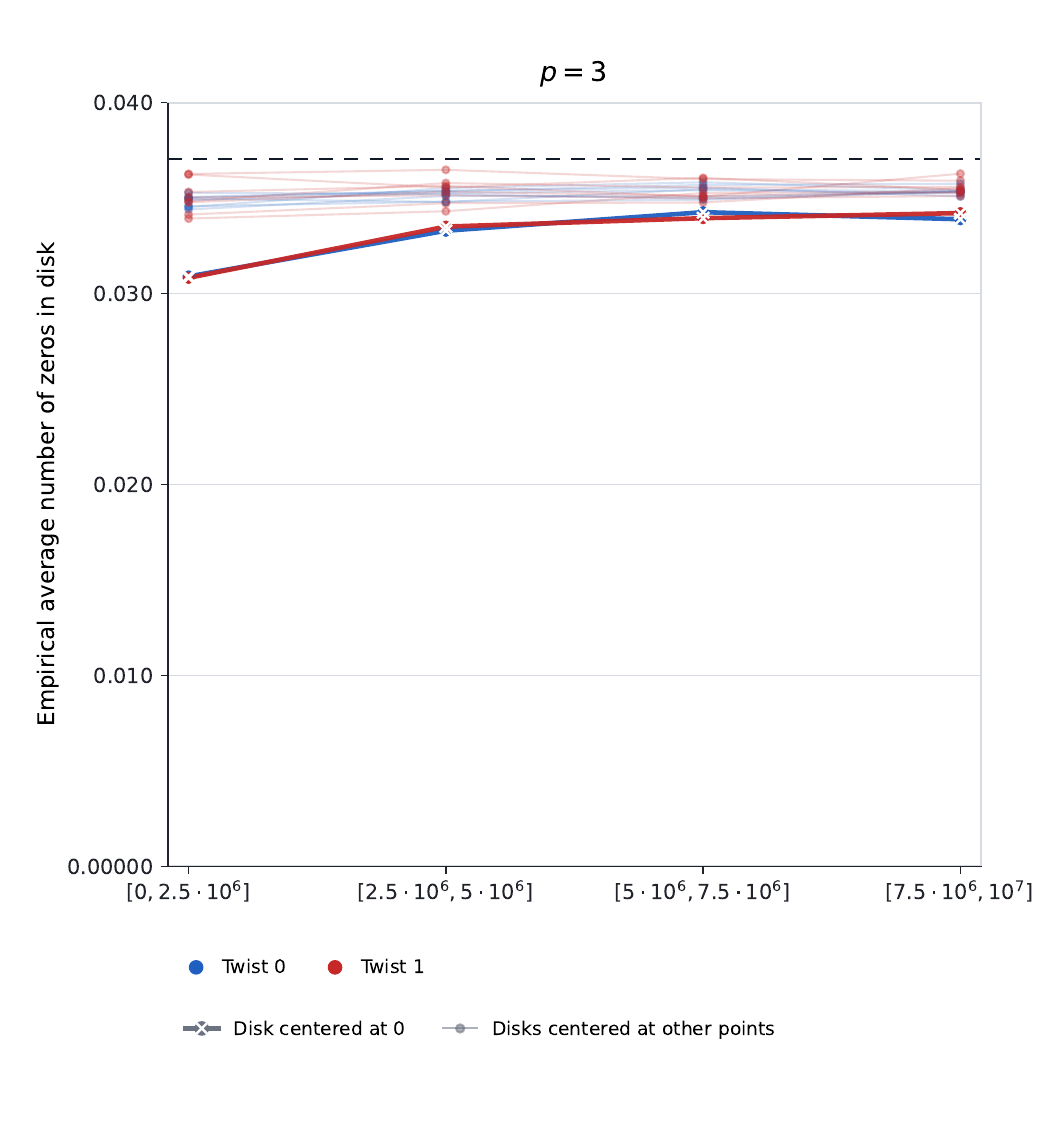}
        \caption{Empirical averages of the number of zeros in small disks in $p\Z_p$ for $p=3$. For each twist, averages are taken over four conductor ranges, shown on the $x$-axis, and the height of a point above a given conductor range is the empirical average number of zeros, while the color of the point corresponds to the twist. Horizontal dotted line shows the random matrix prediction $\frac{1}{p^3}$. Lines are drawn between points corresponding to the same twist and disk at different conductor ranges. The points corresponding to averages over the central disk $p^3\Z_p$ are indicated with crosses and the lines between them are drawn in bold. The other points correspond to averages over shifted disks $ap+p^3\Z_p$, $a \in \{1,\ldots,p^2-1\}$. Here the only twists are $1$ and $0=1+\frac{p-1}{2} \pmod{p-1}$.}
        \label{fig:p3_small_disk_means}
\end{figure}

\begin{figure}[H]
        \centering
\includegraphics[scale=.85]{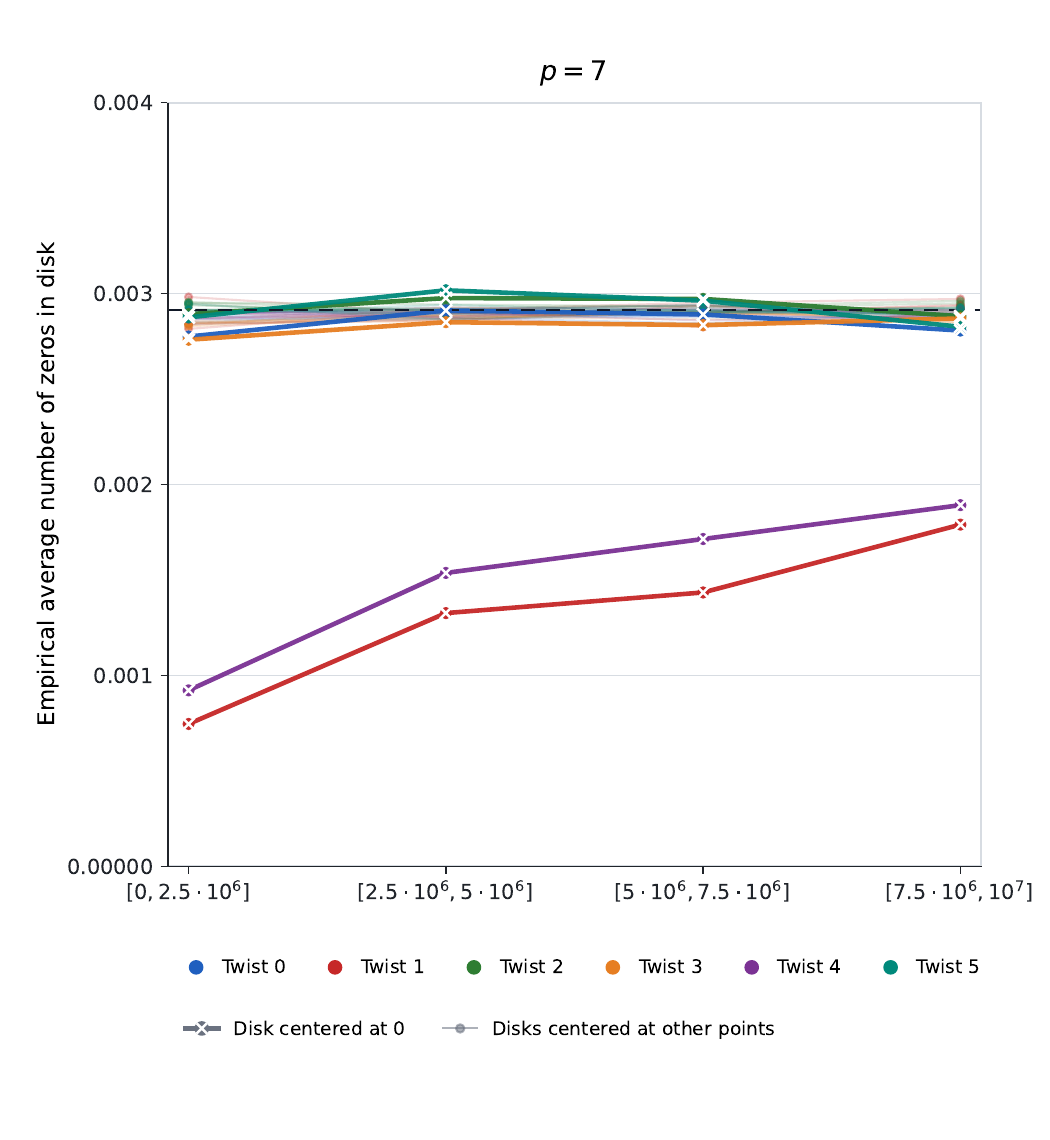}
        \caption{As in \Cref{fig:p3_small_disk_means}, but for $p=7$. As before, points corresponding to averages over the central disk $p^3\Z_p$ are indicated with crosses and the lines between them are drawn in bold. However, because $p^3$ is much larger than in previous plots, we use slightly different averages for the off-center disks. Namely, the other points without crosses correspond to empirical averages in the larger disks $ap+p^2\Z_p$ for $a \in \{1,\ldots,p-1\}$, normalized by $\frac{1}{p}$ so that the random matrix prediction is still $\frac{1}{p^3}$. Each such value is equivalently the average of $p$ different empirical averages in $p$ disjoint shifts of the disk $p^3\Z_p$.}
        \label{fig:p7_small_disk_means}
\end{figure}

The average numbers of zeros in translates of $p^3\Z_p$ are given in \Cref{fig:p5_small_disk_means} for $p=5$, with the analogous figures for $p=3,7$ given later in Figure~\ref{fig:p3_small_disk_means} and~\ref{fig:p7_small_disk_means}. One feature is striking: In each case, the averages are very close to the random matrix prediction, \emph{except} when the disk is centered at $0$ and the twist is $1$ or $1+\frac{p-1}{2}$ (when $p=3$ the latter is twist $0$ since it is taken modulo $p-1=2$). Such disks contain far fewer zeros than predicted. The plots show the average roughly increasing for higher conductors, but it is not at all clear visually whether it will eventually converge to the random matrix prediction or remain bounded away from it. These data initially caused us to doubt that \Cref{conj:distinguished_polynomial_intro} held for such twists. 

However, further optimizing the computations for these specific small-ball averages, we were able to compute them for conductors $D$ of size roughly $10^{10}$. These tell a quite convincing story that the lack of zeros near $0$ only occurs for low conductors and does not affect asymptotics, see \Cref{fig:twist1_smallball_high_low} for the case of twist $i=1$. Due to the number of days of computation required, we ran these only for twist $1$ and the central disk.

\begin{rmk}
    The twists $1$ and $1+\frac{p-1}{2}$ are exactly the ones where the associated character $\chi$ is $\omega$ times a quadratic character. These, in turn, are the ones for which the zeros of $F_{\chi}$ are related to class groups of quadratic imaginary fields, see \Cref{sec:class_numbers}. We do not have any explanation for the origin of this effect, though it may be analogous to the negative secondary terms which are a common and interesting feature in Cohen-Lenstra heuristics, as noted already for the $\lambda$-invariant in \cite{ellenberg2011modeling}. 
\end{rmk}

\begin{figure}[H]
        \centering
        \vspace{5ex}
        \includegraphics[scale=.92]{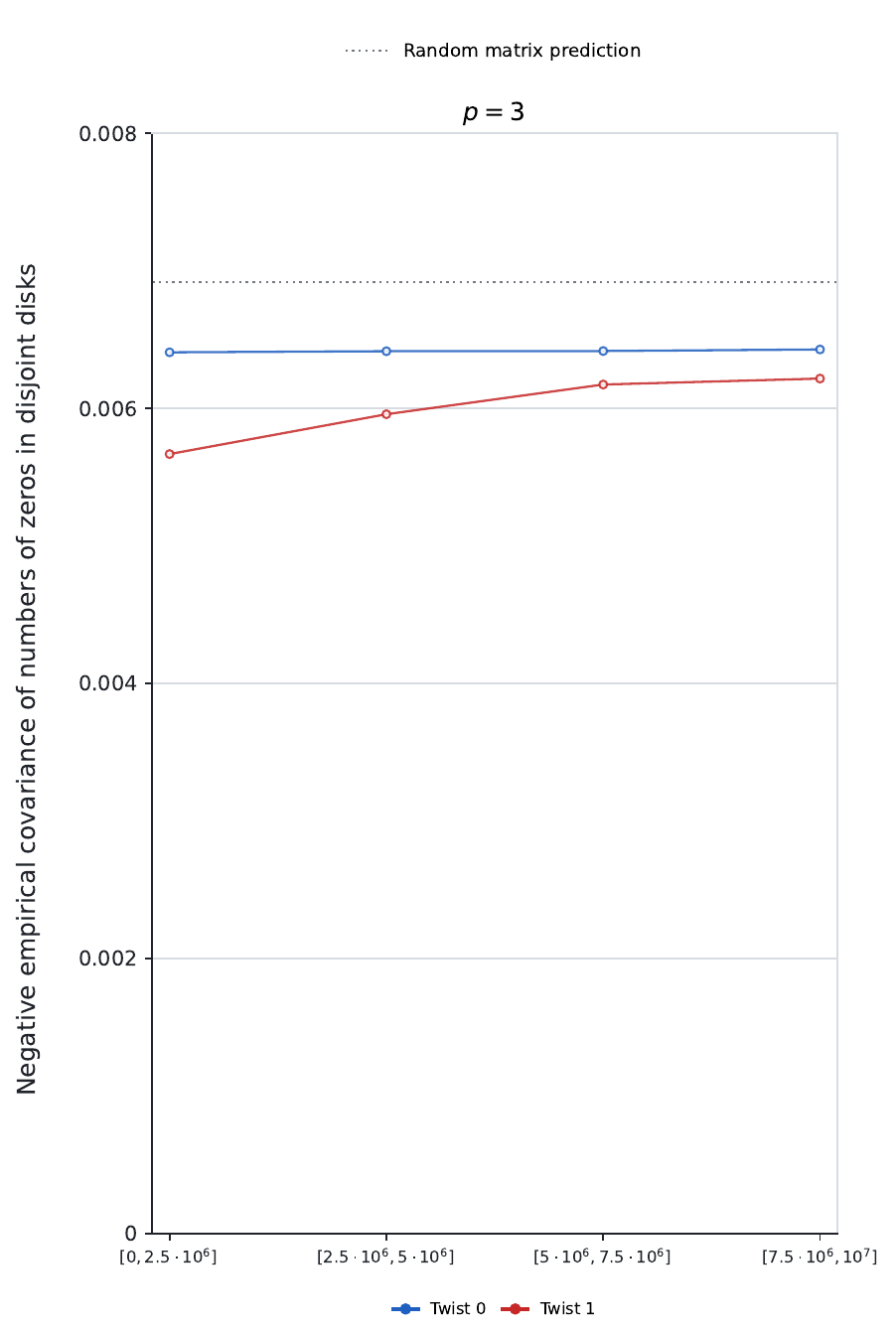}

        \caption{Plot of (negative of) the empirical covariances of the number of roots in $p+p^2\Z_p$ and the number of roots in $2p+p^2\Z_p$, for $p=3$, with random matrix prediction $p^{-4}\theta_3(-\sqrt{p};p^{-3})$ from \Cref{cor:Zp_zeros_intro}.}
        \label{fig:covariances_p3}
\end{figure}

\begin{figure}[H]
        \centering
        \vspace{5ex}
        \includegraphics[scale=.92]{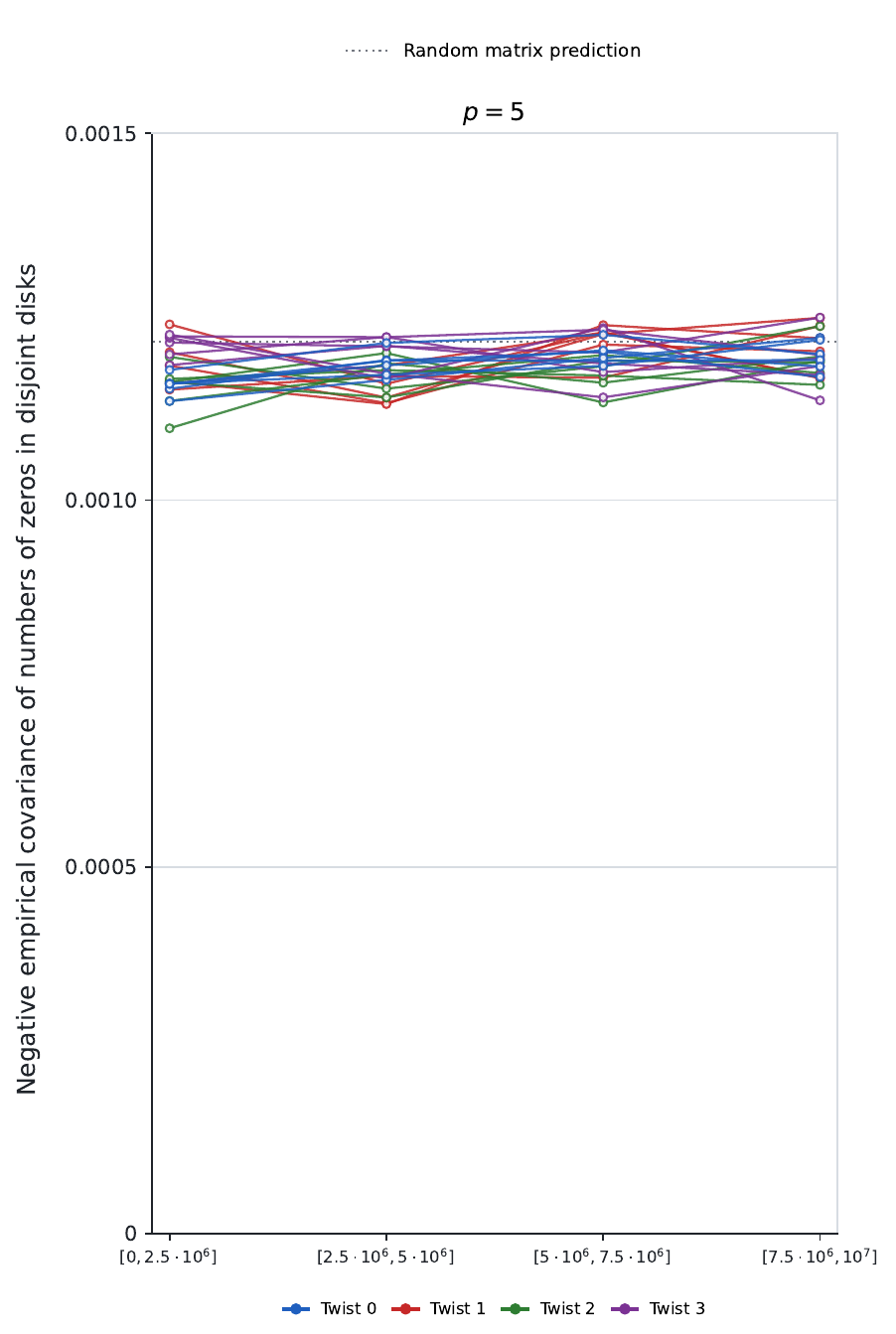}

        \caption{Plot of (negative of) the empirical covariances of the numbers of roots in all $\binom{4}{2}=6$ pairs of disks $ap+p^2\Z_p,a=1,2,3,4$, for $p=5$, with random matrix prediction $p^{-4}\theta_3(-\sqrt{p};p^{-3})$ from \Cref{cor:Zp_zeros_intro}.}
        \label{fig:covariances_p5}
\end{figure}

\begin{figure}[H]
        \centering
        \vspace{5ex}
        \includegraphics[scale=.92]{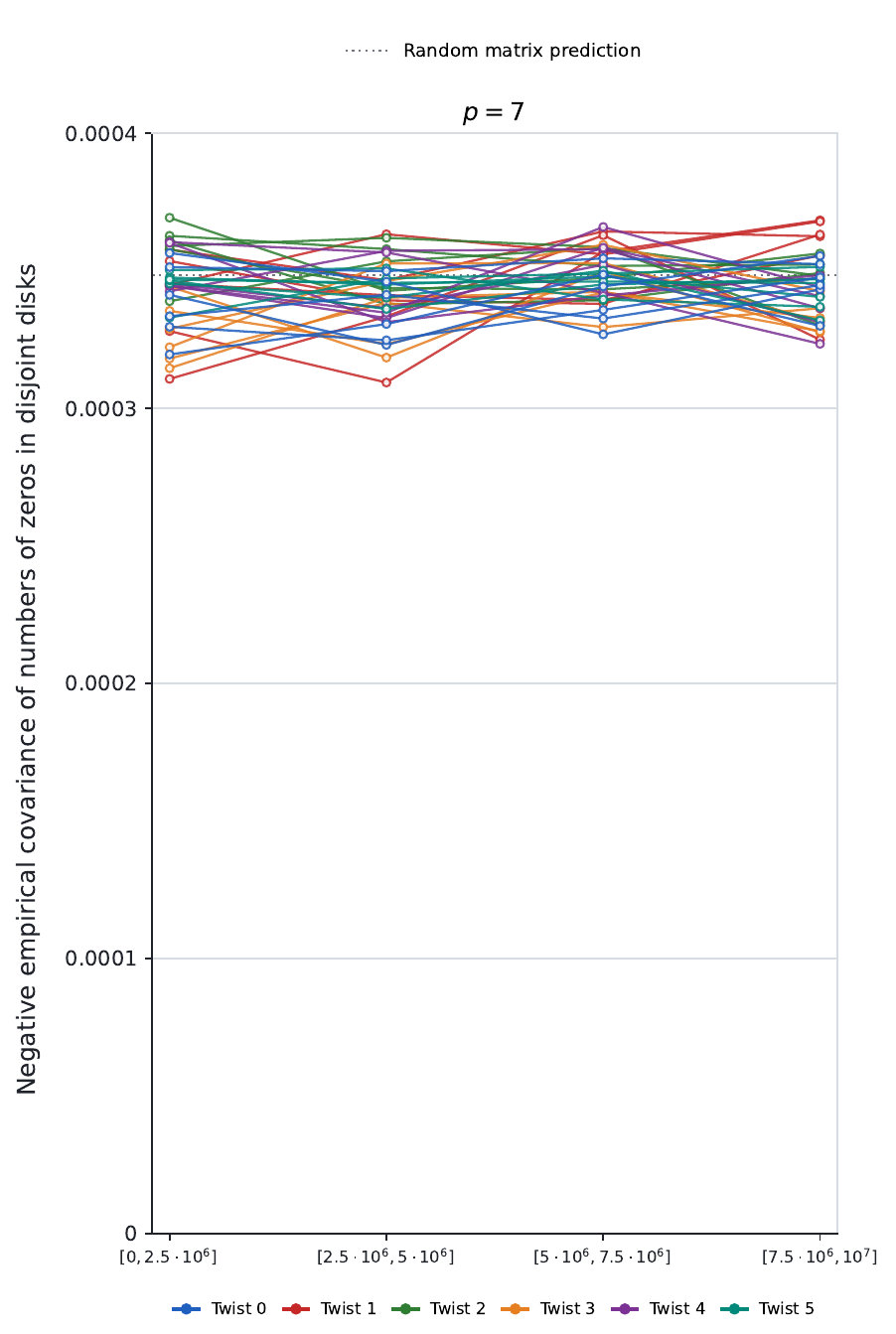}

        \caption{Plot of (negative of) the empirical covariances of the numbers of roots in all $\binom{4}{2}=6$ pairs of disks $ap+p^2\Z_p,a=1,2,5,6$, for $p=7$, with random matrix prediction $p^{-4}\theta_3(-\sqrt{p};p^{-3})$ from \Cref{cor:Zp_zeros_intro}.}
        \label{fig:covariances_p7}
\end{figure}

Tests of the pair-correlation statistics of \Cref{cor:Zp_zeros_intro} are given in Figures~\ref{fig:covariances_p3}, \ref{fig:covariances_p5}, and \ref{fig:covariances_p7}. Note that in each case, we only consider disks which are not centered at the origin; for disks centered at the origin, the anomalous behavior of the mean just discussed will cause similar disagreement with the covariances.

All data discussed thus far is for zeros in $\Z_p$. For nontrivial extensions, our methods do not distinguish which extension a root lies in, only the ramification index of the extension and the $p$-adic valuation of the root. Consequently, we sum the prediction over all extensions of a given degree and ramification. Because the number of such extensions is well known, this gives an equally explicit prediction. In the totally ramified case in \Cref{cor:ramified_shell_intro} the prediction is in fact simpler after summing, which is why we stated it with the sum there.

\begin{figure}[H]
        \centering
        \vspace{5ex}
        \includegraphics[scale=1]{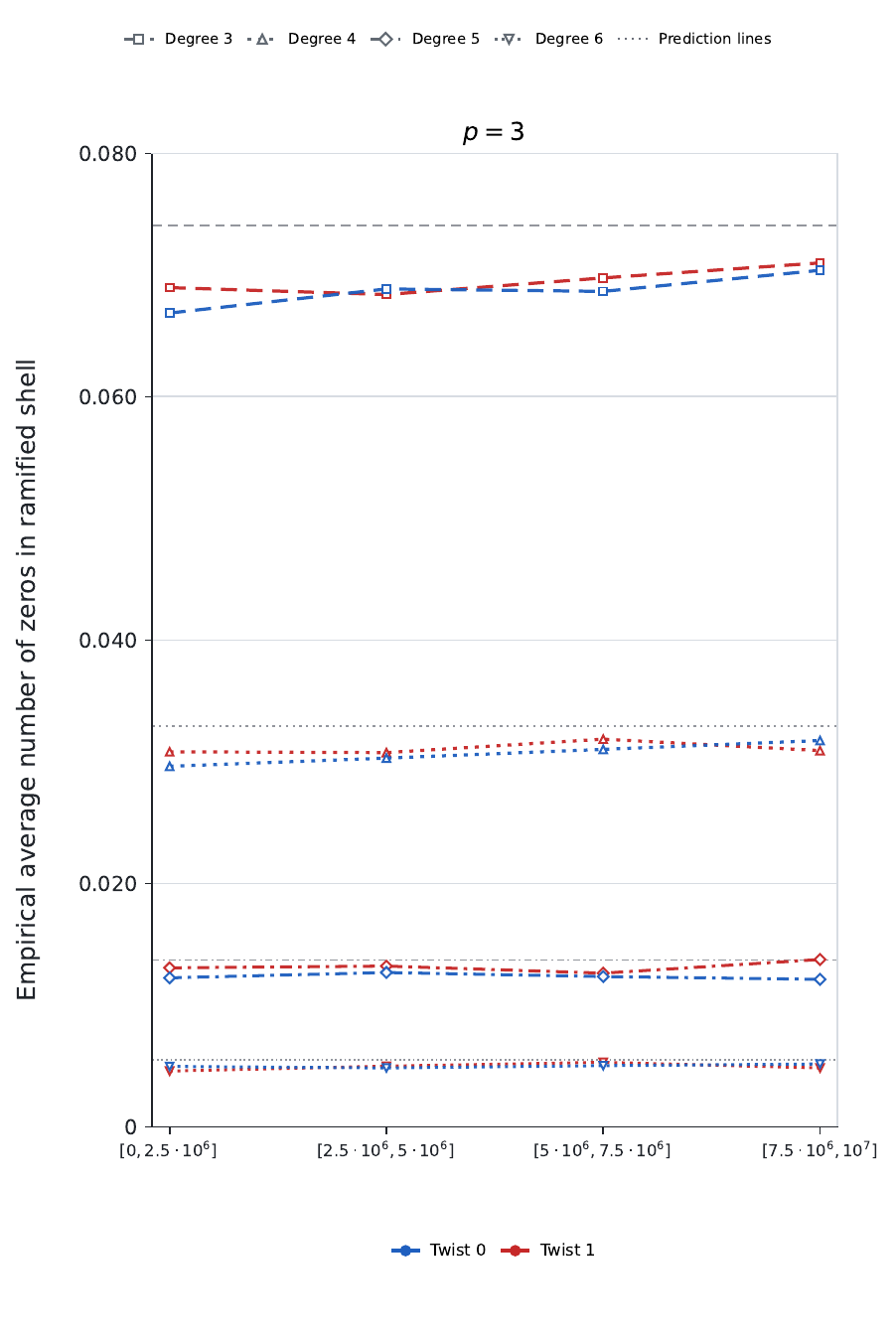}

        \caption{Test of \Cref{cor:ramified_shell_intro}: average number of zeros in $\bigcup_K \pi \mc{O}_K \setminus \pi^2 \mc{O}_K$ where the union is over totally ramified extensions of degree $d=3,4,5,6$, for $p=3$.}
        \label{fig:p3_ramified_shells}
\end{figure}

\begin{figure}[H]
        \centering
        \vspace{5ex}
        \includegraphics[scale=1]{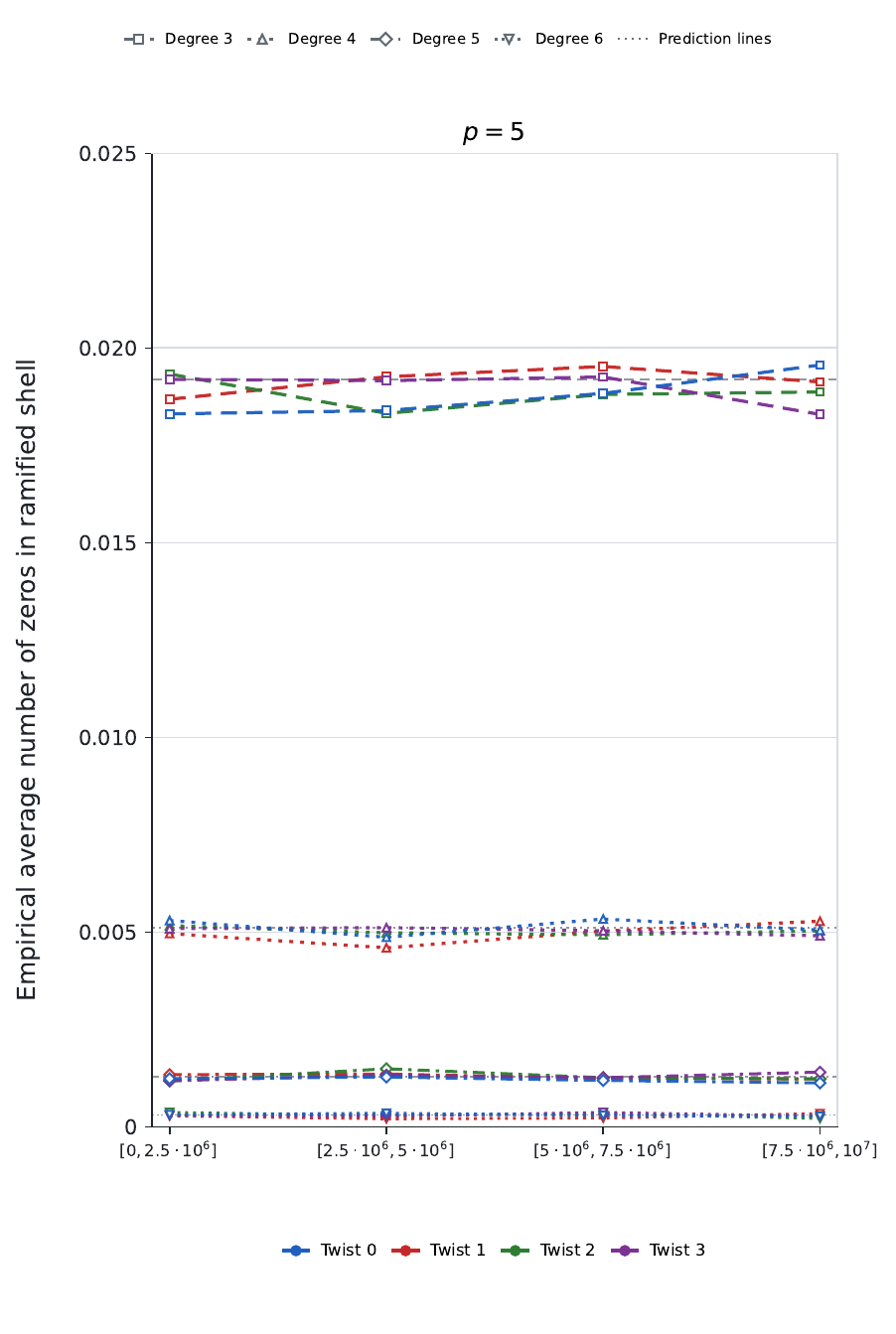}

        \caption{Test of \Cref{cor:ramified_shell_intro}: average number of zeros in $\bigcup_K \pi \mc{O}_K \setminus \pi^2 \mc{O}_K$ where the union is over totally ramified extensions of degree $d=3,4,5,6$, for $p=5$.}
        \label{fig:p5_ramified_shells}
\end{figure}

\begin{figure}[H]
        \centering
        \vspace{5ex}
        \includegraphics[scale=1]{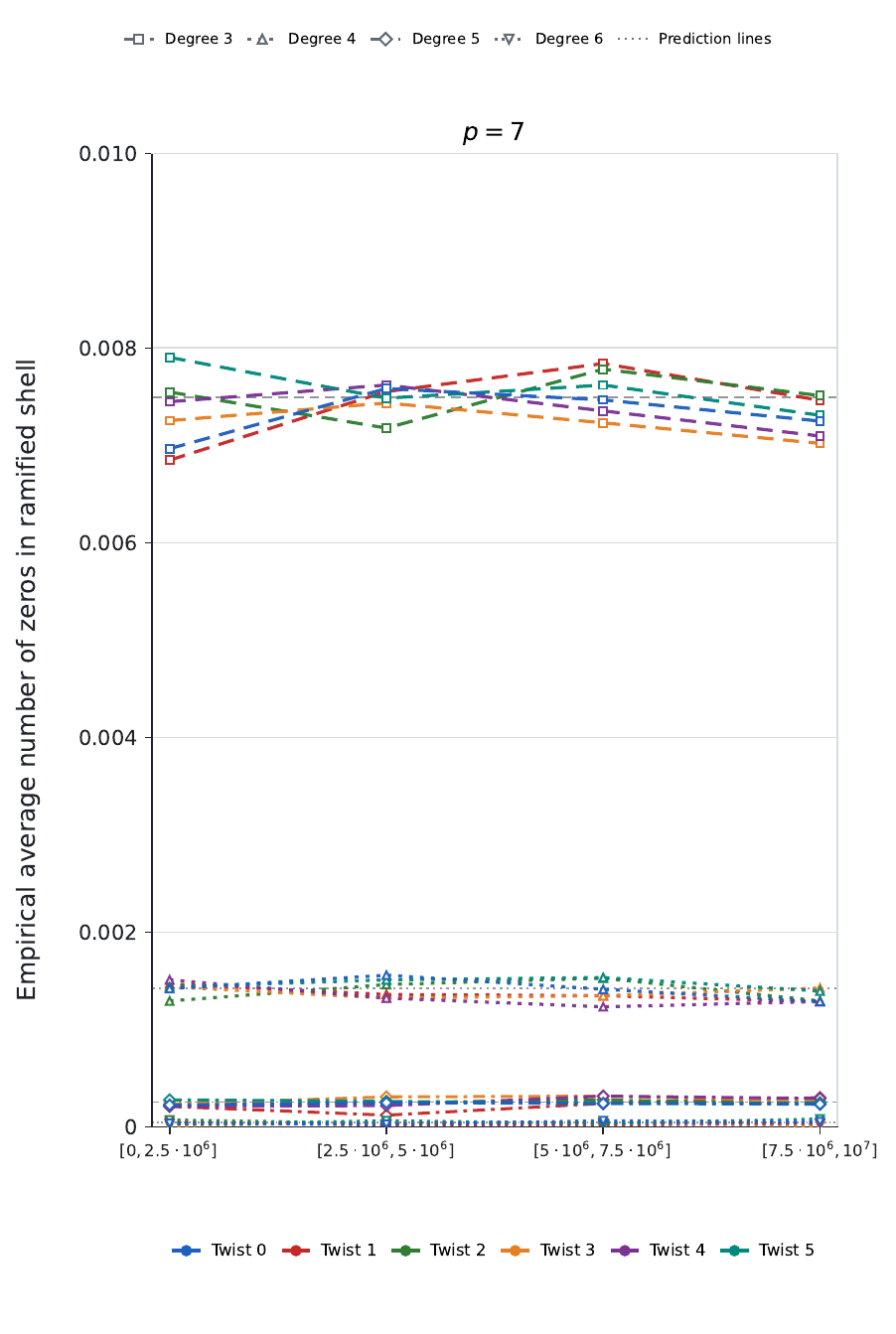}

        \caption{Test of \Cref{cor:ramified_shell_intro}: average number of zeros in $\bigcup_K \pi \mc{O}_K \setminus \pi^2 \mc{O}_K$ where the union is over totally ramified extensions of degree $d=3,4,5,6$, for $p=7$.}
        \label{fig:p7_ramified_shells}
\end{figure}

For quadratic and ramified cubic extensions, we similarly see quite close agreement with the random matrix predictions of \Cref{cor:quadratic_intro} and \Cref{cor:cubic_extensions}. We do not show the case of unramified cubic extensions, even though we have a random matrix prediction in \Cref{cor:cubic_extensions}, because the prediction is so small that the data on zeros is dominated by noise. This can already be seen to a lesser extent for higher-degree extensions at $p=7$ in \Cref{fig:p7_ramified_shells}. Though the prediction for cubic extensions is not exact, the proven upper and lower bounds of \Cref{cor:cubic_extensions} are close enough together to match the same dotted line (at the cost of thickening it a little when $p=3$). For higher-degree extensions, the agreement of the data with the prediction of \Cref{cor:ramified_shell_intro} is shown in Figures~\ref{fig:p3_ramified_shells}, ~\ref{fig:p5_ramified_shells}, and~\ref{fig:p7_ramified_shells} for $p=3,5,7$ respectively.

Last but not least, we must mention the discriminant prediction in \Cref{thm:disc_ejv_intro}, the source of the original doubt in \cite{ellenberg2011modeling}. See \Cref{tab:lambda2-discriminant-low-high-comparison}. We have shown both `low' ($\sim 10^7$) and `high' ($\sim 10^{10}$) conductor ranges to make visible how raising the conductor pushes the values toward the predicted ones, as with the local averages in \Cref{fig:twist1_smallball_high_low}. Similarly to that data, because the computations are much more expensive at high conductor, we only compute twist $1$. 

\begin{table}[htbp]
\centering
\small
\sisetup{group-digits=false}
\setlength{\tabcolsep}{8pt}
\textbf{$p=3$}\par\smallskip
\begin{tabular}{@{}c S[table-format=1.6] S[table-format=1.6] S[table-format=1.6] S[table-format=1.6e-1]@{}}
\toprule
$v_p(\operatorname{Disc} (P_\chi))$ & {\cite{ellenberg2011modeling} data} & {Low conductors} & {High conductors} & {Prediction} \\
\midrule
$1$ & 0.65280 & 0.612304 & 0.594796 & 0.592593 \\
$2$ & 0.22867 & 0.259043 & 0.266798 & 0.271605 \\
$3$ & 0.083467 & 0.083536 & 0.088590 & 0.087791 \\
$4$ & 0.025733 & 0.030525 & 0.031499 & 0.032007 \\
$5$ & 0.008933 & 0.009740 & 0.010956 & 0.010567 \\
$6$ & 0.000267 & 0.003005 & 0.004793 & 0.003624 \\
$7$ & 0.000133 & 0.001267 & 0.001883 & 0.001204 \\
\bottomrule
\end{tabular}
\vspace{5ex}
\par\medskip
\textbf{$p=5$}\par\smallskip
\begin{tabular}{@{}c S[table-format=1.6] S[table-format=1.6] S[table-format=1.6e-1]@{}}
\toprule
$v_p(\operatorname{Disc} (P_\chi))$ & {Low conductors} & {High conductors} & {Prediction} \\
\midrule
$1$ & 0.779747 & 0.760381 & 0.768000 \\
$2$ & 0.175622 & 0.190196 & 0.185600 \\
$3$ & 0.035608 & 0.038485 & 0.036864 \\
$4$ & 0.007573 & 0.008912 & 0.007629 \\
$5$ & 0.001208 & 0.001620 & 0.001524 \\
$6$ & 0.000242 & 0.000203 & 0.000307 \\
$7$ & 0 & 0.000203 & 0.000061 \\
\bottomrule
\end{tabular}
\par\medskip
\vspace{5ex}
\textbf{$p=7$}\par\smallskip
\begin{tabular}{@{}c S[table-format=1.6] S[table-format=1.6] S[table-format=1.6e-1]@{}}
\toprule
$v_p(\operatorname{Disc} (P_\chi))$ & {Low conductors} & {High conductors} & {Prediction} \\
\midrule
$1$ & 0.848910 & 0.842105 & 0.839650 \\
$2$ & 0.133622 & 0.135444 & 0.137443 \\
$3$ & 0.015975 & 0.019507 & 0.019584 \\
$4$ & 0.001194 & 0.002576 & 0.002849 \\
$5$ & 0.000299 & 0.000368 & 0.000407 \\
$6$ & 0 & 0 & 0.000058 \\
$7$ & 0 & 0 & 0.000008 \\
\bottomrule
\end{tabular}
\vspace{5ex}
\caption{Empirical distribution of $v_p(\operatorname{Disc}( P_\chi))$ conditional on $\lambda=2$. Prediction is from \cite{ellenberg2011modeling}, see \Cref{thm:disc_ejv_intro}. `Low conductors' refers to empirical averages over the range $(.75\cdot 10^7,10^7]$, `High conductors' to the range $(10^{10},10^{10}+10^6]$. For $p=3$, data from \cite[p8]{ellenberg2011modeling} at lower conductor range is shown.}
\label{tab:lambda2-discriminant-low-high-comparison}
\end{table}

\textbf{Notes.} When speaking of conductor ranges, we will say `conductor' to refer to the conductor of the underlying character $\theta$, rather than $\chi = \theta \omega^i$. Also, as we discuss in more detail in \Cref{sec:computational_methods}, there is a very small set of characters $\chi$ for which our code does not fully classify the extension type of all of the roots of $F_\chi$ in the unit disk, which leads to upper and lower error bars. In all cases, these error bars have smaller width than the icons plotted in the figure, so we simply plot the icon at the midpoint between the two bounds but do not draw additional error bars. The actual numbers and error bars are contained in the supplement to the arxiv version.

\subsection{What do the positions of the zeros tell us about arithmetic?} \label{subsec:zeros_and_arithmetic_intro} The positions of zeros of $p$-adic $L$-functions also give arithmetic information, specifically about class numbers in cyclotomic towers over an abelian CM extension $K/\Q$. This provides some additional motivation to understand finer statistics of these zeros. For the sake of exposition, let $K$ be quadratic imaginary; the general case is treated as well in \Cref{sec:class_numbers}.

Classical Iwasawa theory is concerned with class groups of fields in the cyclotomic tower $K = K_0 \subset K_1 \subset K_2 \subset \ldots$ over the base field $K$. Since the order of $\text{Gal}(K(\zeta_{p^{n+1}})/K)$ is $dp^n$ where $d \mid (p-1)$, we can define $\Delta \subset \text{Gal}(K(\zeta_{p^{n+1}})/K)$ to be the unique subgroup of order $d$. Then we let $K_n = K(\zeta_{p^{n+1}})^{\Delta}$. 

Let $h_p(K_n)$ be the $p$-part of the class number of $K_n$. Then $h_p(K_n) \mid h_p(K_{n+1})$ for all $n$, and the zeros of an associated $p$-adic $L$-function tell us how the sequence $(h_p(K_n))_{n \geq 0}$ grows. To each quadratic imaginary field $K$ is associated a quadratic Dirichlet character $\chi$, as detailed in \Cref{sec:class_numbers}. Let $F_{\chi \omega}(T)$ be the Iwasawa power series corresponding to the $p$-adic $L$-function $L_p(s,\chi \omega)$, and let $\alpha_1,\ldots,\alpha_\lambda$ be the zeros of $F_{\chi \omega}$ which lie inside the unit disk $B_1(0) \subset \bar{\Q}_p$. Then we have, see \Cref{thm:general-class-number}, that
\begin{equation}
    \label{eq:hquotient_intro}
    \frac{h_p(K_{n+1})}{h_p(K_n)} = u_n \prod_{i=1}^\lambda \prod_\zeta (\alpha_i - (\zeta-1)),
\end{equation}
where $u_n \in \Z_p^\times$ is a unit, and second product is over all primitive $p^{n+1}$-th roots of unity $\zeta$. The left-hand side is always a power of $p$ by definition, and the choice of unit $u_n$ is made to ensure that the right-hand side is as well. The base layer $h_p(K_0)$ is also related to the product of zeros of $F_{\chi \omega}$, see \Cref{thm:base_layer_K_0}.

Generically, each factor in the product in \eqref{eq:hquotient_intro} gives the same contribution, so all that matters is the number of zeros, called the Iwasawa $\lambda$-invariant of $K$. Indeed, 
$$h_p(K_n) = p^{\la n + \nu} \quad \quad \quad \quad \text{ for all sufficiently large $n$,}$$
for another invariant $\nu$ associated to $K$. 

While $\la$ captures the asymptotics, both the secondary term $\nu$ and the fixed-$n$ values of the class numbers $h_p(K_n)$ are controlled by the positions of zeros via \eqref{eq:hquotient_intro}. All these are explained in more detail in \Cref{sec:class_numbers}.

\subsection{Outline} \Cref{sec:plots} gives many pages of additional plots which were referenced above but not included on earlier pages to avoid breaking the flow of text. \Cref{sec:preliminaries} gives the basic definitions of $p$-adic $L$-functions, Iwasawa series, etc. elided over above.

\Cref{sec:class_numbers} gives an exposition of the relation between finer statistics of zeros of $p$-adic $L$-functions and class numbers, expanding upon \Cref{subsec:zeros_and_arithmetic_intro}. This material is all known or almost-known in the literature, but much is scattered among various sources with other aims, and in many cases we must give the proofs here to get exactly the form of the result which seems most enlightening for our purpose. These results are not logically necessary for any others in the paper, but serve as motivation.%

Sections~\ref{sec: Random matrix computations} and~\ref{sec:distributional-consequences} contain the random matrix computations and proofs of Consequences \ref{thm:lambda_intro} through \ref{cor:ramified_shell_intro}, as well as two others (\Cref{cor:Zp_third_moment} and \Cref{cor:cubic_extensions}) which are bounds rather than exact formulas and which we did not state above due to their length and technicality. \Cref{sec: Random matrix computations} contains all of the random matrix computations, while \Cref{sec:distributional-consequences} gathers these to prove the above consequences and some more general versions not stated above.

The logical structure is roughly as follows. The paper \cite{shen2026eigenvalues} proved a formula (\Cref{prop: limiting correlation function} below) for the expected number of eigenvalues in a set, and more generally of pairs and higher $k$-tuples of eigenvalues in the set, as integrals over the set of certain functions analogous to the correlation functions in classical random matrix theory. The correlation functions are not so explicit in general, but were made explicit in \cite{shen2026eigenvalues} for eigenvalues and pairs of eigenvalues in $\Z_p$, and the predictions of \Cref{cor:Zp_zeros_intro} and the variance in \Cref{cor:Zp_higher_moments_intro} come from relatively straightforward computations from this carried out in \Cref{sec: Random matrix computations}. Similarly for quadratic extensions. The totally ramified prediction in \Cref{cor:ramified_shell_intro} also follows smoothly.

For the third moment in $\Z_p$ given in \Cref{cor:Zp_third_moment}, and the case of cubic extensions, no explicit forms or sharp bounds were given in \cite{shen2026eigenvalues}, and the proofs of \Cref{cor:Zp_third_moment} and the cubic version (\Cref{cor:cubic_extensions}) are new to this paper. In the latter case, we must crucially use recent orbital integral results of Cho-Kang-Lee \cite{cho2024stable} as well as \Cref{prop: limiting correlation function}.

\Cref{sec:distributional-consequences} then uses these random matrix results to prove the above consequences and some more general ones. The new content which is not merely collecting results of the previous section is the treatment of the tightness issue for root-counting statistics alluded to above.

\Cref{sec:computational_methods} details our computational methods. The code and raw $p$-adic $L$-function data may also be accessed at \url{https://github.com/knospe/iwasawa-power-series}. Tables compiling the precise numbers plotted in these plots are recorded in a supplement to the arxiv version.

\textbf{Acknowledgments.} We wish to thank Daniel Delbourgo, Jordan Ellenberg, Akshay Venkatesh, and Larry Washington for encouragement, helpful conversations, and comments. We are also grateful to Sonal Jain for exchanges regarding the code and data used in \cite{ellenberg2011modeling}.%

\textbf{Statement on AI use.} The code used to produce the initial data for $p$-adic $L$-functions was written and run by human authors without AI assistance; the version published at the link above was further polished for readability with the aid of AI. The code used to turn this data into plots and tables was entirely written by Codex. The same Codex agent produced the above-mentioned supplementary data table file, with substantial iteration from the authors to improve formatting and readability. To check against errors in these statistics, two authors reproduced all statistics with different AI agents which had access to the same data but not to each others' code, plots or tables. Parts of \Cref{sec: Random matrix computations} and \Cref{sec:distributional-consequences} were written in small chunks by ChatGPT with piece-by-piece prompting by one of the authors and then mostly rewritten; all others are human-written in the first draft. Conversations with AI agents were still helpful at several points in the writing/understanding the material in Sections~\ref{sec:preliminaries} and~\ref{sec:class_numbers}, particularly for the proof of \Cref{thm:correspondence_large_zeros}. AI agents also suggested edits which were then considered and implemented (or not) by the authors. 

\section{Additional plots}\label{sec:plots}

\begin{figure}[H]
        \centering
        \vspace{5ex}
        \includegraphics[scale=.95]{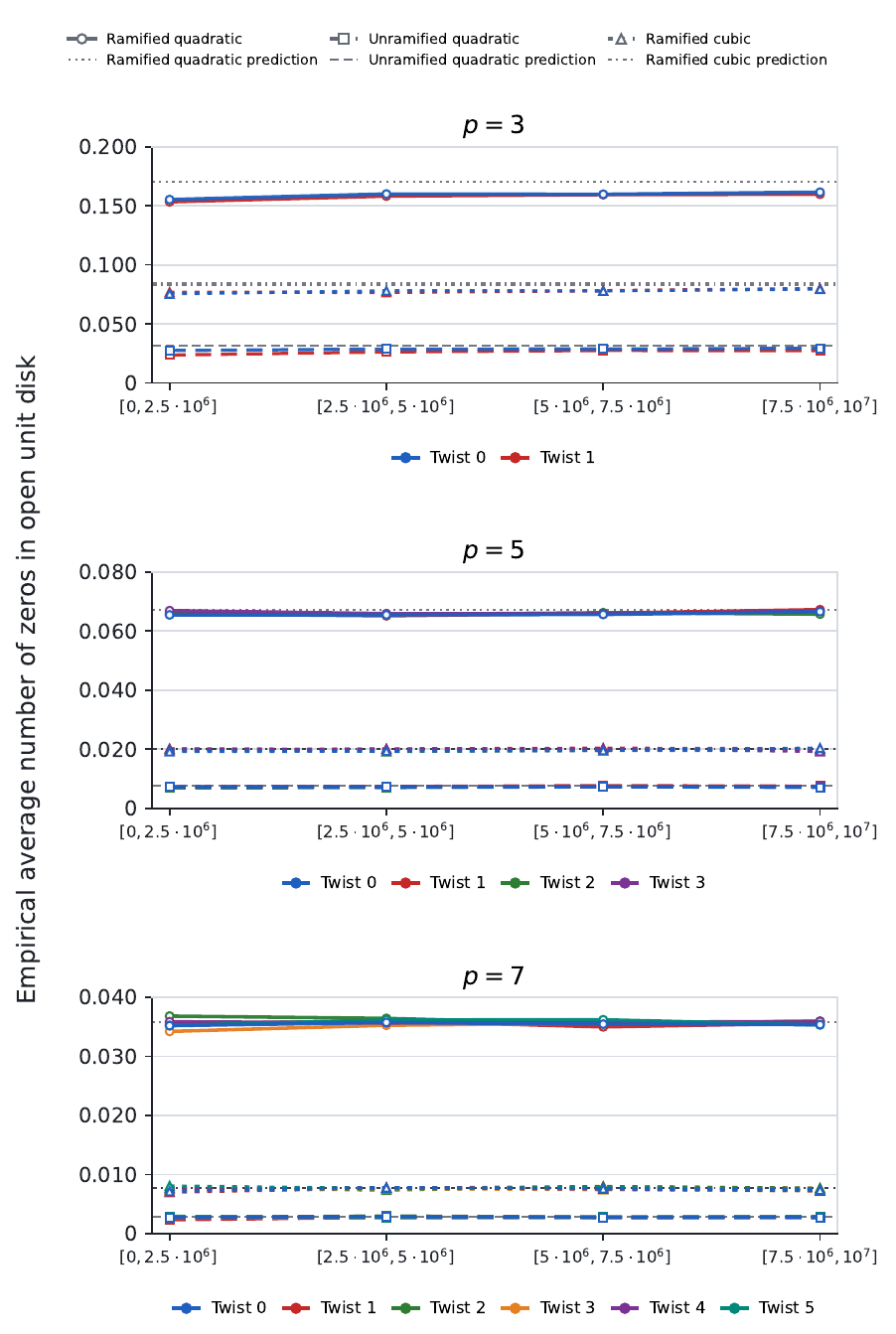}
        \caption{Test of \Cref{cor:quadratic_intro} and \Cref{cor:cubic_extensions}: empirical mean number of zeros in the unit disk which either (1) lie in any ramified quadratic extension (circles), (2) lie in any unramified quadratic extensions (squares), or (3) lie in any ramified cubic extension (triangles). As before, averages are computed over four different conductor ranges given on the $x$-axis, the primes $p=3,5,7$, and all twists $0,\ldots,p-2$. }%
        \label{fig:cubic_quadratic_all_primes}
\end{figure}

\section{Fundamentals from Iwasawa theory}\label{sec:preliminaries}

In this section, $p$ will be an odd prime. We give a more detailed exposition of $p$-adic $L$-functions, their zeros, and relations to arithmetic questions. This is all standard within Iwasawa theory, but experts may still find it worth skimming to glean our notational conventions.

\subsection{$p$-adic $L$-functions associated to Dirichlet characters}

In this article, we consider $p$-adic $L$-functions $L_p(s,\chi)$ associated with a Dirichlet character $\chi$. Recall that a Dirichlet character modulo $d$ is a function $\Z/d\Z \to \bQ$ which is a character (i.e. multiplicative) on $(\Z/d\Z)^\times$ and $0$ otherwise. We call such a character \emph{primitive} if it does not factor through the map $\Z/d\Z \to \Z/m\Z$ for any proper divisor $m$ of $d$. A Dirichlet character $\chi$ is a function on $\Z$ given by precomposing a Dirichlet character modulo $d'$ with the map $\Z \to \Z/d'\Z$, and always possesses a unique positive integer $d \mid d'$ such that $\chi$ factors through a primitive Dirichlet character modulo $d$. We refer to $d$ as the \emph{primitive conductor} of $\chi$.  When multiplying primitive Dirichlet characters, we replace the product with the associated primitive character. 

We fix once and for all embeddings $\bQ \inj \bQ_p$ and $\bQ \inj \C$, so that we can regard the values of $\chi$ (which are roots of unity) as lying either in $\bQ_p$ or $\C$ as appropriate. Here and throughout, $\C_p$ denotes the completion of $\bQ_p$. Viewing the values as complex numbers, we recall that the classical (complex) Dirichlet $L$-function is given by
\begin{equation}
    L(s,\chi) := \sum_{n \geq 1} \frac{\chi(n)}{n^s}
\end{equation}
for $\Re{s} > 1$, and this determines it as a meromorphic function on $\C$. In particular, it takes special values 
\begin{equation}\label{eq:L_bernoulli}
    L(1-k,\chi) = - \frac{B_{k,\chi}}{k}
\end{equation}
at the nonpositive integers $1-k$, where $B_{k,\chi}$ are certain algebraic numbers known as $\chi$-twisted Bernoulli numbers, see e.g. \cite{washingtonbook}. Because they are algebraic, we may view them inside $\bQ_p$, and will do this henceforth.

There are different ways to define the \emph{$p$-adic $L$-function} associated to a Dirichlet character; we take the one which is most relevant to the way we actually compute these functions numerically, see \Cref{sec:special-values}. Specifically, consider the  $L$-values
\begin{equation}
    L_p(1-k,\chi) = (1-\chi(p) p^{k-1}) \cdot L(1-k,\chi), \quad \quad \quad \quad k\geq 1,\ k \equiv 0 \mod (p-1),
    \label{eq:interpol_no_omega}
\end{equation}
where the Euler factor at $p$ is removed. Viewing the $\chi$ in the prefactor as valued in $\bQ_p$, these values have the property that the map $k \mapsto L_p(1-k,\chi)$ is $p$-adically continuous. Furthermore, there exists a unique $p$-adic meromorphic function (analytic if $\chi \neq \bf{1}$) which interpolates the values $L_p(1-k,\chi)$. See e.g. \cite[Theorem 5.11]{washingtonbook}. 

\begin{defi}
    \label{def:p-adic_l-fn}
    We denote by $L_p(s,\chi)$ the unique $p$-adic meromorphic function, defined for $s$ in the open ball in $\C_p$ around $0$ of radius $p^{\frac{p-2}{p-1}}$, for which $L_p(1-k,\chi)$ is given by \eqref{eq:interpol_no_omega} for $k$ as in that equation.
\end{defi}

We have defined $L_p(s,\chi)$ in \eqref{eq:interpol_no_omega} by continuation from one congruence class modulo $p-1$ on the nonpositive integers, but it takes related special values on the other nonpositive integers as well. Recall the (multiplicative) \emph{Teichm\"uller character}
\begin{equation*}
    \omega: (\Z/p\Z)^\times \to \Q_p^\times,
\end{equation*}
where for all $i\in\{1,\ldots,p-1\}$, $\omega(i)$ is the unique element in $i+p\Z_p$ that satisfies $\omega(i)^{p-1}=1$. As with Dirichlet characters, $\omega$ lifts to function $\Z \to \Q_p$ by setting its value on $p\Z$ to $0$, and by abuse of notation we also denote this function by $\omega$. Hence given any Dirichlet character $\chi$, we have a family of `twisted' Dirichlet characters $\chi \omega^r: \Z \to \bQ_p, r=0,1,\ldots,p-2$.

The other special values of the $p$-adic $L$-function are then given in terms of special values of the Dirichlet $L$-functions associated to these twisted characters: for $L_p$ defined as above, we have that
\begin{equation}\label{eq:interpol_general}
    L_p(1-k,\chi) = (1-(\chi\omega^{-k})(p) p^{k-1}) \cdot L(1-k,\chi\omega^{-k})
\end{equation}
for any positive integer $k$, generalizing \eqref{eq:interpol_no_omega}. In other words, the $p$-adic $L$-function intertwines the values of $p-1$ different classical complex Dirichlet $L$-functions associated to the twisted characters $\chi \omega^r$.

For general $\chi$, $L_p(s,\chi)$ may be identically $0$. It is in fact identically zero if and only if $\chi$ is an odd function on $\Z$ (that is, $\chi(-1)=-1$), as the relevant special values of the complex Dirichlet $L$-function used for interpolation vanish in this case. Since $\chi(-1) = \pm 1$ by multiplicativity, all Dirichlet characters are either even or odd, so equivalently $L_p(s,\chi)$ is nontrivial if and only if $\chi$ is even. Henceforth, whenever $F_\chi$ or $P_\chi$ is considered, we assume that $\chi$ is even unless explicitly stated otherwise.

The function $L_p(s,\chi)$ may also have a `trivial zero' at $s=0$. If 
\begin{equation}\label{eq:trivial_zero}
    1-(\chi \omega^{-1})(p) = 0,
\end{equation}
then \eqref{eq:interpol_general} at $k=1$ yields $L_p(0,\chi)=0$. The condition \eqref{eq:trivial_zero} occurs exactly when $\chi = \theta \omega$ and $\theta(p)=1$, which in particular can only occur if $\theta$ has conductor coprime to $p$.

\subsection{The meaning of twists and the family {$X_i(D)$}}

The above discussion was for general $\chi$. In this work, we consider characters of the form $\chi=\theta \omega^i$ where $\theta$ is quadratic. We fix $i$ and consider averages over families $X_i(D)$ of characters of this form, defined above \Cref{conj:distinguished_polynomial_intro}, where $\theta$ ranges over non-trivial characters of conductor at most $D$ satisfying 
\begin{enumerate}
    \item $\theta$ has primitive conductor coprime to $p$;
    \item $\chi$ is even as a function on $\Z$;
    \item if $i=1$, $\theta(p) \neq 1$. 
\end{enumerate}
Now we can explain where these conditions come from. For the second, the $p$-adic $L$-function is trivial if it is not satisfied, as mentioned above. The third condition corresponds to the absence of a trivial zero by \eqref{eq:trivial_zero}. The reason it is only relevant when $i=1$ is that $\omega(p)=0$, so outside of this case, $\chi \omega^{-1}$ automatically vanishes at $p$ and \eqref{eq:trivial_zero} cannot be satisfied.

To explain the condition that $\theta$ must have conductor coprime to $p$, we recall the following standard result:

\begin{lemma}
    \label{thm:dirichlet_character_decomp}
    Let $\chi: \Z \to \bQ_p$ be an arbitrary nontrivial Dirichlet character. Then it has a unique decomposition $\chi = \theta \omega^i \psi$, where
    \begin{itemize}
        \item $\theta$ is a Dirichlet character of conductor $d$ coprime to $p$,
        \item $\omega$ is the Teichm\"uller character where $i \in \{0,\ldots,p-2\}$, and
        \item $\psi$ is either trivial or has conductor $p^{n+1}$ and order $p^n$ for $n \geq 1$. 
    \end{itemize}
\end{lemma}
This is standard and we omit the proof. The characters $\theta \omega^i $ and $\psi$ are called \emph{characters of the first kind} and \emph{characters of the second kind} (with respect to $p$). Our condition on the conductor of $\theta$ simply says that it is of the first kind. As we will see below, multiplying by characters of the second kind gives the same Iwasawa power series, so it is natural to only average over characters of the first kind and their twists.

As mentioned in the introduction, the twists $i=1$ and $i=1+\frac{p-1}{2}$ are special. This is not obvious from \eqref{eq:interpol_general} at this point, since all twists of $\chi$ appear in the interpolation. However, this will become clear in the discussion of the correspondence with the quadratic number fields in \Cref{sec:class_numbers}.

\subsection{The Iwasawa power series}

In a ball around $0$, the $p$-adic $L$-function may also be represented by its \emph{Iwasawa power series}. The following result is due to Iwasawa, see \cite[Theorem 7.10]{washingtonbook}:

\begin{prop}
    \label{thm:iwasawa_series}
    Let $\theta$ be a Dirichlet character of conductor $d\geq 1$ coprime to $p$, $\chi=\theta \omega^i$ where $i\in\{0,1,...,p-2\}$, excluding the case $\theta = \mathbf{1},i=0$ where $\chi$ is trivial. Let $\mc{O}_\chi := \Z_p[\Im(\chi)]$ be the extension of $\Z_p$ defined by the values of $\chi$. Then there exists a unique power series $F_\chi \in \mc{O}_\chi[[T]]$ such that
    \begin{equation}
    {F}_{\chi}\big((1+ dp)^{s} - 1\big) =\;
{L}_p(s, \chi )
\quad\text{for all $s\in\mathbb{C}_p$ with $|s|\! <p^{\frac{p-2}{p-1}}$.}
\label{eq:powerseries}
\end{equation}
\end{prop}

Following \cite{ellenberg2011modeling}, we take the zeros of $F_\chi$ inside the open unit ball in $\C_p$ as the objects of primary interest. Note that the behavior of the $p$-adic exponential map makes the relation between these and the zeros of $L_p(s,\chi)$ somewhat subtle, as we discuss in \Cref{subsec:L_fn_vs_iwasawa_series}. It is however immediate that a trivial zero of $L_p(s,\chi)$ exists (i.e., $L_p(0,\chi)=0$) if and only if $F_\chi(0)=0$. 

\begin{defi}
    \label{def:distinguished_poly}
    Let $\mc{O}_K$ be the ring of integers of a finite extension $K$ of $\Q_p$. Then a monic polynomial $P(T) \in \mc{O}_K[T]$ is \emph{distinguished} if $P(T) \equiv T^{\deg(P)} \pmod{\mf{m}}$, where $\mf{m}$ is the maximal ideal of $\mc{O}_K$.
\end{defi}

\begin{prop}[$p$-adic Weierstrass preparation theorem]
    \label{thm:p-adic_weierstrass}
    Let $\varpi_K$ be a uniformizer of $\mc{O}_K$. With notation as above,
    any nonzero $f(T) \in \mc{O}_K[[T]]$ has a unique factorization
    \begin{equation}
        f(T) = \varpi_K^\mu \cdot u(T) \cdot P(T),
    \end{equation}
    where $\mu \in \Z_{\geq 0}$, $u(T) \in \mc{O}_K[[T]]^\times$ is a unit,
    and $P(T) \in \mc{O}_K[T]$ is a distinguished polynomial.
\end{prop}

By the Ferrero-Washington theorem \cite[Theorem 7.15]{washingtonbook}, the Iwasawa power series $F_{\chi}(T)$ is not divisible by $p$. Hence $F_{\chi}(T)$ can be factorized uniquely as a product of an invertible power series and a distinguished polynomial.

\begin{defi}
    For an even Dirichlet character $\chi$, we denote by $P_\chi(T)$ the distinguished polynomial part of the Iwasawa power series $F_\chi(T)$.
\end{defi}

Any $f(T)$ as in \Cref{thm:p-adic_weierstrass} defines a function on the open unit ball $B_1(0) \subset \C_p$. Since $1/u(T)$ is also a power series and hence also defines a function on $B_1(0)$, $u$ does not have any zeros in $B_1(0)$, hence the zeros of $f$ in $B_1(0)$ are exactly the zeros of the distinguished polynomial $P$. Thus, to study the zeros of the Iwasawa power series $F_\chi$ which lie in the unit ball $B_1(0)$, it suffices to study the zeros of the distinguished polynomial $P_\chi$. 

\begin{rmk}
    \label{rmk:relation_to_ejv_normalization}
    The normalization that we use in \Cref{thm:iwasawa_series} is not the same as the one in \cite{ellenberg2011modeling}, and results in a different Iwasawa power series, but it is the same as theirs up to harmless change of variables, and none of the statistics considered here depend on this choice. 
\end{rmk}

\subsection{Zeros of $p$-adic $L$-functions and Iwasawa power series}\label{subsec:L_fn_vs_iwasawa_series}

Since we follow \cite{ellenberg2011modeling} and study zeros of the Iwasawa power series rather than those of the $p$-adic $L$-function itself, it is worth mentioning how they are related. There is a close relationship, but it becomes necessary to consider $p$-adic $L$-functions of Dirichlet characters twisted by {\em characters of the second kind} (see \Cref{thm:dirichlet_character_decomp}) which have the same Iwasawa power series.

Note that 
$$
    \zeta_{\psi}:= \psi(1+dp)^{-1}
$$
is either $1$ or a root of unity of $p$-power order. Writing $\chi = \theta \omega^i$ as before, we then have
\begin{equation}
    {F}_{\chi}\big( \zeta_{\psi} (1+ dp)^{s} - 1\big) =\;
{L}_p(s, \chi \psi)
\quad\text{for all $s\in\mathbb{C}_p$ with $|s|\! <p^{\frac{p-2}{p-1}}$}
\label{eq:powerseries_general_psi}
\end{equation}
for any character $\psi$ of the second kind. This generalizes \Cref{thm:iwasawa_series} and is also proven in \cite[Theorem 7.10]{washingtonbook}. 

\begin{rmk}
    Note that the same Iwasawa power series $F_\chi$ governs the values of the whole family of $p$-adic $L$-functions ${L}_p(s, \chi \psi)$ as $\psi$ varies.
\end{rmk}

The precise correspondence between zeros of $F_\chi$ and zeros of the $p$-adic $L$-functions is captured in the following two results.

\begin{prop}\label{thm:Lfn_iwasawa_zeros_correspondence} Let $\chi$ be a Dirichlet character of the first kind. Then
there is a one-to-one correspondence between the zeros of $F_{\chi}(T)$ satisfying $\val(T)>\frac{1}{p-1}$, and the zeroes of $L_p(s,\chi)$ on its domain of definition $B_r(0)$ where $r=p^{\frac{p-2}{p-1}}$.
\end{prop}

Hence each zero of $L_p(s,\chi)$ corresponds to a zero of $F_\chi(T)$, but $F_\chi(T)$ may have additional zeros. The treatment of potential zeroes of $F_{\chi}(T)$ with $\val(T)\leq \frac{1}{p-1}$ is more subtle: they may correspond to a zero of $L_p(s,\chi\psi)$, where $\psi$ is a non-trivial character of the second kind, or they may not correspond to the zero of any such $p$-adic $L$-function. The following gives a criterion for when each case occurs in terms of the location of the zero.

\begin{prop}\label{thm:correspondence_large_zeros} Let $\chi$ be an even Dirichlet character of the first kind and let $T_0$ be a zero of $F_{\chi}(T)$ satisfying $\val(T_0)\leq \frac{1}{p-1}$. Then the following are equivalent:
\begin{enumerate}
\item There exists $n \in \Z_{\geq 1}$ such that $\val(T_0) = \frac{1}{(p-1)p^{n-1}}$, and a primitive $p^n$-th root of unity $\zeta$ such that
\begin{equation}
    \label{eq:T_0_close_zeta}
    \val( T_0 - (\zeta -1) ) > \frac{1}{p-1}.
\end{equation} \label{item: T_0 and zeta}
\item There exists a character of the second kind $\psi$ for which $s_0$ defined by
\begin{equation}
    \label{eq:def_s_0}
    s_0 = \frac{\log(1+T_0)}{\log(1+dp)}
\end{equation}
satisfies $\val(s_0) > \frac{1}{p-1} -1$ and
is a zero of $L_p(s,\chi \psi)$. \label{item: existence of psi}
\end{enumerate}
Furthermore, if these conditions are satisfied, then there is a unique such $\zeta$ satisfying \eqref{eq:T_0_close_zeta}, and the character $\psi$ associated to it by \Cref{thm:exists_psi} satisfies the second condition.

\end{prop}

\begin{rmk}
    We explicitly do not rule out the possibility that there are multiple characters $\psi$, not just the one associated to $\zeta$, for which $L_p(s_0,\chi \psi)=0$. It seems likely simply by measure considerations that this never occurs, but we are not sure how to establish it.
\end{rmk}

The rest of this subsection is concerned with proving Propositions \ref{thm:Lfn_iwasawa_zeros_correspondence} and \ref{thm:correspondence_large_zeros}. We note that \Cref{thm:correspondence_large_zeros} is a slight extension of \cite[Proposition 3]{ervanmet1992} and the proof is inspired by the one there. Since we have taken $B_{p^{(p-2)/(p-1)}}(0)$ as the domain of definition of $L_p(s,\chi)$, \Cref{thm:Lfn_iwasawa_zeros_correspondence} follows immediately from \Cref{thm:analytic_bijection} below.

\begin{prop}\label{thm:analytic_bijection} Let $d$ be coprime to $p$. Then, the $p$-adic function $f$ given by
$$ f(s) = (1+dp)^s - 1$$
is an analytic bijection between the open balls $B_r(0)$ and $B_{r'}(0)$ where 
$r=p^{\frac{p-2}{p-1}}$ and $r'= p^{-1/(p-1)}$.%
\end{prop}
\begin{proof} (see \cite[Proposition 1]{ervanmet1992})
    Recall that $\val(\exp(x)-1)=\val(x)$ when $\val(x) > \tfrac{1}{p-1}$, by Taylor expansion. Hence 
    $$\val(f(s)) = \val(\exp(s \log(1+dp))-1) = \val(s \log(1+dp)).$$
    By Taylor expansion, using that $d$ is coprime to $p$, $\val(\log(1+dp)) = 1$. Hence 
    $$\val(f(s)) = \val(s) + 1,$$
    so $f$ is indeed an analytic map from $B_r(0)$ to $B_{r'}(0)$. 
    
    The inverse map is given by 
    $$f^{-1}(T) = \frac{\log(1+T)}{\log(1+dp)},$$
    which is defined for all $T \in B_{r'}(0)$. For such $T$, Taylor expansion similarly gives $\val(\log(1+T)) = \val(T)$, hence 
    $$\val(f^{-1}(T)) = \val(T)-1.$$
    This shows that $f^{-1}(T) \in B_r(0)$ when $T \in B_{r'}(0)$. Hence $f$ is a bijection.

\end{proof}

\begin{lemma}
    \label{thm:exists_psi}
    Let $d$ be coprime to $p$ and $n \geq 1$. Then the map
    $$
        \psi \mapsto \psi(1+dp)^{-1}
    $$
    defines a bijection between the set of characters of the second kind of order $p^n$ and the set of primitive $p^n$-th roots of unity.
\end{lemma}
\begin{proof}
    A character $\psi$ as in the statement is a lift to $\Z$ of a multiplicative character $(\Z/p^{n+1}\Z)^\times \to \bQ^\times$ of order $p^n$. It is a standard fact easily verifiable by induction on $r$ that for units $u \in \Z_p^\times$, $(1+up)^{p^r} = 1+u'p^{r+1}$ for another unit $u'$. Hence
    $$
        (1+dp)^{p^{n}} \equiv 1 \pmod{p^{n+1}}
    $$
    but
    $$
        (1+dp)^{p^{n-1}} \not \equiv 1 \pmod{p^{n+1}},
    $$
    so $1+dp$ has order $p^{n}$ and generates the unique $p$-Sylow subgroup of $(\Z/p^{n+1}\Z)^\times$. Thus every character of $(\Z/p^{n+1}\Z)^\times$ of order $p^n$ sends $1+dp$ to a primitive $p^n$-th root of unity, and conversely any such assignment defines a character of order $p^n$, so the statement holds.
\end{proof}

\begin{proof}[Proof of {\Cref{thm:correspondence_large_zeros}}] 

\eqref{item: T_0 and zeta} $\Rightarrow$ \eqref{item: existence of psi}: Let $\psi$ be the unique character of the second kind of order $p^n$ such that $\psi(1+dp)^{-1} = \zeta$ guaranteed by \Cref{thm:exists_psi}. Since $F_\chi(T_0)=0$, it suffices to show for 
\begin{equation*}
    s_0 = \frac{\log(1+T_0)}{\log(1+dp)}
\end{equation*}
as defined in \eqref{eq:def_s_0} that
\begin{equation*}
    L_p(s_0,\chi \psi) = F_\chi(T_0),
\end{equation*}
which we will do using the generalized relation between the $p$-adic $L$-function and Iwasawa power series given in \eqref{eq:powerseries_general_psi}. For \eqref{eq:powerseries_general_psi} to apply, we must check that 
\begin{equation}
    \label{eq:1to2_suffices_1}
    T_0 = \zeta(1+dp)^{s_0}-1 
\end{equation}
and 
\begin{equation}\label{eq:s0_small}
    \val(s_0) > \frac{1}{p-1}-1,
\end{equation}
the latter of which ensures that $s_0$ lies in the domain of $L_p(s,\chi \psi)$.

To do this, we first change variables to $x_0 = \zeta^{-1}(1+T_0) - 1$. Our $x_0$ serves as a proxy for $T_0$, and the desired equality \eqref{eq:1to2_suffices_1} is trivially equivalent to 
\begin{equation}
    \label{eq:1to2_suffices_2}
    1+x_0 = (1+dp)^{s_0}.
\end{equation}
The advantage of \eqref{eq:1to2_suffices_2} over \eqref{eq:1to2_suffices_1} is that it no longer features the prefactor $\zeta$, allowing us to use facts about bijectivity of the logarithm function near $1$.

Using that $\val(\zeta)=0$ and our hypothesis \eqref{eq:T_0_close_zeta}, we have
$$ \val (x_0)= \val(\zeta x_0) = \val(T_0 - (\zeta-1)) > \frac{1}{p-1}.$$
Series expansion of the logarithm then yields
\begin{equation}\label{eq:vplogx0}
    \val(\log(1+x_0)) > \frac{1}{p-1}.
\end{equation}
Furthermore,
\begin{equation}
    \label{eq:x0_and_T0}
    \log(1+T_0)=\log(1+x_0)
\end{equation}
since $1+x_0 = \zeta^{-1}(1+T_0)$ and $\log(\zeta)=0$. Starting with formula \eqref{eq:def_s_0} for $s_0$, the denominator has valuation $\val(\log(1+dp)) = 1$, while for the numerator we convert $x_0$ to $T_0$ by \eqref{eq:x0_and_T0} and then use the bound \eqref{eq:vplogx0}. Combining these yields the desired bound \eqref{eq:s0_small}.

We now have
\begin{align}
    \begin{split}
        (1+dp)^{s_0} &= \exp(s_0 \log(1+dp)) \\ 
        &= \exp(\log(1+T_0)) \\ 
        &= \exp(\log(1+x_0)) \\ 
        &= 1+x_0,
    \end{split}
\end{align}
where the first equality comes from the bound \eqref{eq:s0_small}, the second is the definition \eqref{eq:def_s_0} of $s_0$, the third is \eqref{eq:x0_and_T0}, and the last follows by the bound \eqref{eq:vplogx0} similarly to the first. This proves \eqref{eq:1to2_suffices_2} and hence \eqref{eq:1to2_suffices_1}, completing the proof of \eqref{item: T_0 and zeta} $\Rightarrow$ \eqref{item: existence of psi}.

\eqref{item: existence of psi} $\Rightarrow$ \eqref{item: T_0 and zeta}: 
Suppose that
\begin{equation}\label{eq:def_s_0_2}
    s_0 = \frac{\log(1+T_0)}{\log(1+dp)}
\end{equation}
is a zero of $L_p(s,\chi \psi)$ for some character $\psi$ of the second kind, where $s_0$ lies in the domain of definition of $L_p(s,\chi \psi)$, i.e.,
\begin{equation}
    \val(s_0) > \frac{1}{p-1}-1.
\end{equation}
Equivalently by \eqref{eq:def_s_0_2}, $\val(\log(1+T_0)) > \frac{1}{p-1}$. 
Then $\exp(\log(1+T_0)) = 1+y_0$, where $y_0 \in \C_p$ satisfies
\begin{equation}
    \label{eq:vpy0}
    \val(y_0) > \frac{1}{p-1}.
\end{equation}
Furthermore, $\log(1+T_0)=\log(1+y_0)$. Because the kernel of $\log$ on $B_1(1)$ is exactly the set of $p$-power roots of unity, there is a $p^n$-th root of unity $\zeta$ for some $n \geq 0$ such that 
\begin{equation} 1+T_0 = \zeta ( 1 + y_0). 
\label{eq:1plust0}
\end{equation}

Now, \eqref{eq:1plust0} and \eqref{eq:vpy0} yield
\begin{equation*}
    \val( T_0 - (\zeta -1)) > \frac{1}{p-1},
\label{eq:valT0}
\end{equation*}
which is the desired \eqref{eq:T_0_close_zeta}. 

Since $\val(T_0) \leq \frac{1}{p-1}$ and $\val(y_0) > \frac{1}{p-1}$, $T_0$ and $y_0$ cannot be equal, so the root of unity in \eqref{eq:1plust0} must be non-trivial, i.e., $n\geq 1$. Hence we have 
\begin{equation*}
    \val(\zeta - 1) = \frac{1}{(p-1)p^{n-1}} \leq \frac{1}{p-1}.
\end{equation*}
Inequality \eqref{eq:T_0_close_zeta} and the strong triangle inequality then imply
 $$ \val(T_0) = \val ( \zeta - 1) = \frac{1}{(p-1)p^{n-1}} .$$
This completes the proof of \eqref{item: existence of psi} $\Rightarrow$ \eqref{item: T_0 and zeta}.

For the additional `Furthermore' of the statement, assume for the sake of contradiction that $\zeta'$ is another $p$-power root of unity satisfying (\ref{eq:T_0_close_zeta}). Then 
\begin{equation*}
    \val(\zeta-\zeta') = \val((T_0-(\zeta'-1)) - (T_0 - (\zeta-1))) > \frac{1}{p-1}
\end{equation*}
by \eqref{eq:T_0_close_zeta} and the strong triangle inequality, but also
\begin{equation*}
    \val(\zeta - \zeta') = \val(1- \zeta' \zeta^{-1} ) \leq \frac{1}{p-1}
\end{equation*}
because $\zeta'\zeta^{-1}$ is a $p$-power root of unity which is not $1$. This is a contradiction, hence $\zeta$ is unique. The proof of \eqref{item: T_0 and zeta} $\Rightarrow$ \eqref{item: existence of psi} gives that the character $\psi$ for which $\psi(1+dp)^{-1} = \zeta$ satisfies \eqref{item: existence of psi}.

\end{proof}

\section{Zeros of $p$-adic $L$-functions and class numbers}\label{sec:class_numbers}

This section explains the arithmetic significance of the zeros of $p$-adic
$L$-functions and Iwasawa power series introduced above, which is not logically necessary for the rest of the paper, but serves as motivation. It also explains why the twists $i=1$ and $i=1+\frac{p-1}{2}$, for which we observe unusually few zeros near $0$ in the associated Iwasawa power series when the conductor is low, are exactly those associated with quadratic imaginary fields.

\subsection{The general case} For general abelian CM fields, the $p$-adic $L$-functions encode the $p$-primary part of class numbers
and their growth along cyclotomic $\Z_p$-extensions. We will make this
relationship explicit in terms of the distinguished polynomials $P_\chi$. We first need some elementary facts about Dirichlet characters, which may be found in \cite[Chapter 3]{washingtonbook}.

Let $n$ be a positive integer. Class field theory gives a canonical isomorphism between $(\Z/n\Z)^{\times}$ and the Galois group $G=\text{Gal}(\Q(\zeta_n)/\Q)$.
A Dirichlet character of primitive conductor $n$ can then be viewed canonically as a character of $G$, since it comes from lifting a Dirichlet character $\chi : (\Z/n\Z)^{\times} \rightarrow \bQ$ modulo $n$ to $\Z$. We will make heavy use of this identification without comment in this section.

\begin{defi}
    Let $\chi: \Gal(\Q(\zeta_n)/\Q) \to \bQ^\times$ be a group character and $H=\chi^{-1}(1) \subset \Gal(\Q(\zeta_n)/\Q)$ be its kernel. We refer to the fixed field $K$ of $H$ as the {\em field associated with $\chi$} or {\em the field cut out $\chi$}. 
    
    Similarly, for a subgroup $X$ of Dirichlet characters modulo $n$, we define the associated field of $X$ to be the fixed field of the kernel $H=\{g \in \Gal(\Q(\zeta_n)/\Q): \chi(g)=1 \text{ for all }\chi \in X\}$. 
\end{defi}

Conversely, if $K/\Q$ is an abelian field extension, then the Kronecker-Weber theorem implies that $K$ is contained in a cyclotomic field $\Q(\zeta_n)$ for some $n\in \Z_{\geq 1}$. The smallest such $n$ defines the {\em conductor} of $K$. Now the characters of $\text{Gal}(K/\Q)$ can be viewed as Dirichlet characters modulo $n$ which vanish on $\text{Gal}(\Q(\zeta_n)/K)$. These characters form a subgroup $X$ of the group all Dirichlet characters modulo $n$, and $X$ is called the group of Dirichlet characters associated to $K$. 

Note that by viewing $\Gal(K/\Q)$ as a quotient $\Gal(\Q(\zeta_n)/\Q)/H$, $X$ is identified with the character group of $\Gal(K/\Q)$. It follows immediately that $\#X = [K:\Q]$.

\begin{example}\label{ex:quadratic_characters_and_fields}
    By the above, quadratic fields correspond to (nontrivial) quadratic Dirichlet characters, as these are the only characters generating groups of order $2$. Quadratic characters which are even (resp. odd) correspond to fields which are totally real (resp. pure imaginary). 

    To see this, first suppose $\chi$ is even. Then $\chi(-1)=1$, so the element $\zeta_n \mapsto \zeta_n^{-1}$ of $\Gal(\Q(\zeta_n)/\Q)$ (which is just complex conjugation) lies in the kernel of $\chi$. Consequently, the corresponding fixed field is totally real. On the other hand, if $\chi$ is odd, then the kernel does not include complex conjugation, so the fixed field is not totally real, hence (because it is quadratic) it is pure imaginary. $\square$
    
\end{example}

Let $K$ be a CM field, i.e., an imaginary quadratic extension of a totally real field, with class number $h(K)$. Let $K^+$ be the maximal real subfield of $K$ and  associated class number $h^+(K)$. By base changing the Hilbert class field of $K^+$ to $K$, we have $h^+(K) \mid h(K)$, and we then define
\begin{equation}
    \label{eq:hminus}
    h^-(K) = \frac{h(K)}{h^+(K)}. 
\end{equation}

In the following, we consider {\em abelian} CM fields $K$. Note that an abelian field is either totally real or CM, since the complex conjugation is central in its Galois group. We assume that the conductor $\mc{D}$ of $K$ is either $\mc{D}=d$ or $\mc{D}=dp$, where $d$ is coprime to $p$. This assumption on the conductor implies that all Dirichlet characters associated to $K$ are of the first kind.

\begin{prop} \label{thm:base_layer_K_0} Let $K$ be an abelian CM field of conductor not divisible by $p^2$. Then $h^-$ defined in \eqref{eq:hminus} satisfies

\begin{equation}\label{eq:class_number_formula} h^-(K) \sim w \prod_{\chi \text{ odd}} B_{1,\chi} \sim \prod_{\substack{\chi \text{ odd} \\ \chi \omega \neq \bf{1}}} B_{1,\chi},
\end{equation}
where $\sim$ denotes equality of $p$-adic valuations, the product runs over the odd Dirichlet characters $\chi$ associated to $K$, and
\begin{equation}
    w = \begin{cases}
        p & \text{if }K \text{ contains a $p$-th root of unity} \\ 
        1 & \text{if }K \text{ does not contain a $p$-th root of unity}
    \end{cases}
\end{equation}
is the number of $p$-power roots of unity contained in $K$.

If in addition, $\chi(p) \neq 1$ (i.e., $L_p(s,\chi\omega)$ has no trivial zero)  for all odd characters $\chi$ associated to $K$, then 
\begin{equation}\label{eq:after_class_number}
    h^-(K) \sim  \prod_{\substack{\chi \text{ odd} \\ \chi \omega \neq \bf{1}} } (1-\chi(p))^{-1} L_p(0,\chi \omega) \sim  \prod_{\substack{\chi \text{ odd} \\ \chi \omega \neq \bf{1}}} (1-\chi(p))^{-1} P_{\chi \omega}(0).
\end{equation}
Furthermore, if $p$ does not divide the order of $\chi$ or if $p$ divides the conductor of $\chi$, then the factor $(1-\chi(p))^{-1}$ is a $p$-adic unit.  
\end{prop}
\begin{proof} 
Note that the value of $w$ can only be $1$ or $p$ because our assumption on the conductor excludes $K$ which contain $p^2$-th roots of unity. Then the first equality of \eqref{eq:class_number_formula} is the analytic class number for $h^-(K)$ (up to a $p$-adic unit), see \cite[Theorem 4.17]{washingtonbook}. If $K$ does not contain a $p$-th root of unity then both products in \eqref{eq:class_number_formula} are identical. Otherwise, we have $v_p(B_{1,\omega^{-1}}) = -1$ which gives the right-hand side of \eqref{eq:class_number_formula}.

For the first equality of \eqref{eq:after_class_number}, recall that by the $k=1$ case of \eqref{eq:L_bernoulli} we have $L(0,\chi) = - B_{1,\chi}$. For $k=1$, the interpolation factor of $L_p(0,\chi \omega)$ in \eqref{eq:interpol_general} is $1-\chi(p)$, which is non-zero by assumption. If
$\chi(p)$ is a $d$-th root of unity for some $d$ coprime to $p$ or if 
$\chi(p)=0$, then
$1-\chi(p)$ has valuation $0$ in $\bQ_p$. Hence $B_{1,\chi} = -(1-\chi(p))^{-1} L_p(0,\chi \omega)$, which gives the first equality.

For the second equality of \eqref{eq:after_class_number}, \Cref{thm:iwasawa_series} ensures that  $L_p(0,\chi \omega)$ is replaced by $F_{\chi \omega}(0)$. Taking the Weierstrass factorization (\Cref{thm:p-adic_weierstrass}) of $F_{\chi \omega}$, the unit term does not change the valuation, and $\mu=0$ by the Ferrero-Washington theorem as mentioned, so $\val(F_{\chi \omega}(0))=\val(P_{\chi \omega}(0))$ and this gives the result.
\end{proof}

Hence, {\em the number and the valuations} of zeros of the Iwasawa power series of the corresponding characters directly determine the minus class number, i.e., more zeros close to the origin yield larger $p$-class groups. 

A similar formula gives the exact growth of the $p$-class numbers in the cyclotomic $\Z_p$ extension $K_{\infty}$ of $K$. Let $K$ be an abelian field of conductor not divisible by $p^2$ as above, and let $K_n$ be the $n$-th layer of the extension, which we recall is defined as follows. The order of $\text{Gal}(K(\zeta_{p^{n+1}})/K)$ is always $dp^n$ where $d \mid (p-1)$, and we let $\Delta \subset \text{Gal}(K(\zeta_{p^{n+1}})/K)$ be the unique subgroup of order $d$. Then we let $K_n = K(\zeta_{p^{n+1}})^{\Delta}$.

$\text{Gal}(K_n/K)$ can be viewed as the multiplicative group $(1+p \Z_p)/(1+p^{n+1}\Z_p)$ (canonically), which is isomorphic to $\Z/p^n\Z$ (noncanonically). In this way, the characters of $\text{Gal}(K_n/K)$ become Dirichlet characters of the second kind of conductor $p^{n+1}$ and order $p^n$.
The full Galois group $\text{Gal}(K_n/\Q)$ decomposes into $\text{Gal}(K/\Q) \times \text{Gal}(K_n/K)$ because $K \cap \Q_n = \Q$. The characters of $\text{Gal}(K_n/\Q)$ have a corresponding decomposition into a product of a character $\chi$ of the first kind and a character $\psi$ of the second kind.

\begin{theorem}\label{thm:growth-general} Let $K$ be an abelian CM field and $K_n$ be as above. Let $n\geq 1$. Then
\begin{equation}
\label{eq:growth-general}
    \frac{h^-(K_n)}{h^-(K)} \sim
    \prod_{\substack{\chi \text{ odd} \\ \chi\omega \neq \bf{1}}} \ \prod_{\psi \neq \bf{1}} L_p ( 0, \chi \omega \psi) \sim
\prod_{\substack{\chi \text{ odd} \\ \chi\omega \neq \bf{1}}} \ \prod_{\psi \neq \bf{1}} P_{\chi \omega}(\zeta_{\psi} - 1),
\end{equation} 
where the products are over the nontrivial odd characters $\chi$ of $K$ and the non-trivial characters $\psi$ of the second kind of order dividing $p^n$. 
\label{thm:general-class-number}
\end{theorem}
Note that the above formula \eqref{eq:growth-general} also holds in the trivial zero case, i.e., when $P_{\chi \omega}(0)=0$.
\begin{proof} 

We adapt \cite[Theorem 7.13 and Theorem 7.14]{washingtonbook} (including errata), where the case $K=\Q(\zeta_{dp})$ and
$K_n=\Q(\zeta_{dp^{n+1}})$ is treated. 

The analytic class number formula for $h^-(K_n)$ gives
$$ h^-(K_n) \sim w_n \prod_{\chi \text{ odd,   } \psi} B_{1,\chi \psi},$$
where $w_n$ is the number of roots of unity in $K_n$. The product runs over all odd characters $\chi$ of $K_n$ of the first kind and all characters of the second kind $\psi$ of order dividing $p^n$, by the discussion above the theorem statement. Here, we use that characters of the second kind are always even so $\chi \psi$ is odd if and only if $\chi$ is odd. 

This product can then be decomposed into two factors: 
$$ \left( w_0 \prod_{\chi \text{ odd}} B_{1,\chi }  \right)  \cdot\left(\frac{w_n}{w_0} \prod_{\substack{\chi \text{ odd} \\ \psi \neq \bf{1}}}  B_{1,\chi \psi} \right). $$
The first factor equals $h^-(K)$ (up to a $p$-adic unit). In the following, we investigate the second factor, which can be rewritten as
\begin{equation}
\label{eq:secondfactor}
\frac{w_n}{w_0}
\left(
\prod_{\psi \neq \bf{1}} B_{1,\omega^{-1} \psi} \right)
\left(
\prod_{\substack{\chi \text{ odd, } \chi\omega \neq \bf{1} \\ \psi \neq \bf{1}}}  B_{1,\chi \psi} \right) 
\end{equation}
The first product of \eqref{eq:secondfactor} is over characters $\psi \neq \bf{1}$ of the second kind, and is only present if $\chi = \omega^{-1}$ is a character of $K$, i.e., if $K$ contains a non-trivial $p$-th root of unity $\zeta_p$.

\noindent {\bf Case 1:} $\zeta_p \notin K$.

Then $K_n$ does not contain non-trivial $p$-th roots of unity, so that $w_0 \sim 1$ and $w_n \sim 1$. Since $\chi \omega \neq \bf{1}$, we have
\begin{equation}
    B_{1,\chi \psi} = - L_p(0, \chi\omega \psi) \sim P_{\chi \omega}(\zeta_{\psi} -1 ),
    \label{eq:bernoulli-general}
\end{equation}
where in the first equality we note that $\psi(p)=0$ so the Euler factor disappears, and the second equality is by \eqref{eq:powerseries_general_psi} and \Cref{thm:p-adic_weierstrass}. This implies the assertion.

\noindent {\bf Case 2:} $\zeta_p \in K$.

Then the cyclotomic extension field $K_n$ contains all $p^{n+1}$-th roots of unity, and so
$w_0 = p$ and $w_n = p^{n+1}$. We need to consider the first product in equation \eqref{eq:secondfactor} which can be written as
$$ \prod_{\psi \neq \bf{1}} L_p(0,\psi), $$
since the interpolation factor $(1-\psi \omega^{-1}(p))$ is trivial.
Here the corresponding character of the first kind is trivial, and in this situation, the $p$-adic $L$-function $L_p(s,\psi)$ is meromorphic rather than holomorphic. However, it can still be expressed as a {\em quotient} of power series, see \cite[Theorem 7.10]{washingtonbook} and the text before that theorem:
$$ L_p(s,\psi) = \frac{g(\zeta_{\psi}(1+p)^s-1, \bf{1})}{h(\zeta_{\psi}(1+p)^s-1, \bf{1})}. $$
Therefore, we have
$$ L_p(0,\psi) = \frac{g(\zeta_{\psi} -1, \bf{1})}{h(\zeta_{\psi} -1, \bf{1})} .$$
It is shown in \cite[Lemma 7.12]{washingtonbook} and the proof of  \cite[Theorem 7.13]{washingtonbook} that $g(\zeta_{\psi} -1, {\bf 1})$ is a unit and 
$h(\zeta_{\psi}(1+p)^s-1, {\bf 1}) \equiv 1-\zeta_{\psi}^{-1} \mod p$. From this, we obtain $\val(L_p(0,\psi)) = - \frac{1}{(p-1)p^{n-1}}$ for characters $\psi$ of order $p^n$. The sum of valuations over all non-trivial characters $\psi$ of the second kind of order dividing $p^n$ gives $-n$. Since $\val(\frac{w_n}{w_0}) = n$, we obtain 
\begin{equation*}
\frac{w_n}{w_0}
\left(\prod_{\psi \neq \bf{1}} B_{1,\omega^{-1} \psi} \right) \sim
\frac{w_n}{w_0}\left(\prod_{\psi \neq \bf{1}} L_p(0,\psi)  \right)
\sim 1.
\end{equation*}
For the second product in equation \eqref{eq:secondfactor}, we can again use equation \eqref{eq:bernoulli-general}. This completes the proof.
\end{proof}
In particular, this means that the positions of the zeros of $P_{\chi \omega}$ control the growth of class numbers $h^-(K_n)$, recovering the classical $\lambda$-invariant growth formula when $n$ is large.

\begin{cor}
    In the setting of \Cref{thm:growth-general},
    \begin{equation}
        \val(h^-(K_{n+1}))-\val(h^-(K_n))
        = \sum_{\substack{\chi \text{ odd} \\ \chi \omega \neq \bf{1}}} \sum_\alpha \sum_\zeta \val(\alpha - (\zeta-1))
    \end{equation}
    where $\chi$ ranges over nontrivial odd characters of $K$, $\alpha$ ranges over all roots of $P_{\chi \omega}$, and $\zeta$ ranges over all primitive $p^{n+1}$-th roots of unity. 
\end{cor}
\begin{proof}
    This is immediate from \Cref{thm:growth-general}.
\end{proof}

In particular, let $\la_p(\chi \omega) = \deg(P_{\chi \omega})$ be  the $\lambda$-invariant associated to this character, and $\la$ be the sum of these. If
    \begin{equation}
        \val(\alpha) > \frac{1}{(p-1)p^n}
    \end{equation}
    for all $\alpha$, which will be true whenever $n$ is large enough, then
    \begin{equation}
        \val \left(\frac{h^-(K_{n+1})}{h^-(K_n)}\right) = \sum_{\substack{\chi \text{ odd} \\ \chi \omega \neq \bf{1}}} \la_p(\chi \omega). 
    \end{equation}
    This recovers the classical Iwasawa growth formula \cite[Theorem 11]{iwasawa1959gamma}, which in the present abelian setting takes the form
    $$
    \val\bigl(h^-(K_n)\bigr)=\lambda n+\nu
    $$
    for all sufficiently large $n$. The $\mu$-invariant vanishes by 
    the theorem of Ferrero-Washington \cite{ferrerowashington}. The $\nu$-invariant may also be written in terms of contributions from the zeros and roots of unity $\zeta$ via \Cref{thm:base_layer_K_0} and \Cref{thm:growth-general}.

\subsection{The quadratic imaginary case} Let $K$ be an imaginary quadratic field and let $\chi$ be the associated Dirichlet character, which is odd by \Cref{ex:quadratic_characters_and_fields}. Depending on whether the conductor of $K$ is divisible by $p$ or not, we either have $\chi=\theta$ or $\chi = \theta \omega^{\frac{p-1}{2}}$, where $\theta$ is a quadratic character with conductor coprime to $p$. In the first case, $\theta$ must be odd, and in the latter case the product of $\theta$ and $\omega^{\frac{p-1}{2}}$ is odd. By specializing \Cref{thm:growth-general} to this case,
\begin{equation} \label{eq:quadratic_example_hquotient}
    \frac{h^-(K_{n+1})}{h^-(K_n)} \sim \frac{h(K_{n+1})}{h(K_n)} \sim \prod_\zeta P_{\chi \omega}(\zeta-1),
\end{equation}
where the product is over primitive $p^{n+1}$-th roots of unity $\zeta$. The first equality follows since $K_n^+=\Q_n$, where $\Q_n$ denotes
the $n$th layer of the cyclotomic $\Z_p$-extension of $\Q$, whose $p$-class group is trivial by \cite{iwasawa1956note}; the second
follows from \Cref{thm:growth-general}.

The characters appearing here are exactly the ones with twists $1$ or $1+\frac{p-1}{2}$ in our terminology. As seen earlier, these are the ones for which the zeros of $P_{\chi \omega}$ are anomalously sparse close to $0$ when the conductor is low.

If in addition $\chi(p) \neq 1$ (i.e., no trivial zero) then at the base layer we also have by \Cref{thm:base_layer_K_0} that
\begin{equation}
    h^-(K) = h(K) \sim P_{\chi \omega}(0).
\end{equation}
Here we note that because $p$ is odd and $\chi$ has order $2$, $p$ does not divide the order of $\chi$, so the factor $(1-\chi(p))^{-1}$ is a unit. Consequently, in this case
\begin{equation}
    h^-(K_n) \sim h(K_n) \sim \prod_\zeta P_{\chi \omega}(\zeta-1),
\end{equation}
where the product is over all (not necessarily primitive) $p^n$-th roots of unity. 

Denote by $\alpha_1,\ldots,\alpha_\lambda$ the zeros of $P_{\chi \omega}$, the above is
\begin{equation}
    \label{eq:hquotient_not_intro}
    \frac{h_p(K_{n+1})}{h_p(K_n)} \sim \prod_{i=1}^\lambda \prod_\zeta (\alpha_i - (\zeta-1)),
\end{equation}
given earlier as \eqref{eq:hquotient_intro}. Let us first see how this formula recovers the fact that the number of zeros in the unit disk, which is the $\lambda$-invariant studied, controls the asymptotic growth of $h_p(K_n)$. Recalling the elementary fact that $v_p(\zeta-1) = \frac{1}{(p-1)p^n}$, we see that at least for large $n$, one typically has
\begin{equation}
    \label{eq:val_typical_alpha_zeta}
    v_p(\alpha_i) > \frac{1}{(p-1)p^n}
\end{equation}
and hence 
\begin{equation}
    v_p(\alpha_i-(\zeta-1)) = v_p(\zeta-1) = \frac{1}{(p-1)p^n}.
\end{equation}
If at a given $n$, \eqref{eq:val_typical_alpha_zeta} holds true for  all roots $\alpha_i$ of the Iwasawa series in the unit disk, then 
\begin{equation}\label{eq:no_weird_zeros_recurrence}
    \frac{h_p(K_{n+1})}{h_p(K_n)} = p^{\lambda},
\end{equation}
where we recall that $\lambda$ is the number of roots of $F_{\chi \omega}$ in the unit disk. For any given $\chi$, for sufficiently large $n$ \eqref{eq:val_typical_alpha_zeta} will hold for all $i$. Hence \eqref{eq:no_weird_zeros_recurrence} implies that 
\begin{equation}\label{eq:lambda_and_nu_intro}
    h_p(K_n) = p^{\lambda n + \nu}
\end{equation}
for all sufficiently large $n$, for some $\nu$.%

While the $\lambda$-invariant gives the main term in this asymptotic growth, finer information on the positions of zeros controls the class groups in the lower layers of the cyclotomic tower, and hence the $\nu$-invariant. More specifically, if \eqref{eq:val_typical_alpha_zeta} does not hold for a given $\alpha_i$, then there are two possibilities:
\begin{enumerate}
    \item If $v_p(\alpha_i) < \frac{1}{(p-1)p^n}$, then $v_p(\alpha_i-(\zeta-1)) = v_p(\alpha_i) < \frac{1}{(p-1)p^n}$.
    \item If $v_p(\alpha_i) = \frac{1}{(p-1)p^n}$, then it may happen that $v_p(\alpha_i-(\zeta-1)) > \frac{1}{(p-1)p^n}$ due to cancellations.
\end{enumerate}
In other words, $|\alpha_i|$ may be unusually large and cause $|\alpha_i-(\zeta-1)|$ to be unusually large, or $\alpha_i$ may be unusually close to a (shifted) $p$-power root of unity $\zeta-1$ and cause $|\alpha_i-(\zeta-1)|$ to be unusually small. The value of the $\nu$-invariant in \eqref{eq:lambda_and_nu_intro} is a sum of contributions from both of these events and the base layer $h_p(K_0)$.

These are details of the positions of zeros within the disk $B_1(0)$ which the $\lambda$-invariant does not capture, so a full understanding of the local behavior of zeros would yield an understanding of the joint distribution of class numbers in the whole cyclotomic tower, not only their asymptotic growth.%

\section{Random matrix computations}\label{sec: Random matrix computations}

Our \Cref{conj:distinguished_polynomial_intro} predicts that the distinguished polynomials associated with $p$-adic $L$-functions converge in distribution to those arising from large random matrices over $\Z_p$. To obtain concrete predictions from this conjecture, we need 
\begin{enumerate}
\item a precise description of the limiting random polynomial, and
\item formulas for the expected numbers of eigenvalues in prescribed subsets of finite extensions of $\Q_p$.
\end{enumerate}
These are the two purposes of this section.

In \Cref{subsec:limiting distinguished polynomial}, we recall the existence and almost sure squarefreeness of the limiting distinguished polynomial $P_{A_\infty}$, using a result of the second author \cite[Lemma 4.4]{shen2026resultant}. We then relate this limit to the original multiplicative Haar model of Ellenberg-Jain-Venkatesh \cite[Theorem 4.1]{ellenberg2011modeling}. This gives a common limiting polynomial for the models appearing in the conjecture and in the earlier heuristic. The rest of the section concerns local eigenvalue statistics. Our starting point is the limiting correlation functions established in \cite{shen2026eigenvalues}, including the comparison between additive Haar measure on $\Mat_n(\Z_p)$ and multiplicative Haar measure on $\GL_n(\Z_p)$. We recall the results needed here and use them to compute the eigenvalue counts and correlations relevant to our numerical tests.

Unless otherwise stated, matrices in $\Mat_n(\Z_p)$ (resp. $\GL_n(\Z_p)$) are distributed according to additive (resp. multiplicative) Haar probability measure. Throughout this section we study the limit as $n\to\infty$.

\subsection{Existence of the limiting distinguished polynomial}\label{subsec:limiting distinguished polynomial}

We first recall the existence of the limiting random polynomial $P_{A_\infty}$ appearing in \Cref{conj:distinguished_polynomial_intro}. Recall from \Cref{sec: Introduction} that $\mathscr P=\bigsqcup_{d\geq0}\mathscr P_d$ denotes the space of distinguished polynomials, where
$$
\mathscr P_d=\{x^d+a_{d-1}x^{d-1}+\cdots+a_0:a_0,\ldots,a_{d-1}\in p\Z_p\},
$$
with $\mathscr P_0=\{1\}$. For $A\in\Mat_n(\Z_p)$, let $P_A\in\mathscr P$ be the distinguished factor of $\det(xI_n-A)$, whose roots are precisely the eigenvalues of $A$ lying in the open unit disk, counted with multiplicity. The following result shows that these random polynomials admit a limiting distribution as $n\to\infty$.

\begin{thm}[{\cite[Lemma 4.4]{shen2026resultant}}]\label{thm:existence of limiting distinguished polynomial}
For every $n\ge 1$, let $A_n$ be a random element of $\Mat_n(\Z_p)$ distributed according to additive Haar probability measure, and let $P_{A_n}$ be its distinguished polynomial, a $\mathscr P$-valued random variable. Then there exists $\mathscr{P}$-valued random variable $P_{A_\infty}$ such that $P_{A_n}$ converges weakly to $P_{A_\infty}$ as $n\to\infty$. Moreover, $P_{A_\infty}$ is almost surely squarefree.
\end{thm}

\begin{proof}
This is the special case $S=\{1\}$ of \cite[Lemma 4.4]{shen2026resultant}. In the notation of that paper, $F_1(x)=x$, so the corresponding lifted subspace is the open unit disk. Consequently, $\mathscr P^{\{1\}}=\mathscr P$, and the polynomial $P_{A_n}^{\{1\}}$ defined there is precisely $P_{A_n}$ in our notation. The claimed weak convergence and almost sure squarefreeness follow directly.
\end{proof}

We next consider the original random matrix model of Ellenberg-Jain-Venkatesh \cite{ellenberg2011modeling}, in which the matrices are sampled from $\GL_n(\Z_p)$ according to multiplicative Haar measure\footnote{Ellenberg-Jain-Venkatesh also considered a random matrix model based on $\GSp_\alpha(2n,\Z_p)$, but showed that its distinguished polynomials converge weakly to the same limiting distribution as those of the $\GL_n(\Z_p)$ model, provided that $\alpha\not\equiv 1\pmod p$; see \cite[Theorem 4.1]{ellenberg2011modeling}.}. Their model concerns eigenvalues close to $1$, rather than $0$. The following theorem shows that, after shifting the eigenvalues by $1$, the distinguished polynomials converge to the same limiting random polynomial $P_{A_\infty}$.

\begin{thm}[cf. {\cite[Theorem 4.1]{ellenberg2011modeling}}]\label{thm:existence of limiting distinguished polynomial GL}
For every $n\ge 1$, let $G_n\in\GL_n(\Z_p)$ be a random matrix distributed according to multiplicative Haar probability measure, and let $P_{G_n-I_n}$ be the distinguished factor of its characteristic polynomial $\det(xI_n-(G_n-I_n))$. Then, as $n\to\infty$, we have
$$
P_{G_n-I_n} \to P_{A_\infty} \quad \quad \quad \quad \text{in distribution,}
$$
 where $P_{A_\infty}$ is the limiting random polynomial defined in \Cref{thm:existence of limiting distinguished polynomial}.
\end{thm}

\begin{proof}
Fix $d\ge 0$ and $m\ge 1$, and let $S\subset\mathscr P_d$ be the preimage of a subset of the monic degree-$d$ polynomials in $(\Z/p^m\Z)[x]$ under coefficientwise reduction modulo $p^m$. By \cite[Theorem 4.1]{ellenberg2011modeling}, we have
\begin{equation}
    \label{eq:Gn_to_Gd}
    \lim_{n\to\infty}\Pr(P_{G_n-I_n}\in S)
=\left(\prod_{j=1}^\infty(1-p^{-j})\right)\Pr(P_{G_d-I_d}\in S).
\end{equation}
Since $S\subset\mathscr P_d$, for any $A \in \Mat_d(\Z_p)$, if $P_{A}\in S$ then $A$ is nilpotent modulo $p$, and consequently $I_d+A \in \GL_d(\Z_p)$. Recall also that if a random matrix in $\Mat_d(\Z_p)$ with additive Haar distribution is conditioned to lie in $\GL_d(\Z_p)$, then it has multiplicative Haar distribution. Applying this to a random matrix $A_d$ distributed by the additive Haar measure on $\Mat_d(\Z_p)$, we therefore have
\begin{align}
\begin{split}\label{eq:Gd_to_Ad}
    \Pr(P_{G_d-I_d}\in S) &= \Pr(P_{A_d} \in S | I_d+A_d \in \GL_d(\Z_p)) \\ 
    &= \frac{\Pr(P_{A_d} \in S \text{ and } I_d+A_d \in \GL_d(\Z_p))}{\Pr(I_d + A_d \in \GL_d(\Z_p))} \\ 
    &=\frac{\Pr(P_{A_d}\in S)}{\prod_{j=1}^d(1-p^{-j})}.
    \end{split}
\end{align}

The last equality comes since $I_d + A_d \in \GL_d(\Z_p)$ if and only if $I_d + A_d \pmod{p} \in \GL_d(\F_p)$, and $I_d+A_d \pmod{p}$ is a uniformly random element of $\Mat_d(\F_p)$, so the probability is a simple explicit computation.

We next compute $\Pr(P_{A_n}\in S)$ for $n\ge d$, following the counting argument of \cite[Section 4.3]{ellenberg2011modeling}. Write $R=\Z/p^m\Z$ and let $B$ be an endomorphism of $R^n$. For sufficiently large $N$, the kernels and images of successive powers of $B$ stabilize, since they are increasing (respectively, decreasing) sequences of submodules of a finite module $R^n$. Setting $X=\ker(B^N)$ and $Y=\operatorname{im}(B^N)$, we obtain a unique $B$-invariant decomposition $R^n=X\oplus Y$, where $B|_X$ is nilpotent and $B|_Y$ is invertible. Both $X$ and $Y$ are free $R$-modules. Moreover, the distinguished factor of the characteristic polynomial of $B$ modulo $p^m$ is precisely the characteristic polynomial of $B|_X$.

It follows that $P_{A_n}\in S$ if and only if, for the reduction $B$ of $A_n$ modulo $p^m$, the corresponding $X$ has rank $d$ and the characteristic polynomial of $B|_X$ belongs to the reduction of $S$ modulo $p^m$. By the orbit-stabilizer theorem, there are
$$
\frac{|\GL_n(R)|}{|\GL_d(R)|\,|\GL_{n-d}(R)|}
$$
ordered decompositions $R^n=X\oplus Y$ with $\operatorname{rank}(X)=d$. For each decomposition, there are $|\Mat_d(R)|\Pr(P_{A_d}\in S)$ choices for $B|_X$, since the required condition already forces $B|_X$ to be nilpotent modulo $p$, and hence nilpotent over $R$. Independently, there are $|\GL_{n-d}(R)|$ choices for $B|_Y$. Since the decomposition is uniquely determined by $B$, dividing the resulting count by $|\Mat_n(R)|$ gives
\begin{align}
\label{eq:An_to_Ad}
\begin{split}
\Pr(P_{A_n}\in S)
&=\frac{|\GL_n(R)|}{|\GL_d(R)|\,|\GL_{n-d}(R)|}
\frac{|\Mat_d(R)|\,|\GL_{n-d}(R)|}{|\Mat_n(R)|}
\Pr(P_{A_d}\in S)\\
&=\frac{\prod_{j=1}^n(1-p^{-j})}{\prod_{j=1}^d(1-p^{-j})}
\Pr(P_{A_d}\in S),
\end{split}
\end{align}
where we used that $|\GL_k(R)|/|\Mat_k(R)|=\prod_{j=1}^k(1-p^{-j})$.

Combining \eqref{eq:Gn_to_Gd} and \eqref{eq:Gd_to_Ad} yields
$$
\lim_{n\to\infty}\Pr(P_{G_n-I_n}\in S)=\frac{\prod_{j=1}^\infty(1-p^{-j})}{\prod_{j=1}^d(1-p^{-j})}\Pr(P_{A_d}\in S).
$$
Applying \eqref{eq:An_to_Ad} to the right-hand side, we then have
$$\lim_{n\to\infty}\Pr(P_{G_n-I_n}\in S) = \frac{\prod_{j=1}^\infty(1-p^{-j})}{\prod_{j=1}^n(1-p^{-j})}\Pr(P_{A_n}\in S)$$
where $n\ge d$ on the right-hand side is arbitrary. Taking the limit $n \to \infty$ on the right-hand side thus yields
\begin{align*}
    \begin{split}
        \lim_{n\to\infty}\Pr(P_{G_n-I_n}\in S) &= \lim_{n \to \infty}\frac{\prod_{j=1}^\infty(1-p^{-j})}{\prod_{j=1}^n(1-p^{-j})}\Pr(P_{A_n}\in S) \\ 
        &= \Pr(P_{A_\infty} \in S)
    \end{split}
\end{align*}
where the last equality follows from \Cref{thm:existence of limiting distinguished polynomial}. Since this holds for every $d\ge 0$, $m\ge 1$, and every subset $S\subset\mathscr P_d$ determined by reduction modulo $p^m$, the claimed weak convergence follows.
\end{proof}

The preceding results identify the limiting random distinguished polynomial for both the additive and multiplicative Haar models. This might at first seem surprising, but recall that we are considering only the distinguished part of the characteristic polynomial, and hence already throwing away most of the full characteristic polynomial. The multiplicative model has the additional restriction that the characteristic polynomial cannot have roots in the open unit disk around $-1$, but this disk only contains $O(1)$ eigenvalues\footnote{Meant here in a probabilistic sense, i.e., the number of eigenvalues in this disk converges to an integer-valued random variable with no rescaling.} in the additive model, while the matrix of course has $n$ eigenvalues where $n \to \infty$. So it is not unreasonable that this difference between the models does not really affect the asymptotics in the unit disk around $0$.

We now turn to the eigenvalues themselves. Here, explicit expectations for the number of eigenvalues in a set, or for the number of tuples of eigenvalues in several sets, are given in terms of certain \emph{correlation functions} defined below. These are the quantities needed to test the local predictions of \Cref{conj:distinguished_polynomial_intro}.

\subsection{Correlation functions and expected eigenvalue counts}

We fix once and for all an algebraic closure $\bar{\Q}_p/\Q_p$, and when we speak of an extension $K/\Q_p$ we implicitly also fix an embedding into this algebraic closure. For a finite extension $K/\Q_p$ with ring of integers $\mc{O}_K$, we write $\mc{O}_K^{\new}:=\{x\in\mc{O}_K:\Q_p(x)=K\}$ for the elements of $\mc{O}_K$ that do not lie in a proper subextension of $K/\Q_p$. Note that $\mc{O}_K^{\new}$ is a measurable subset with full measure inside $\mc{O}_K$.

\begin{defi}\label{defi: tuples of eigenvalues}
Let $K_1,\ldots,K_m$ be extensions of $\Q_p$ and $r := \sum_{i=1}^m [K_i:\Q_p]$. Let $U \subset \mc{O}_{K_1}^{\new} \times \cdots \times \mc{O}_{K_m}^{\new}$ be any measurable set. For all $n\ge r$, we denote by
$$Z_{U,n}=Z_{U,n}(A),\quad A\in\Mat_n(\Z_p)$$
the number of $m$-tuples $(x_1,\ldots,x_m)$ of distinct, pairwise non-conjugate eigenvalues of $A$ that lie in $U$. This gives a random variable over $\Mat_n(\Z_p)$.
\end{defi}

For example, given $A\in\Mat_n(\Z_p)$, $Z_{\Z_p,n}=Z_{\Z_p,n}(A)$ counts the number of eigenvalues inside $\Z_p$, and $Z_{\Z_p\times\Z_p,n}$ counts the number of pairs of distinct eigenvalues inside $\Z_p$.

The correlation functions are the densities of the expected counting measures associated with these random variables. More precisely, for fixed $K_1,\ldots,K_m$ and $n$, the finite-size correlation function $\rho_{K_1,\ldots,K_m}^{(n)}$ is characterized by
$$
\E[Z_{U,n}]=\int_U\rho_{K_1,\ldots,K_m}^{(n)}(x_1,\ldots,x_m)dx_1\cdots dx_m
$$
for every measurable $U\subset\mc{O}_{K_1}^{\new}\times\cdots\times\mc{O}_{K_m}^{\new}$, where each $dx_i$ is the additive Haar probability measure on $\mc{O}_{K_i}$. Here $\rho_{K_1,\ldots,K_m}^{(n)}$ is not a probability density (it does not in general integrate to $1$), and should instead be viewed as a local intensity of a point process. For $m=1$ its integral over $U$ is the expected number of eigenvalues in $U$, while for $m=2$ it gives the expected number of ordered pairs of distinct, non-conjugate eigenvalues in $U$. These functions $\rho_{K_1,\ldots,K_m}^{(n)}$ are the $p$-adic analogues of the correlation functions used in classical random matrix theory, see \cite[Section 1]{shen2026eigenvalues} for more explanation of this.

\begin{defi}\label{def:delta_sigma}
    Let $x_1,\ldots,x_m \in \bar{\Q}_p$, and for each $i$ let $r_i$ be the number of Galois conjugates of $x_i$ (including itself) in $\bar{\Q}_p$ and let $\{\sigma_{i,l_i}(x_i):1\le l_i\le r_i\}$ denote the full set of Galois conjugates. Then we write
    $$\Delta_\sigma(x_1,\ldots,x_m)=\Delta(\sigma_{1,1}(x_1),\ldots,\sigma_{1,r_1}(x_1),\ldots,\sigma_{m,1}(x_m),\ldots,\sigma_{m,r_m}(x_m))$$
    where $\Delta(z_1,\ldots,z_r)=\prod_{1\le i<j\le r}(z_i-z_j)$ is the usual Vandermonde determinant. 
\end{defi}

Note that the absolute value of the usual discriminant is related to $\Delta_\sigma$ by $$|\Disc(x_1,\ldots,x_m)|=|\Delta_\sigma(x_1,\ldots,x_m)|^2.$$

\begin{defi}\label{defi: Den and distance}
Let $K/\Q_p$ be an extension of degree $r$, $\mc{O}_K$ its ring of integers, and $x \in \mc{O}_K^{\new}$. Then we define
$$\Mod_{\Z_p[x]}=\{M\subset K\mid M \text{ is a $\Z_p[x]$-module, also a $\Z_p$-lattice of rank $r$}\},$$
equivalently the set of fractional ideals of the ring $\Z_p[x]$.

We have a canonical action of the group $K^\times$ on $\Mod_{\Z_p[x]}$ by multiplication, and similarly for its subgroup $\Lambda = \{\pi^n: n \in \Z\}$ where $\pi$ is a uniformizer of $K$. This latter action has finitely many orbits, and we denote the number by $\#(\Lambda \backslash \Mod_{\Z_p[x]})$. 

Finally, we define a function $\Den: \mc{O}_K^{\new} \to \R_{\ge 0}$, which depends implicitly on $K$, by
$$\Den(x)=\frac{1}{1-p^{-r/e}}\cdot |\Delta_\sigma(x)|\cdot \#(\L\backslash \Mod_{\Z_p[x]})$$
where $e$ is the ramification index of $K/\Q_p$.
\end{defi}

For example, when $x\in\Z_p$, we have $\#(\L\backslash \Mod_{\Z_p[x]})=|\Delta_\sigma(x)|=1$, $r=e=1$, and therefore $\Den(x)=\frac{1}{1-p^{-1}}$.

For the calculations below, we need the large-matrix limits of these local intensities. The second and third author  \cite[Theorems 6.1 and 6.2]{shen2026eigenvalues} show that, for every fixed tuple of extensions $K_1,\ldots,K_m$, the finite-size correlation functions converge as $n\to\infty$ to limiting correlation functions $\rho_{K_1,\ldots,K_m}^{(\infty)}$.

\begin{thm}[{\cite[Theorem 6.1 and Theorem 6.2]{shen2026eigenvalues}}]\label{prop: limiting correlation function}
Let $K_1,\ldots,K_m$ be extensions of $\Q_p$ and $r := \sum_{i=1}^m [K_i:\Q_p]$. Let $U \subset \mc{O}_{K_1}^{\new} \times \cdots \times \mc{O}_{K_m}^{\new}$ be any measurable set. Then, we have
$$\lim_{n\rightarrow\infty}\E[Z_{U,n}]=\int_U\rho^{(\infty)}_{K_1,\ldots,K_m}(x_1,\ldots,x_m)dx_1\cdots dx_m,$$
where
\begin{equation}
\rho_{K_1,\ldots,K_m}^{(\infty)}(x_1,\ldots,x_m)= |\Delta_\sigma(x_1,\ldots,x_m)| \cdot \left(\prod_{i=1}^m\Den(x_i) \right) \cdot \lim_{n\rightarrow\infty}\E|\det(Z(A))|.
\end{equation}
Here, the integration (and definition of measurability of $U$) is with respect to the Haar probability measure on $\mc{O}_{K_1}\times \cdots \times \mc{O}_{K_m}$, $Z(x)=\prod_{i=1}^m Z_i(x)$ for $Z_i\in\Z_p[x]$ the minimal polynomial of $x_i$, and $A\in\Mat_{n-r}(\Z_p)$ is distributed by the additive Haar measure.
\end{thm}

The comparison in \Cref{thm:existence of limiting distinguished polynomial GL} concerns the limiting distinguished polynomial itself. To compare the expected counts of eigenvalues, we also need to compare the correlation functions of the two matrix models. The following result shows that the limiting correlation functions for multiplicative Haar measure on $\GL_n(\Z_p)$ agree with those for additive Haar measure on $\Mat_n(\Z_p)$ whenever all the eigenvalues are $p$-adic units.

\begin{thm}[{\cite[Theorem 7.5]{shen2026eigenvalues}}]\label{thm: limiting correlation function GL}
Let $K_1,\ldots,K_m$ be finite extensions of $\Q_p$ and $r:=\sum_{i=1}^m[K_i:\Q_p]$. Let $U\subset\mc{O}_{K_1}^{\new}\times\cdots\times\mc{O}_{K_m}^{\new}$ be any measurable set. For $n\ge r$, let $Z_{U,n}^{\GL}$ denote the number of $m$-tuples of distinct, pairwise non-conjugate eigenvalues in $U$ of a multiplicatively Haar-distributed matrix in $\GL_n(\Z_p)$. Then we have
$$
\lim_{n\rightarrow\infty}\E[Z_{U,n}^{\GL}]
=\int_U\rho_{K_1,\ldots,K_m}^{(\infty),\GL}(x_1,\ldots,x_m)dx_1\cdots dx_m,
$$
where
$$
\rho_{K_1,\ldots,K_m}^{(\infty),\GL}(x_1,\ldots,x_m)
=
\begin{cases}
\rho_{K_1,\ldots,K_m}^{(\infty)}(x_1,\ldots,x_m),
& x_i\in\mc{O}_{K_i}^{\times}\text{ for all }i,\\
0,&\text{otherwise}.
\end{cases}
$$
\end{thm}

In particular, this comparison applies to the eigenvalues close to $1$ considered by Ellenberg-Jain-Venkatesh. Indeed, if $|x_i|<1$ for every $i$, then $1+x_i\in\mc{O}_{K_i}^{\times}$. Combining \Cref{prop: limiting correlation function}, \Cref{thm: limiting correlation function GL} with translation invariance of additive Haar measure on $\Mat_n(\Z_p)$ therefore gives
$$
\rho_{K_1,\ldots,K_m}^{(\infty),\GL}(1+x_1,\ldots,1+x_m)
=\rho_{K_1,\ldots,K_m}^{(\infty)}(x_1,\ldots,x_m).
$$
Therefore, the additive Haar model used in the calculations below yields the same limiting local correlation functions as the shifted multiplicative Haar model in the open unit disk.

In general, we do not expect the expression for the limiting correlation functions in \Cref{prop: limiting correlation function} to admit a substantial further simplification. There is, however, an important special case in which these contributions simplify remarkably.

\begin{thm}[{\cite[Theorem 6.6]{shen2026eigenvalues}}]\label{thm: correlation for generator}
Let $K$ be a finite extension over $\Q_p$, and $x_1\in\mathcal{O}_K$ such that $\mathcal{O}_K=\Z_p[x_1]$. Then, we have
$$\rho_{K}^{(\infty)}(x_1)=|\Disc_{K/\Q_p}|.$$
\end{thm}

We now apply these correlation functions to the expected eigenvalue counts used in our numerical comparisons.

\subsection{Eigenvalues in $\Z_p$}

\begin{thm}[{\cite[Theorem 9.1]{shen2026eigenvalues}}]\label{thm:1-pt_cor_fn_Zp}
For all positive integers $n$ and $x_1\in\Z_p$, we have $\rho_{\Q_p}^{(n)}(x_1)=1$.
\end{thm}

In light of \Cref{thm:1-pt_cor_fn_Zp}, integrating $\rho_{\Q_p}^{(n)}$ over a ball such as $p^l \Z_p$ immediately yields the following:

\begin{cor}\label{thm: 1pt expectation in plZp}
We have $\E[Z_{p^l\Z_p,n}]=p^{-l}$ for all $n\ge 1$ and $l\ge 0$.
\end{cor}

\begin{thm}[{\cite[Theorem 9.3]{shen2026eigenvalues}}]\label{thm: 2pt in Zp}
Suppose $x_1,x_2\in\Z_p$, such that $|x_1-x_2|=p^{-m}$. Then, we have
$$\rho_{\Q_p\times\Q_p}^{(\infty)}(x_1,x_2)=p^{-m}\sum_{k\ge 0}(-1)^kp^{-\frac{(2m+1)k^2+(4m+1)k}{2}}(1+p^{-k-1}).$$
\end{thm}

\begin{thm}\label{thm: 2pt expectation in plZp}
Let $l\ge 0$. Then, we have 
$$\lim_{n\rightarrow\infty}\E[Z_{p^l\Z_p\times p^l\Z_p,n}]=(1-p^{-1})\sum_{k\ge0}(-1)^k(1+p^{-k-1})
\frac{p^{-\frac{k^2+k}{2}}\;p^{-l(k^2+2k+3)}}{1-p^{-(k^2+2k+2)}}.$$
\end{thm}

\begin{proof}
Notice that for all $m\ge l$, the subset
$$\mathcal{U}_{m,l}:=\{(x_1,x_2)\in p^l\Z_p\times p^l\Z_p:|x_1-x_2| = p^{-m}\}\subset p^l\Z_p\times p^l\Z_p$$
has measure $\mu_{\Q_p\times\Q_p}(\mathcal{U}_{m,l})=(1-p^{-1})p^{-l-m}$. Therefore, following \Cref{prop: limiting correlation function} and the correlation function given in \Cref{thm: 2pt in Zp}, we have
\begin{align}
\begin{split}
\lim_{n\rightarrow\infty}\E[Z_{p^l\Z_p\times p^l\Z_p,n}]&=\int_{p^l\Z_p\times p^l\Z_p}\rho_{\Q_p\times\Q_p}^{(\infty)}(x_1,x_2)dx_1dx_2\\
&=\sum_{m\ge l}\mu_{\Q_p\times\Q_p}(\mathcal{U}_{m,l})\rho_{\Q_p\times\Q_p}^{(\infty)}(0,p^m)\\
&=\sum_{m\ge l}(1-p^{-1})p^{-l-2m}\sum_{k\ge 0}(-1)^kp^{-\frac{(2m+1)k^2+(4m+1)k}{2}}(1+p^{-k-1})\\
&=(1-p^{-1})\sum_{k\ge 0}(-1)^k(1+p^{-k-1})p^{-\frac{k^2+k}{2}}\cdot \sum_{m\ge l}p^{-(k^2+2k+2)m-l}\\
&=(1-p^{-1})\sum_{k\ge0}(-1)^k(1+p^{-k-1})
\frac{p^{-\frac{k^2+k}{2}}\;p^{-l(k^2+2k+3)}}{1-p^{-(k^2+2k+2)}}.
\end{split}
\end{align}
This completes the proof.
\end{proof}

\begin{thm}
Let $l\ge 0$. Then, we have
$$\lim_{n\rightarrow\infty}\Var[Z_{p^l\Z_p,n}]=p^{-l}-p^{-2l}+(1-p^{-1})\sum_{k\ge0}(-1)^k(1+p^{-k-1})
\frac{p^{-\frac{k^2+k}{2}}\;p^{-l(k^2+2k+3)}}{1-p^{-(k^2+2k+2)}}.$$
\end{thm}

\begin{proof}
Following the previous computations in \Cref{thm: 1pt expectation in plZp} and \Cref{thm: 2pt expectation in plZp}, we have
\begin{align}
\begin{split}
\lim_{n\rightarrow\infty}\Var[Z_{p^l\Z_p,n}]&=\lim_{n\rightarrow\infty}\E[Z_{p^l\Z_p,n}^2]-\lim_{n\rightarrow\infty}\E[Z_{p^l\Z_p,n}]^2\\
&=\lim_{n\rightarrow\infty}\E[Z_{p^l\Z_p\times p^l\Z_p,n}]+\lim_{n\rightarrow\infty}\E[Z_{p^l\Z_p,n}]-p^{-2l}\\
&=p^{-l}-p^{-2l}+\lim_{n\rightarrow\infty}\E[Z_{p^l\Z_p\times p^l\Z_p,n}]\\
&=p^{-l}-p^{-2l}+(1-p^{-1})\sum_{k\ge0}(-1)^k(1+p^{-k-1})
\frac{p^{-\frac{k^2+k}{2}}\;p^{-l(k^2+2k+3)}}{1-p^{-(k^2+2k+2)}}.
\end{split}
\end{align}
\end{proof}

We next turn to three-point correlations. In contrast with the two-point correlation function in \Cref{thm: 2pt in Zp}, we do not obtain an exact formula, but the following bounds will be sufficient for our applications. For the rest of the paper, we use the standard $q$-Pochhammer notation $$(a;q)_n:=\prod_{j=0}^{n-1}(1-aq^j),\quad n\in\mathbb Z_{\geq 0},$$ 
with $(a;q)_0:=1$, and write $(a;q)_\infty:=\prod_{j=0}^\infty(1-aq^j)$.

\begin{thm}\label{thm: 3pt correlation in pZp}
Suppose $x_1,x_2,x_3$ are distinct elements in $p\Z_p$. Then, we have
$$\frac{\rho_{\Q_p\times\Q_p\times\Q_p}^{(\infty)}(x_1,x_2,x_3)}{|(x_1-x_2)(x_2-x_3)(x_3-x_1)|}\in\left[\frac{(p^{-1};p^{-1})_\infty}{(1-p^{-1})^3},\frac{(p^{-1};p^{-1})_3}{(1-p^{-1})^3}\right].$$
\end{thm}

The following lemma will be useful for our proof of \Cref{thm: 3pt correlation in pZp}.

\begin{lemma}\label{lem: expected determinant 3pt}
Suppose $x_1,x_2,x_3$ are distinct elements in $p\Z_p$. For each $n\ge 1$, let $A\in\Mat_n(\Z_p)$ be Haar distributed. Then
$$(p^{-1};p^{-1})_\infty \leq \lim_{n\rightarrow\infty}\E|\det((x_1I_n-A)(x_2I_n-A)(x_3I_n-A))| \leq (p^{-1};p^{-1})_3.$$
\end{lemma}

\begin{proof}
For the lower bound, notice that for all fixed $n\ge 1$,
\begin{align}
\begin{split}
\E|\det((x_1I_n-A)(x_2I_n-A)(x_3I_n-A))|&\ge\Pr(A\in\GL_n(\Z_p))\\
&=\Pr(\overline A\in\GL_n(\F_p))\\
&=(p^{-1};p^{-1})_n.
\end{split}
\end{align}
Here, $\overline A\in\Mat_n(\F_p)$ is the reduction of $A$ mod $p$. Taking $n\rightarrow\infty$, we deduce that
$$\lim_{n\rightarrow\infty}\E|\det((x_1I_n-A)(x_2I_n-A)(x_3I_n-A))|\ge\lim_{n\rightarrow\infty}(p^{-1};p^{-1})_n=(p^{-1};p^{-1})_\infty.$$

For the upper bound, applying \cite[Corollary 8.5]{shen2026eigenvalues}\footnote{In the notation of \cite[Corollary 8.5]{shen2026eigenvalues}, we are taking $K=\Q_p$, so that $\mathcal{O}_K=\Z_p,r=1$, and $q=p$. Moreover, we are taking $k=3$.}, we have
$$\E|\det(A^3)|=\frac{(p^{-1};p^{-1})_3}{(p^{-n-1};p^{-1})_3}.$$
for any fixed $n$. By the H\H{o}lder inequality, we have that
\begin{align}
\begin{split}
\E|\det((x_1I_n-A)(x_2I_n-A)(x_3I_n-A))|&\le \prod_{i=1}^3(\E|\det(x_iI_n-A)^3|)^{1/3}\\
&=\E|\det(A^3)|\\
&=\frac{(p^{-1};p^{-1})_3}{(p^{-n-1};p^{-1})_3}.
\end{split}
\end{align}
Here, the second line holds because $A-x_iI_n$ has the same distribution as $A$. Taking $n\rightarrow\infty$, we complete the proof.
\end{proof}

\begin{proof}[Proof of \Cref{thm: 3pt correlation in pZp}]
In this case, we have $V(x_1)=V(x_2)=V(x_3)=\frac{1}{1-p^{-1}}$, and $$|\Delta_\sigma(x_1,x_2,x_3)|=|(x_1-x_2)(x_2-x_3)(x_3-x_1)|.$$ 
Following \Cref{prop: limiting correlation function}, we have
\begin{align}
\begin{split}
\rho_{\Q_p\times\Q_p\times\Q_p}^{(\infty)}(x_1,x_2,x_3)&=|\Delta_\sigma(x_1,x_2,x_3)|V(x_1)V(x_2)V(x_3)\\
&\times\lim_{n\rightarrow\infty}\E|\det((x_1I_n-A)(x_2I_n-A)(x_3I_n-A))|\\
&=\frac{|(x_1-x_2)(x_2-x_3)(x_3-x_1)|}{(1-p^{-1})^3}\\
&\times\lim_{n\rightarrow\infty}\E|\det((x_1I_n-A)(x_2I_n-A)(x_3I_n-A))|.
\end{split}
\end{align}
By the estimate in \Cref{lem: expected determinant 3pt}, we complete the proof.
\end{proof}

\begin{thm}\label{thm: 3pt zeros in pZp}
Let $l\ge 1$. Then, we have
\begin{equation}
\lim_{n\rightarrow\infty}\E[Z_{p^l\Z_p\times p^l\Z_p\times p^l\Z_p,n}]\in \left[\frac{p^{-6l}(1-p^{-1}+p^{-2})(p^{-1};p^{-1})_\infty}{(1-p^{-1})(1-p^{-2})(1-p^{-5})},\frac{p^{-6l}(1-p^{-1}+p^{-2})(p^{-1};p^{-1})_3}{(1-p^{-1})(1-p^{-2})(1-p^{-5})}\right].
\end{equation}
\end{thm}

\begin{lemma}\label{lem: integral of 3pt distance}
We have
$$
\int_{\Z_p^3}\bigl|(x_1-x_2)(x_2-x_3)(x_3-x_1)\bigr|dx_1dx_2dx_3
=
\frac{1-p^{-1}+p^{-2}}{(1+p^{-1})(1+p^{-1}+p^{-2}+p^{-3}+p^{-4})}.
$$
\end{lemma}

\begin{proof}
Consider the continuous group automorphism of $\Z_p^3$ given by
$$
(x_1,x_2,x_3)\mapsto (y_1,y_2,y_3):=(x_1-x_3,\;x_2-x_3,\;x_3).
$$
Its inverse is
$$
(y_1,y_2,y_3)\mapsto (y_1+y_3,\;y_2+y_3,\;y_3).
$$
This change of variables preserves Haar measure by translation invariance. Therefore,
$$
\int_{\Z_p^3}\bigl|(x_1-x_2)(x_2-x_3)(x_3-x_1)\bigr|dx_1dx_2dx_3=\int_{\Z_p^3} |y_1y_2(y_1-y_2)|dy_1dy_2dy_3.
$$
Since the integrand on the right-hand side is independent of $y_3$ and $\mu_{\Q_p}(\Z_p)=1$, this becomes
$$
\int_{\Z_p^3}\bigl|(x_1-x_2)(x_2-x_3)(x_3-x_1)\bigr|dx_1dx_2dx_3=\int_{\Z_p^2} |y_1y_2(y_1-y_2)|dy_1dy_2.
$$
For $m\ge 0$ let
$$
\mathcal{U}_m:=\{y\in \Z_p: \val(y)=m\},
$$
and note that 
$$
\mu_{\Q_p}(\mathcal{U}_m)=(1-p^{-1})p^{-m}.
$$
We decompose the integral according to whether $\val(y_1)\neq \val(y_2)$ or
$\val(y_1)=\val(y_2)$. Denote by $I_\mathrm{off}$ the contribution from case $\val(y_1)\neq \val(y_2)$, and $I_{\diag}$ the contribution from case $\val(y_1)=\val(y_2)$.

First suppose $m:=\val(y_1)<k:=\val(y_2)$. Write $y_1=p^m u$, $y_2=p^k w$ with
$u,w\in\Z_p^\times$. Then
$$
y_1-y_2=p^m\bigl(u-p^{k-m}w\bigr),
$$
and since $k-m\ge 1$ we have $u-p^{k-m}w\in\Z_p^\times$, hence
$\val(y_1-y_2)=m$ and
$$
|y_1y_2(y_1-y_2)|=p^{-m}p^{-k}p^{-m}=p^{-(2m+k)}.
$$
Therefore, the contribution from $m<k$ equals
$$
\sum_{m\ge 0}\sum_{k\ge m+1} \mu_{\Q_p}(\mathcal{U}_m)\mu_{\Q_p}(\mathcal{U}_k)p^{-(2m+k)}=(1-p^{-1})^2\sum_{m\ge 0}\sum_{k\ge m+1} p^{-(3m+2k)}.
$$
By symmetry the region $k<m$ contributes the same amount, so 
\begin{align*}
I_{\mathrm{off}}
&=2(1-p^{-1})^2\sum_{m\ge 0}p^{-3m}\sum_{k\ge m+1}p^{-2k}\\
&=2(1-p^{-1})^2\sum_{m\ge 0}p^{-3m}\cdot \frac{p^{-2(m+1)}}{1-p^{-2}}\\
&=\frac{2(1-p^{-1})^2p^{-2}}{(1-p^{-2})(1-p^{-5})}.
\end{align*}

Next suppose $m:=\val(y_1)=\val(y_2)$. Write $y_1=p^m u$, $y_2=p^m w$ with
$u,w\in\Z_p^\times$. Then
$$
|y_1y_2(y_1-y_2)|=p^{-m}p^{-m}\cdot p^{-m}|u-w|=p^{-3m}|u-w|.
$$
Hence the diagonal contribution is
$$
I_{\diag}
=\sum_{m\ge 0} \mu_{\Q_p}(\mathcal{U}_m)^2p^{-3m}C,
\qquad
C:=\frac{1}{(1-p^{-1})^2}\int_{\Z_p^\times\times \Z_p^\times} |u-w|dudw.
$$
To compute $C$, fix $u\in\Z_p^\times$. For $t\ge 1$,
$$
\mu_{\Q_p}\bigl(\{w\in\Z_p^\times: \val(u-w)\ge t\}\bigr)=p^{-t},
$$
since this is exactly one residue class modulo $p^t$ inside $\Z_p^\times$. Thus for $t\ge 1$,
$$
\mu_{\Q_p}\bigl(\{w\in\Z_p^\times: \val(u-w)=t\}\bigr)=p^{-t}-p^{-(t+1)}=(1-p^{-1})p^{-t},
$$
and the remaining mass is
$$
\mu_{\Q_p}\bigl(\{w\in\Z_p^\times: \val(u-w)=0\}\bigr)=(1-p^{-1})-p^{-1}=1-2p^{-1}.
$$
Therefore,
\begin{align*}
\int_{\Z_p^\times} |u-w|dw
&=(1-2p^{-1})\cdot 1+\sum_{t\ge 1} p^{-t}\cdot (1-p^{-1})p^{-t}\\
&=1-2p^{-1}+(1-p^{-1})\sum_{t\ge 1} p^{-2t}\\
&=\frac{1-p^{-1}-p^{-2}}{1+p^{-1}}.
\end{align*}
Integrating over $u\in\Z_p^\times$ gives
$$
\int_{\Z_p^\times\times \Z_p^\times} |u-w|dudw=(1-p^{-1})\cdot \frac{1-p^{-1}-p^{-2}}{1+p^{-1}},
$$
hence
$$
C=\frac{1-p^{-1}-p^{-2}}{(1+p^{-1})(1-p^{-1})}=\frac{1-p^{-1}-p^{-2}}{1-p^{-2}}.
$$
Consequently, 
$$
I_{\diag}=\sum_{m\ge 0} (1-p^{-1})^2p^{-2m}\cdot p^{-3m}\cdot C=\frac{(1-p^{-1})^2}{1-p^{-5}}\cdot \frac{1-p^{-1}-p^{-2}}{1-p^{-2}}.
$$
Finally, we have
\begin{align*}
I&=I_{\mathrm{off}}+I_{\diag}\\
&=\frac{(1-p^{-1})^2}{(1-p^{-2})(1-p^{-5})}
\Bigl(2p^{-2}+1-p^{-1}-p^{-2}\Bigr)\\
&=\frac{(1-p^{-1})^2}{(1-p^{-2})(1-p^{-5})}(1-p^{-1}+p^{-2})\\
&=\frac{1-p^{-1}+p^{-2}}{(1+p^{-1})(1+p^{-1}+p^{-2}+p^{-3}+p^{-4})},
\end{align*}
as claimed.
\end{proof}

\begin{proof}[Proof of \Cref{thm: 3pt zeros in pZp}]
On the one hand, by \Cref{prop: limiting correlation function}, we have
\begin{align}
\begin{split}
\lim_{n\rightarrow\infty}\E[Z_{p^l\Z_p\times p^l\Z_p\times p^l\Z_p,n}]&=\int_{p^l\Z_p^3}\rho_{\Q_p\times\Q_p\times\Q_p}^{(\infty)}(x_1,x_2,x_3)dx_1dx_2dx_3\\
&\ge\frac{(p^{-1};p^{-1})_\infty}{(1-p^{-1})^3}\int_{p^l\Z_p^3}|(x_1-x_2)(x_2-x_3)(x_3-x_1)|dx_1dx_2dx_3\\
&=\frac{p^{-6l}(p^{-1};p^{-1})_\infty}{(1-p^{-1})^3}\int_{\Z_p^3}|(x_1-x_2)(x_2-x_3)(x_3-x_1)|dx_1dx_2dx_3\\
&=\frac{p^{-6l}(p^{-1};p^{-1})_\infty}{(1-p^{-1})^3}\cdot\frac{1-p^{-1}+p^{-2}}{(1+p^{-1})(1+p^{-1}+p^{-2}+p^{-3}+p^{-4})}\\
&=\frac{p^{-6l}(1-p^{-1}+p^{-2})(p^{-1};p^{-1})_\infty}{(1-p^{-1})(1-p^{-2})(1-p^{-5})}.
\end{split}
\end{align}
Here, the second line comes from the lower bound of $\rho_{\Q_p\times\Q_p\times\Q_p}^{(\infty)}(x_1,x_2,x_3)$ given in \Cref{thm: 3pt correlation in pZp}, and the fourth line comes from \Cref{lem: integral of 3pt distance}. On the other hand, if we instead use the upper bound of $\rho_{\Q_p\times\Q_p\times\Q_p}^{(\infty)}(x_1,x_2,x_3)$ in \Cref{thm: 3pt correlation in pZp}, we can obtain the upper bound of the limiting expectation $\lim_{n\rightarrow\infty}\E[Z_{p^l\Z_p\times p^l\Z_p\times p^l\Z_p,n}]$ analogously.
\end{proof}

\subsection{Eigenvalues in quadratic extensions}

In this subsection, we are only interested in the case $p>2$. Under this assumption, there is exactly one unramified quadratic extension of $\Q_p$, namely $K=\Q_p[\sqrt{d}]/\Q_p$ where $d\in\Z$ such that $d\pmod p$ is a non-square in $\F_p:=\Z/p\Z$. There are two ramified quadratic extensions, which are $\Q_p[\sqrt{p}]/\Q_p$ and $\Q_p[\sqrt{dp}]/\Q_p$. 

\begin{thm}[{\cite[Theorem 1.8]{shen2026eigenvalues}}]\label{thm: quadratic correlation function}
Let $K/\mathbb Q_p$ be an unramified quadratic extension, and choose $x_1\in\mathcal O_K$ such that $\mathcal O_K=\mathbb Z_p[x_1]$. Then the limiting density function is given by
    \begin{equation}
        \rho_K^{(\infty)}(a_0+a_1x_1) = \frac{p^{-m}(1+p^{-1}-2p^{-m-1})}{1-p^{-1}}\sum_{k \geq 0} p^{-\frac{(2m+1)k^2+(4m+1)k}{2}}(1-p^{-k-1}),
    \end{equation}
    for any $a_0,a_1 \in \Z_p$ with $|a_1|=p^{-m}>0$. 

    In the case where $K$ is ramified, choosing $x_1 \in K$ with $|x_1|=p^{-1/2}$, we instead have 
    \begin{equation}\label{eq: correlation function_quadratic ramified}
        \rho_K^{(\infty)}(a_0+a_1x_1) = |\Disc_{K/\Q_p}|\frac{p^{-m}(1-p^{-m-1})}{1-p^{-1}}\sum_{k \geq 0} p^{-(m+1)k^2 - (2m+1)k}(1-p^{-2k-2})
    \end{equation}
    for any $a_0,a_1 \in \Z_p$ with $|a_1|=p^{-m}>0$. 
\end{thm}

Note that \Cref{thm: quadratic correlation function} is a purely random matrix statement, valid for every prime $p$. The Ellenberg-Jain-Venkatesh heuristic, however, only concerns the case where $p$ is odd. In this case, if $K/\Q_p$ is ramified quadratic, then $K/\Q_p$ is tamely ramified and the discriminant factor $|\Disc_{K/\Q_p}|$ in \eqref{eq: correlation function_quadratic ramified} is just $p^{-1}$.

\begin{thm}\label{thm: expected number in quadratic extensions}
Let $K/\Q_p$ be a quadratic extension. For the unramified case, for all $l\ge 0$, we have
\begin{multline}\label{eq: eigenvalues in pl_unramified}
\lim_{n\rightarrow\infty}\E[Z_{p^l\mathcal{O}_K^{\new},n}]=\sum_{k\ge 0}(1-p^{-k-1})p^{-\frac{k^2+k}{2}}\\
\cdot\frac{
(1+p^{-1})p^{-l(k^2+2k+3)}\bigl(1-p^{-(k^2+2k+3)}\bigr)
-2p^{-1}p^{-l(k^2+2k+4)}\bigl(1-p^{-(k^2+2k+2)}\bigr)
}{
\bigl(1-p^{-(k^2+2k+2)}\bigr)\bigl(1-p^{-(k^2+2k+3)}\bigr)
}.
\end{multline}

For the ramified case, choosing $x_1\in K$ to be generator with $|x_1|=p^{-1/2}$, we instead have
\begin{multline}\label{eq: eigenvalues in xi2l_ramified}
\lim_{n\rightarrow\infty}\E[Z_{p^l\mathcal{O}_K^{\new},n}]=|\Disc_{K/\Q_p}| \cdot p^{-l}\sum_{k\ge0}(1-p^{-2k-2})p^{-(k^2+k)}\\
\cdot\frac{p^{-l(k^2+2k+2)}\bigl(1-p^{-(k^2+2k+3)}\bigr)-p^{-1}p^{-l(k^2+2k+3)}\bigl(1-p^{-(k^2+2k+2)}\bigr)}{\bigl(1-p^{-(k^2+2k+2)}\bigr)\bigl(1-p^{-(k^2+2k+3)}\bigr)},
\end{multline}
\begin{multline}\label{eq: eigenvalues in xi2l+1_ramified}
\lim_{n\rightarrow\infty}\E[Z_{(x_1^{2l+1}\mathcal{O}_K)\cap \mathcal{O}_K^{\new},n}]=|\Disc_{K/\Q_p}| \cdot p^{-l-1}\sum_{k\ge0}(1-p^{-2k-2})p^{-(k^2+k)}\\
\cdot\frac{p^{-l(k^2+2k+2)}\bigl(1-p^{-(k^2+2k+3)}\bigr)-p^{-1}p^{-l(k^2+2k+3)}\bigl(1-p^{-(k^2+2k+2)}\bigr)}{\bigl(1-p^{-(k^2+2k+2)}\bigr)\bigl(1-p^{-(k^2+2k+3)}\bigr)}.
\end{multline}
\end{thm}

\begin{proof}
Notice that for all $m\ge l$, the subset 
$$
\mathcal{U}_{m,l}:=\{b_0+b_1x_1:b_0,b_1\in\Z_p,\val(b_0)\ge l,\val(b_1)=m\}\subset p^l\mathcal{O}_K
$$
has Haar probability measure $\mu_K(\mathcal{U}_{m,l})=(1-p^{-1})p^{-m-l}$.
Moreover, by \Cref{thm: quadratic correlation function}, the limiting correlation function $\rho_K^{(\infty)}$ takes value $\rho_K^{(\infty)}(p^mx_1)$ in $\mathcal{U}_{m,l}$. For the unramified case, following \Cref{prop: limiting correlation function} and the correlation function given in \Cref{thm: quadratic correlation function}, we have

\begin{align}
\begin{split}
\lim_{n\rightarrow\infty}\E[Z_{p^l\mathcal{O}_K^{\new},n}]&=\int_{p^l\mathcal{O}_K}\rho_{K}^{(\infty)}(x)dx\\
&=\sum_{m\ge l}\mu_K(\mathcal{U}_{m,l})\rho_K^{(\infty)}(p^mx_1)\\
&=\sum_{m\ge l}p^{-l}(1+p^{-1}-2p^{-m-1})\sum_{k\ge 0}(1-p^{-k-1})p^{-\frac{(2m+1)k^2+(4m+1)k}{2}-2m}\\
&=\sum_{k\ge 0}(1-p^{-k-1})p^{-\frac{k^2+k}{2}}\sum_{m\ge l}p^{-(k^2+2k+2)m-l}(1+p^{-1}-2p^{-m-1})\\
&=\text{RHS}\eqref{eq: eigenvalues in pl_unramified}.
\end{split}
\end{align}

Similarly, when $K/\Q_p$ is ramified, we have 
\begin{align}
\begin{split}
\lim_{n\rightarrow\infty}\E[Z_{p^l\mathcal{O}_K^{\new},n}]&=\int_{p^l\mathcal{O}_K}\rho_{K}^{(\infty)}(x)dx\\
&=\sum_{m\ge l}\mu_K(\mathcal{U}_{m,l})\rho_K^{(\infty)}(p^mx_1)\\
&=|\Disc_{K/\Q_p}|\sum_{m\ge l}p^{-l}(1-p^{-m-1})\sum_{k\ge 0}(1-p^{-2k-2})p^{-(m+1)k^2-(2m+1)k-2m}\\
&=|\Disc_{K/\Q_p}|\sum_{k\ge 0}(1-p^{-2k-2})p^{-(k^2+k)}\sum_{m\ge l}p^{-(k^2+2k+2)m-l}(1-p^{-m-1})\\
&=\text{RHS}\eqref{eq: eigenvalues in xi2l_ramified}.
\end{split}
\end{align}

Moreover,
\begin{align}
\begin{split}
\lim_{n\rightarrow\infty}\E[Z_{(x_1^{2l+1}\mathcal{O}_K)\cap \mathcal{O}_K^{\new},n}]&=\int_{x_1^{2l+1}\mathcal{O}_K}\rho_{K}^{(\infty)}(x)dx\\
&=\sum_{m\ge l}\mu_K(\mathcal{U}_{m,l+1})\rho_K^{(\infty)}(p^mx_1)\\
&=p^{-1}\cdot\sum_{m\ge l}\mu_K(\mathcal{U}_{m,l})\rho_K^{(\infty)}(p^mx_1)\\
&=p^{-1}\cdot\lim_{n\rightarrow\infty}\E[Z_{p^l\mathcal{O}_K^{\new},n}]\\
&=\text{RHS}\eqref{eq: eigenvalues in xi2l+1_ramified}.
\end{split}
\end{align}
This completes the proof.
\end{proof}

\subsection{Eigenvalues in cubic extensions}

Throughout this subsection, finite extensions of $\mathbb Q_p$ are regarded as subfields of a fixed algebraic closure $\bar{\Q}_p$. Thus conjugate non-Galois subfields are counted separately when we sum over extensions. There is exactly one unramified cubic extension of $\Q_p$, namely $\Q_p[\zeta]$, where $\zeta:=\zeta_{p^3-1}$ is a primitive $(p^3-1)^{th}$ root of unity. Moreover, for the ramified case,

\begin{enumerate}
\item When $p\equiv 2\pmod 3$, there are three ramified cubic extensions of $\Q_p$, namely 
$$\Q_p[p^{1/3}]/\Q_p,\Q_p[\xi p^{1/3}]/\Q_p,\Q_p[\xi^2 p^{1/3}]/\Q_p,\quad \xi^3=1,\xi\ne 1.$$ 
Moreover, these extensions are non-Galois, generated by the same polynomial $x^3-p$, tamely ramified, and have discriminant $p^{-2}$.
\item When $p\equiv 1\pmod 3$, there are three ramified cubic extensions of $\Q_p$, namely
$$\Q_p[p^{1/3}]/\Q_p,\Q_p[(gp)^{1/3}]/\Q_p,\Q_p[(g^2p)^{1/3}]/\Q_p.$$
Here, $g\in\Z$ such that $g\pmod p$ does not lie in $(\F_p^\times)^3$. Moreover, these extensions are Galois, tamely ramified, and have discriminant $p^{-2}$.
\item When $p=3$, there are nine isomorphism classes of ramified cubic extensions of $\Q_3$. This follows by specializing the general list of degree-$p$ extensions of $\Q_p$ for odd $p$ given by Jones-Roberts \cite[Proposition 2.3.1]{JONES200680}: their classification lists all such extensions explicitly, and when plugging in $p=3$ this list gives exactly the nine ramified cubic extensions of $\Q_3$. For a separate explicit list in the $p=3$ case, see also \cite[Proposition 3.1]{dhar2024classifying}. We list the generating polynomials of these ramified extensions here:
$$x^3+6x+3,x^3+3x+3,\quad (|\Disc_{K/\Q_p}|=3^{-3}),$$
$$x^3-3x^2+21,x^3-3x^2+3,x^3-3x^2+12,x^3+3x^2+3,\quad (|\Disc_{K/\Q_p}|=3^{-4}),$$
$$x^3+3,x^3+21,x^3+12,\quad (|\Disc_{K/\Q_p}|=3^{-5}).$$
Moreover, the isomorphism classes generated by 
$$x^3-3x^2+3,x^3-3x^2+12,x^3-3x^2+21$$ 
are Galois extensions, so each of those polynomials generates one extension; the isomorphism classes generated by 
$$x^3+3x^2+3,x^3+6x+3,x^3+3x+3,x^3+3,x^3+12,x^3+21$$
are non-Galois extensions, so each of those polynomials generates three extensions. Together we have $21$ cubic extensions.
\end{enumerate}

The following result is essentially a restatement of an orbital integral calculation in \cite{cho2024stable}, as we will explain in the proof.

\begin{thm}\label{thm: cubic orbital integral}
Let $K=\Q_p[\zeta]/\Q_p$ be the unramified cubic extension, where $\zeta=\zeta_{p^3-1}$. Let $a,b,c\in\Z_p$ with $(b,c)\neq(0,0)$, and set $m:=\min\{\val(b),\val(c)\}$. Then, we have
$$\#(\Lambda\backslash\Mod_{\Z_p[a+b\zeta+c\zeta^2]})=\frac{p^{3m}}{(1-p^{-1})(1-p^{-2})}\left(\frac{(p+1)(p^2+p+1)}{p^3}-\frac{3(p^2+p+1)}{p^{m+3}}+\frac{3}{p^{3m+3}}\right).$$

In the case where $K/\Q_p$ is instead a ramified cubic extension, choosing $\pi\in K$ to be a generator with $|\pi|=p^{-1/3}$, we instead have
$$\#(\Lambda\backslash\Mod_{\Z_p[a+b\pi+c\pi^2]})=\begin{cases}
\frac{p^{3m}}{(1-p^{-1})(1-p^{-2})}\left(\frac{p+1}{p}-\frac{2p+1}{p^{m+2}}+\frac{1}{p^{3m+3}}\right) & \val(b)=m,\val(c)\ge m\\
\frac{p^{3m+1}}{(1-p^{-1})(1-p^{-2})}\left(\frac{p+1}{p}-\frac{p+2}{p^{m+2}}+\frac{1}{p^{3m+4}}\right) 
& \val(b)>m,\val(c)=m \\
\end{cases}.$$
\end{thm}

\begin{lemma}\label{lem: order of quotient module}
Let $K=\Q_p[\zeta]/\Q_p$ be the unramified cubic extension, where $\zeta=\zeta_{p^3-1}$. Let $a,b,c\in\Z_p$ with $(b,c)\neq(0,0)$, and set $m:=\min\{\val(b),\val(c)\}$. Then, we have
$$\#(\mathcal{O}_K/\Z_p[a+b\zeta+c\zeta^2])=p^{3m}.$$
In the case where $K/\Q_p$ is ramified cubic, choosing $\pi\in K$ to be a generator with $|\pi|=p^{-1/3}$, we instead have
$$\#(\mathcal{O}_K/\Z_p[a+b\pi+c\pi^2])=\begin{cases}
p^{3m} & \val(b)=m,\val(c)\ge m\\
p^{3m+1} & \val(b)>m,\val(c)=m
\end{cases}.$$
\end{lemma}

\begin{proof}
For the unramified case, we take $u\in\mathcal{O}_K$ such that 
$$a+p^mu=a+b\zeta+c\zeta^2,\quad m:=\min\{\val(b),\val(c)\}.$$
In this case, we have
$$\mathcal{O}_K=\Z_p[u]=\Z_p+\Z_pu+\Z_p u^2,\quad\Z_p[a+b\zeta+c\zeta^2]=\Z_p[p^mu]=\Z_p+\Z_p\cdot p^mu+\Z_p \cdot p^{2m}u^2,$$
so that $\#(\mathcal{O}_K/\Z_p[a+b\zeta+c\zeta^2])=p^{3m}$. For the ramified case, we take $u\in\mathcal{O}_K$ such that 
$$a+p^mu=a+b\pi+c\pi^2,\quad m:=\min\{\val(b),\val(c)\}.$$
When $\val(b)=m,\val(c)\ge m$, we have $\val(u)=1/3,\val(u^2)=2/3$. Therefore,
$$\mathcal{O}_K=\Z_p+\Z_pu+\Z_p u^2,\quad\Z_p[a+b\pi+c\pi^2]=\Z_p[p^mu]=\Z_p+\Z_p\cdot p^mu+\Z_p \cdot p^{2m}u^2,$$
so that $\#(\mathcal{O}_K/\Z_p[a+b\pi+c\pi^2])=p^{3m}$. When $\val(b)>m,\val(c)=m$, we have $\val(u)=2/3,\val(p^{-1}u^2)=1/3$. Therefore,
$$\mathcal{O}_K=\Z_p+\Z_pu+\Z_p \cdot p^{-1}u^2,\quad\Z_p[a+b\pi+c\pi^2]=\Z_p[p^mu]=\Z_p+\Z_p\cdot p^mu+\Z_p \cdot p^{2m}u^2,$$
so that $\#(\mathcal{O}_K/\Z_p[a+b\pi+c\pi^2])=p^{3m+1}$. 
\end{proof}

\begin{proof}[Proof of {\Cref{thm: cubic orbital integral}}]
Before showing the details, we briefly explain how to extract the desired lattice count from the computations in
Yun \cite{yun2013orbital} and Cho-Kang-Lee \cite{cho2024stable}. The overall
strategy is as follows. Yun defines an orbital integral $O_\gamma$ in \cite[Section 4.2]{yun2013orbital} and proves in \cite[(4.1)]{yun2013orbital} that it is equal to the number of orbits of a certain lattice set. In our setting, this lattice set is precisely $\Mod_{\mathbb Z_p[a+b\zeta+c\zeta^2]}$ in the unramified case and $\Mod_{\mathbb Z_p[a+b\pi+c\pi^2]}$ in the ramified case, so Yun's formula identified $O_\gamma$ with our desired lattice counts
$$
\#\bigl(
\Lambda\backslash
\Mod_{\mathbb Z_p[a+b\zeta+c\zeta^2]}
\bigr)
\qquad \text{and} \qquad
\#\bigl(
\Lambda\backslash
\Mod_{\mathbb Z_p[a+b\pi+c\pi^2]}
\bigr).
$$
The same orbital integral $O_\gamma$ was denoted by $\SO_{\gamma,d\mu}$ in
\cite{cho2024stable}. Moreover, Cho-Kang-Lee explicitly computed a modified version of this orbital integral, which they called the normalized stable orbital integral and denoted by $\SO_\gamma$. In \cite[Section 2.2]{cho2024stable}, they compared the two versions $\SO_{\gamma,d\mu}$ and
$\SO_\gamma$, which are connected by the explicit normalization factor in
\cite[Proposition 2.4]{cho2024stable}. Thus, after identifying Yun's lattice
count with $\SO_{\gamma,d\mu}$, we will use \cite[Proposition 2.4]{cho2024stable} to express $\SO_\gamma$ in terms of this
lattice count and therefore solve for the lattice count using the explicit formula for $\SO_\gamma$ from \cite[Theorem 6.5]{cho2024stable}.

Now let us give the details. We take a matrix $\gamma\in\Mat_3(\Z_p)$ whose eigenvalues are the three Galois conjugates of $a+b\zeta+c\zeta^2$ (resp. $a+b\pi+c\pi^2$) in the unramified (resp. ramified) case. In the notation of \cite[Section 4.2]{yun2013orbital}, Yun considers 
$$
X_\gamma=\{g\in G(F)/G(\mathcal O_F): g^{-1}\gamma g\in \mathfrak g(\mathcal O_F)\}
$$
for a nonarchimedean local field $F$. Here, $G=\GL_n$, and $\mathfrak g$ is the corresponding Lie algebra which we denote
$\Mat_n$. Notice that the quotient $G(F)/G(\mathcal O_F)$ can be identified with the set of $\mathcal O_F$-lattices in $F^n$ by sending $g$ to $g\mathcal O_F^n$. In particular, under this identification, the condition $g^{-1}\gamma g\in \mathfrak g(\mathcal O_F)$ says exactly that the lattice $g\mathcal O_F^n$ is stable under the action of $\gamma$. In our situation,
$F=\Q_p$, $\mathcal O_F=\mathbb Z_p$, and $n=3$, so $G=\GL_3$ and
$\mathfrak g=\Mat_3$. Moreover, after identifying $\Q_p[\gamma]=K$, the action of $\gamma$ becomes a multiplication by the corresponding element of $K$. Therefore, $X_\gamma$ in our case is exactly the set denoted earlier by $\Mod_{\mathbb Z_p[a+b\zeta+c\zeta^2]}$ (resp. $\Mod_{\mathbb Z_p[a+b\pi+c\pi^2]}$) in the unramified (resp. ramified) case.

We now use the lattice interpretation of orbital integrals. It is shown in \cite[(4.1)]{yun2013orbital} that
\begin{equation}\label{eq: Yun lattice formula}
O_\gamma=\#(\Lambda_\gamma\backslash X_\gamma),
\end{equation}
where $\Lambda_\gamma$ there is exactly our $\Lambda$, and $O_\gamma$ is a certain orbital integral, the definition of which is given in \cite{yun2013orbital} but is not relevant for our purposes. In \cite{cho2024stable}, the same orbital integral $O_\gamma$ was denoted by
$\SO_{\gamma,d\mu}$. Hence, using the notation of \cite{cho2024stable}, the equality
\eqref{eq: Yun lattice formula} specializes to the form
\begin{equation}\label{eq: Yun lattice formula in CKL notation}
\SO_{\gamma,d\mu}=
\begin{cases}
\#\bigl(
\Lambda\backslash
\Mod_{\mathbb Z_p[a+b\zeta+c\zeta^2]}
\bigr)
& \text{$K/\Q_p$ unramified}\\[4pt]
\#\bigl(
\Lambda\backslash
\Mod_{\mathbb Z_p[a+b\pi+c\pi^2]}
\bigr)
& \text{$K/\Q_p$ ramified}
\end{cases}.
\end{equation}
So it remains to give a formula for $\SO_{\gamma,d\mu}$, which comes from converting the explicit formula for $\SO_\gamma$ established
in \cite[Theorem 6.5]{cho2024stable} into one for $\SO_{\gamma,d\mu}$, and then using
the identification \eqref{eq: Yun lattice formula in CKL notation}. Let us begin with the unramified case, where $K=\Q_p[\zeta]$, $\zeta=\zeta_{p^3-1}$, and $m=\min\{\val(b),\val(c)\}$. In this case the notation $F,\kappa,q,F_1,n,d'$ of \cite{cho2024stable} specializes to ours via
$$
F\mapsto\Q_p,\qquad
\kappa\mapsto\F_p,\qquad
q\mapsto p,\qquad
F_1\mapsto K,\qquad
n\mapsto3,\qquad
d'\mapsto m.
$$
In this case, $|\Delta_\gamma|$ in the notation of \cite{cho2024stable} specializes to our square of discriminant factor $|\Delta_\sigma(a+b\zeta+c\zeta^2)|^2$, and their field-discriminant factor $|\Delta_{F_1/F}|$ is $|\Disc_{K/\Q_p}|$ in our notation. Let $T_\gamma(\Q_p)$ be the centralizer of $\gamma$ inside $\GL_3(\Q_p)$. 

\cite{cho2024stable} treats a different orbital integral $\SO_{\gamma}$, which is related to ours in \cite[Proposition 2.4]{cho2024stable} through
\begin{equation}\label{eq:ckl_prop_2.4}
    \SO_\gamma = |\Delta_\sigma(a+b\zeta+c\zeta^2)| \cdot \frac{\#\GL_3(\F_p) p^{-\dim(\GL_3)}}{\#T_\gamma(\F_p) p^{-\dim(T_\gamma)}} |\Disc_{K/\Q_p}|^{-1/2} \cdot \SO_{\gamma,d\mu},
\end{equation}
and which is computed explicitly in \cite[Theorem 6.5]{cho2024stable} to be
\begin{equation}\label{eq:sogamma_unramified}
    \SO_{\gamma} = \frac{(p+1)(p^2+p+1)}{p^3}-
\frac{3(p^2+p+1)}{p^{m+3}}+\frac{3}{p^{3m+3}}. 
\end{equation}
So it remains to compute the other factors in \eqref{eq:ckl_prop_2.4}. 

We claim that after identifying $\Q_p[\gamma]$ with $K$, an element of $T_\gamma(\Q_p)$ is the same as multiplication by an element of $K^\times$. To see this, it suffices to show that the set of matrices in $\Mat_3(\Q_p)$ that commute with $\gamma$ is exactly $\Q_p[\gamma]$. On the one hand, it is clear that the set polynomials over $\gamma$ all commute with $\gamma$; on the other hand, select an arbitrary nonzero vector $v\in\Q_p^3$. Then, the vectors $v,\gamma v,\gamma^2v$ form a basis of the vector space $\Q_p^3$. Thus, for every matrix $B\in \Mat_3(\Q_p)$ there is a unique polynomial $f\in \Q_p[x]$ of degree
at most $2$ such that $Bv=f(\gamma)v$. Now assume in addition that $B\gamma=\gamma B$. Then 
$$
B\gamma v=\gamma Bv=\gamma f(\gamma)v=f(\gamma)\gamma v,\qquad
B\gamma^2 v=\gamma^2 Bv=\gamma^2 f(\gamma)v=f(\gamma)\gamma^2 v.
$$
Since $v,\gamma v,\gamma^2v$ is a basis, the two linear maps $B$ and
$f(\gamma)$ agree on a basis of $\Q_p^3$. Hence $B=f(\gamma)$. This verifies our claim. Now, since $K/\Q_p$ is unramified cubic, we have $\mathcal O_K/p\mathcal O_K\simeq \F_{p^3}$. Therefore, the reduction of $T_\gamma$ has $\F_p$-points identified with
$$
(\F_p\otimes_{\Z_p}\mathcal{O}_K)^\times\simeq(\mathcal O_K/p\mathcal O_K)^\times\simeq \F_{p^3}^\times,
$$
which has size $p^3-1$. Moreover, the dimension of $T_\gamma$ is exactly
$[K:\Q_p]=3$. Therefore, we have
$$
\#T_\gamma(\F_p)p^{-\dim T_\gamma}
=
(p^3-1)p^{-3}.
$$

We next simplify the discriminant factor. By \cite[Proposition 2.5]{cho2024stable}, we have
$$
|\Disc_{K/\Q_p}|^{1/2}=
\#(\mathcal O_K/\mathbb Z_p[a+b\zeta+c\zeta^2])
|\Delta_\sigma(a+b\zeta+c\zeta^2)|.
$$
By \Cref{lem: order of quotient module}, we have $\#(\mathcal O_K/\Z_p[a+b\zeta+c\zeta^2])=p^{3m}$, so that the above equality can be furthermore written as
$$
|\Disc_{K/\Q_p}|^{1/2}=p^{3m}|\Delta_\sigma(a+b\zeta+c\zeta^2)|.
$$
Combining these computations with \eqref{eq: Yun lattice formula in CKL notation}, \eqref{eq:ckl_prop_2.4} gives the lattice count 
\begin{align}
    \begin{split}
        \#\bigl(\Lambda\backslash
\Mod_{\mathbb Z_p[a+b\zeta+c\zeta^2]}
\bigr)&=p^{3m}\frac{(p^3-1)p^{-3}}
{\#\GL_3(\F_p)p^{-9}}\SO_\gamma
    \end{split}
\end{align}
and substituting in $\#\GL_3(\F_p)p^{-9}=(1-p^{-1})(1-p^{-2})(1-p^{-3})$ and \eqref{eq:sogamma_unramified} yields
\begin{equation}
    \label{eq:unramified_final_count}
\#\bigl(
\Lambda\backslash
\Mod_{\mathbb Z_p[a+b\zeta+c\zeta^2]}
\bigr)=
\frac{p^{3m}}{(1-p^{-1})(1-p^{-2})}
\left(
\frac{(p+1)(p^2+p+1)}{p^3}-\frac{3(p^2+p+1)}{p^{m+3}}+
\frac{3}{p^{3m+3}}
\right),
\end{equation}
which proves the desired formula in the unramified case.

We now treat the ramified case, where we again specialize
$$
F=\Q_p,\qquad\kappa=\F_p,\qquad q=p,\qquad F_1=K,\qquad n=3.
$$
As in the unramified case, let $T_\gamma$ be the group of invertible matrices in $\GL_3$ that commute with $\gamma$. As before, an element of $T_\gamma(\Q_p)$ can be regarded as multiplication by an element of $K^\times$. However, $K/\Q_p$ is now totally ramified cubic, so the invertible elements after reduction become
$$
(\F_p\otimes_{\Z_p}\mathcal{O}_K)^\times\simeq(\mathcal O_K/p\mathcal O_K)^\times=(\mathcal O_K/\pi^3\mathcal O_K)^\times,
$$
whose cardinality is $(p-1)p^2$. Moreover, the dimension of $T_\gamma$ is again $3$, and therefore
$$
\#T_\gamma(\F_p)p^{-\dim T_\gamma}=(p-1)p^2\cdot p^{-3}=(p-1)p^{-1}.
$$
Thus, \cite[Proposition 2.4]{cho2024stable} and \eqref{eq: Yun lattice formula in CKL notation} together give
\begin{equation}\label{eq: CKL normalization ramified lattice}
\SO_\gamma=|\Delta_\sigma(a+b\pi+c\pi^2)|
\frac{\#\GL_3(\F_p)p^{-9}}{(p-1)p^{-1}}|\Disc_{K/\Q_p}|^{-1/2} \#\bigl(
\Lambda\backslash
\Mod_{\mathbb Z_p[a+b\pi+c\pi^2]}
\bigr).
\end{equation}
By \cite[Proposition 2.5]{cho2024stable}, under our notation,
$$
|\Disc_{K/\Q_p}|^{1/2}=\#(\mathcal O_K/\mathbb Z_p[a+b\pi+c\pi^2])
|\Delta_\sigma(a+b\pi+c\pi^2)|.
$$
Regarding the values of $b$ and $c$, we need to further divide into two cases. First, suppose that $\val(b)=m$ and $\val(c)\ge m$. By
\Cref{lem: order of quotient module}, we have $\#(\mathcal O_K/\mathbb Z_p[a+b\pi+c\pi^2])=p^{3m}$. Substituting this into \eqref{eq: CKL normalization ramified lattice} gives the lattice
count 
$$
\#\bigl(
\Lambda\backslash
\Mod_{\mathbb Z_p[a+b\pi+c\pi^2]}
\bigr)=p^{3m}\frac{(p-1)p^{-1}}
{\#\GL_3(\F_p)p^{-9}}
\SO_\gamma=\frac{p^{3m}}{(1-p^{-2})(1-p^{-3})}\SO_\gamma.
$$
\cite[Theorem 6.5]{cho2024stable} also contains a formula for the ramified cubic case, which reads
$$
\SO_{\gamma}=\frac{(p+1)(p^2+p+1)}{p^3}-\frac{(2p+1)(p^2+p+1)}{p^{m+4}}+
\frac{p^2+p+1}{p^{3m+5}}.
$$
Therefore
$$
\#\bigl(\Lambda\backslash\Mod_{\mathbb Z_p[a+b\pi+c\pi^2]}\bigr)=\frac{p^{3m}}{(1-p^{-1})(1-p^{-2})}
\left(\frac{p+1}{p}-\frac{2p+1}{p^{m+2}}+\frac{1}{p^{3m+3}}\right).
$$

Next, suppose that $\val(b)>m$ and $\val(c)=m$. By
\Cref{lem: order of quotient module}, we have $\#(\mathcal O_K/\mathbb Z_p[a+b\pi+c\pi^2])=p^{3m+1}$. Therefore, the same normalization comparison gives
$$
\#\bigl(
\Lambda\backslash
\Mod_{\mathbb Z_p[a+b\pi+c\pi^2]}
\bigr)=p^{3m+1}
\frac{(p-1)p^{-1}}
{\#\GL_3(\F_p)p^{-9}}
\SO_\gamma=\frac{p^{3m+1}}{(1-p^{-2})(1-p^{-3})}\SO_\gamma.
$$
Again by \cite[Theorem 6.5]{cho2024stable}, in this case
$$
\SO_{\gamma}=\frac{(p+1)(p^2+p+1)}{p^3}-\frac{(p+2)(p^2+p+1)}{p^{m+4}}+\frac{p^2+p+1}{p^{3m+6}}
$$
and therefore
$$
\#\bigl(
\Lambda\backslash
\Mod_{\mathbb Z_p[a+b\pi+c\pi^2]}
\bigr)=\frac{p^{3m+1}}{(1-p^{-1})(1-p^{-2})}\left(\frac{p+1}{p}-\frac{p+2}{p^{m+2}}+\frac{1}{p^{3m+4}}
\right).
$$
This completes the proof of both ramified cases.
\end{proof}

The following proposition provides an estimate of the expected determinant term in \Cref{prop: limiting correlation function}.

\begin{prop}\label{prop: estimate of expectated det}
Let $K=\Q_p[\zeta]/\Q_p$ be the unramified cubic extension, where $\zeta=\zeta_{p^3-1}$. Take $a,b,c\in\Z_p$ with $\min\{\val(b),\val(c)\}=:m>0$, and let $Z$ be the minimal polynomial of $a+b\zeta+c\zeta^2$. Suppose $A\in\Mat_n(\Z_p)$ is Haar distributed. Then
$$\left|\lim_{n\rightarrow\infty}\E|\det(Z(A))|-(p^{-1};p^{-1})_3\right|\le(1-(p^{-1};p^{-1})_\infty)p^{-3m}.$$
In the case where $K/\Q_p$ is ramified cubic, choosing $\pi\in K$ to be a generator with $|\pi|=p^{-1/3}$. Take $a,b,c\in\Z_p$ with $\val(c)\ge\val(b)=:m>0$, and let $Z$ be the minimal polynomial of $a+b\pi+c\pi^2$. Then
$$\lim_{n\rightarrow\infty}\E|\det(Z(A))|-(p^{-1};p^{-1})_3\in[-(1-(p^{-1};p^{-1})_\infty)p^{-3m-1},(1-(p^{-1};p^{-1})_\infty)p^{-3m-1}].$$
If instead $\val(b) > \val(c) =: m \geq 0$, then
$$\lim_{n\rightarrow\infty}\E|\det(Z(A))|-(p^{-1};p^{-1})_3\in[-(1-(p^{-1};p^{-1})_\infty)p^{-3m-2},(1-(p^{-1};p^{-1})_\infty)p^{-3m-2}].$$
\end{prop}

\begin{proof}
We first prove the unramified case. For any fixed $n$, the two matrices $A-aI_n$ and $A-(a+b\zeta+c\zeta^2)I_n$ differ only in their diagonal entries, and each entry of the difference has norm $|b\zeta+c\zeta^2|\le p^{-m}$. Thus, by the non-archimedean triangle inequality, we obtain
$$\left| |\det(A-aI_n)|-|\det(A-(a+b\zeta+c\zeta^2)I_n)|\right|\le p^{-m}.$$
In particular, when $|\det(A-aI_n)|>p^{-m}$, we have 
$|\det(A-(a+b\zeta+c\zeta^2)I_n)|=|\det(A-aI_n)|$. Otherwise, when $|\det(A-aI_n)|\le p^{-m}$, we have $|\det(A-(a+b\zeta+c\zeta^2)I_n)|\le p^{-m}$. Therefore,
\begin{align}
\begin{split}
|\det(Z(A))|-|\det((A-aI_n)^3)|&=|\Nm_{K/\Q_p}\det(A-(a+b\zeta+c\zeta^2)I_n)|\\
&-|\det((A-aI_n)^3)|\\
&=|\det(A-(a+b\zeta+c\zeta^2)I_n)|^3-|\det(A-aI_n)|^3\\
&\in[-p^{-3m},p^{-3m}].
\end{split}
\end{align}
Therefore, applying \cite[Corollary 8.5]{shen2026eigenvalues}, we have that for fixed $n\ge 1$,
\begin{align}
\begin{split}
\left|\E|\det(Z(A))|-\frac{(p^{-1};p^{-1})_3}{(p^{-n-1};p^{-1})_3}\right|&=\left|\E|\det(Z(A))|-\E|\det((A-aI_n)^3)|\right|\\
&\le\Pr(|\det(A-aI_n)|\le p^{-m})p^{-3m}\\
&\le\Pr(\overline {A-aI_n}\notin\GL_n(\F_p))p^{-3m}\\
&=(1-(p^{-1};p^{-1})_n)p^{-3m}.
\end{split}
\end{align}
Letting $n$ go to infinity, we complete the proof.

For the ramified case, for any fixed $n$, suppose $\val(c)\ge \val(b)=m>0$. Since $\det(A-aI)\in\mathbb Q_p$, its valuation is integral, so $|\det(A-aI)|>p^{-m-1}$ is equivalent to $|\det(A-aI)|\ge p^{-m}$. Similarly to the unramified case, we have $|\det(A-aI_n)|>p^{-m-1}$ implies $|\det(A-(a+b\pi+c\pi^2)I_n)|=|\det(A-aI_n)|$, and $|\det(A-aI_n)|\le p^{-m-1}$ implies $|\det(A-(a+b\pi+c\pi^2)I_n)|\le p^{-m-1/3}$. Otherwise, suppose $\val(b)>\val(c)=m\ge 0$, we have $|\det(A-aI_n)|>p^{-m-1}$ implies $|\det(A-(a+b\pi+c\pi^2)I_n)|=|\det(A-aI_n)|$, and $|\det(A-aI_n)|\le p^{-m-1}$ implies $|\det(A-(a+b\pi+c\pi^2)I_n)|\le p^{-m-2/3}$. The rest of the proof follows by the same argument as in the unramified case.
\end{proof}

\begin{prop}\label{prop: estimate of cubic correlation function}
Let $K=\Q_p[\zeta]/\Q_p$ be the unramified cubic extension, where $\zeta=\zeta_{p^3-1}$. Denote by $m:=\min\{\val(b),\val(c)\}$. If $m=0$, we have
$\rho_K^{(\infty)}(a+b\zeta+c\zeta^2)=1$. If $m>0$, then $\rho_K^{(\infty)}(a+b\zeta+c\zeta^2)$ has lower bound and upper bound given by
$$\frac{(p^{-1};p^{-1})_3\pm(1-(p^{-1};p^{-1})_\infty)p^{-3m}}{(p^{-1};p^{-1})_3\cdot p^{3m}}\left(\frac{(p+1)(p^2+p+1)}{p^3}-\frac{3(p^2+p+1)}{p^{m+3}}+\frac{3}{p^{3m+3}}\right).$$

In the case where $K/\Q_p$ is ramified cubic, choosing $\pi\in K$ to be a generator with $|\pi|=p^{-1/3}$. If $\val(b)=0$, we have $\rho_K^{(\infty)}(a+b\pi+c\pi^2)=|\Disc_{K/\Q_p}|$. If $\val(c)\ge\val(b)=m>0$, then $\rho_K^{(\infty)}(a+b\pi+c\pi^2)$ has lower bound and upper bound given by
$$\frac{|\Disc_{K/\Q_p}|((p^{-1};p^{-1})_3\pm(1-(p^{-1};p^{-1})_\infty)p^{-3m-1})}{(1-p^{-1})^2(1-p^{-2})p^{3m}}\left(\frac{p+1}{p}-\frac{2p+1}{p^{m+2}}+\frac{1}{p^{3m+3}}\right).$$

If $\val(b)>\val(c)=m\ge 0$, then $\rho_K^{(\infty)}(a+b\pi+c\pi^2)$ has lower bound and upper bound given by
$$\frac{|\Disc_{K/\Q_p}|((p^{-1};p^{-1})_3\pm(1-(p^{-1};p^{-1})_\infty)p^{-3m-2})}{(1-p^{-1})^2(1-p^{-2})p^{3m+1}}\left(\frac{p+1}{p}-\frac{p+2}{p^{m+2}}+\frac{1}{p^{3m+4}}\right).$$
In the above, replacing the sign $\pm$ by $-$ (resp. $+$) gives the lower (resp. upper) bound.
\end{prop}

\begin{proof}
For the unramified case, when $m=0$, we have $\mathcal{O}_K=\Z_p[a+b\zeta+c\zeta^2]$. Therefore, by \Cref{thm: correlation for generator}, we have $\rho_K^{(\infty)}(a+b\zeta+c\zeta^2)=|\Disc_{K/\Q_p}|=1$. When $m>0$, denote by $Z$ the minimal polynomial of $a+b\zeta+c\zeta^2$. Applying \Cref{prop: limiting correlation function}, we have
\begin{align}
\begin{split}
\rho_K^{(\infty)}(a+b\zeta+c\zeta^2)&=|\Delta_\sigma(a+b\zeta+c\zeta^2)| \cdot \Den(a+b\zeta+c\zeta^2) \cdot \lim_{n\rightarrow\infty}\E|\det(Z(A))|\\
&=\frac{|\Delta_\sigma(a+b\zeta+c\zeta^2)|^2}{1-p^{-3}}\#(\Lambda\backslash\Mod_{\Z_p[a+b\zeta+c\zeta^2]})\cdot \lim_{n\rightarrow\infty}\E|\det(Z(A))|\\
&=\frac{|\Disc_{K/\Q_p}|\#(\Lambda\backslash\Mod_{\Z_p[a+b\zeta+c\zeta^2]})}{(1-p^{-3})\#(\mathcal{O}_K/\Z_p[a+b\zeta+c\zeta^2])^2}\cdot \lim_{n\rightarrow\infty}\E|\det(Z(A))|\\
&=\frac{p^{-6m}}{1-p^{-3}}\#(\Lambda\backslash\Mod_{\Z_p[a+b\zeta+c\zeta^2]})\cdot \lim_{n\rightarrow\infty}\E|\det(Z(A))|.
\end{split}
\end{align}
Here, the third line comes from \cite[Proposition 2.5]{cho2024stable} (or see also \cite[Page 10]{caruso2022zeroes}), and the last line comes from \Cref{lem: order of quotient module}. Now, applying the explicit form of $\#(\Lambda\backslash\Mod_{\Z_p[a+b\zeta+c\zeta^2]})$ given in \Cref{thm: cubic orbital integral} and the estimate of $\lim_{n\rightarrow\infty}\E|\det(Z(A))|$ given in \Cref{prop: estimate of expectated det}, we complete the proof. The ramified case follows similarly.
\end{proof}

We now integrate the cubic correlation estimates over small balls around $0$. Note that for the ramified case, we do not include here the region
$(\pi\mathcal{O}_K\backslash \pi^2\mathcal{O}_K)\cap\mathcal{O}_K^{\new}$
because this range is covered by the arbitrary-degree totally ramified calculation, which gives an exact formula rather than just upper and lower bounds; see \Cref{thm: totally ramified 1/m} later.

\begin{thm}\label{thm:cubic_small_ball_expectation}
Let $K=\Q_p[\zeta]/\Q_p$ be the unramified cubic extension, where $\zeta=\zeta_{p^3-1}$. Then, for all integers $l\ge 1$, the limiting expectation $\lim_{n\rightarrow\infty}\E[Z_{p^l\mathcal{O}_K^{\new},n}]$ has lower bound and upper bound given by
\begin{multline}\label{eq:cubic_unramified_small_ball}
\sum_{m\ge l}\frac{p^{-5m-l}}{(1-p^{-1})(1-p^{-3})}\left(\frac{(p+1)(p^2+p+1)}{p^3}-\frac{3(p^2+p+1)}{p^{m+3}}+\frac{3}{p^{3m+3}}\right)\\
\times((p^{-1};p^{-1})_3\pm(1-(p^{-1};p^{-1})_\infty)p^{-3m}).
\end{multline}

In the case where $K/\Q_p$ is ramified cubic, choosing $\pi\in K$ to be a generator with $|\pi|=p^{-1/3}$. For all integers $l\ge 1$, the limiting expectation $\lim_{n\rightarrow\infty}\E[Z_{(\pi^{3l-1}\mathcal{O}_K)\cap \mathcal{O}_K^{\new},n}]$ has lower and upper bound given by
\begin{multline}\label{eq:cubic_ramified_3lminus1_small_ball} 
\sum_{m\ge l}\frac{|\Disc_{K/\Q_p}|p^{-5m-l}}{(1-p^{-1})(1-p^{-2})}\left(\frac{p+1}{p}-\frac{2p+1}{p^{m+2}}+\frac{1}{p^{3m+3}}\right)\left((p^{-1};p^{-1})_3\pm\frac{1-(p^{-1};p^{-1})_\infty}{p^{3m+1}}\right)\\
+\sum_{m\ge l-1}\frac{|\Disc_{K/\Q_p}|p^{-5m-l-2}}{(1-p^{-1})(1-p^{-2})}\left(\frac{p+1}{p}-\frac{p+2}{p^{m+2}}+\frac{1}{p^{3m+4}}\right)
\left((p^{-1};p^{-1})_3\pm\frac{1-(p^{-1};p^{-1})_\infty}{p^{3m+2}}\right),
\end{multline}
the limiting expectation $\lim_{n\rightarrow\infty}\E[Z_{p^l\mathcal{O}_K^{\new},n}]$ has lower and upper bound given by
\begin{multline}
\sum_{m\ge l}\frac{|\Disc_{K/\Q_p}|p^{-5m-l}}{(1-p^{-1})(1-p^{-2})}\left(\frac{p+1}{p}-\frac{2p+1}{p^{m+2}}+\frac{1}{p^{3m+3}}\right)\left((p^{-1};p^{-1})_3\pm\frac{1-(p^{-1};p^{-1})_\infty}{p^{3m+1}}\right)\\
+\sum_{m\ge l}\frac{|\Disc_{K/\Q_p}|p^{-5m-l-2}}{(1-p^{-1})(1-p^{-2})}\left(\frac{p+1}{p}-\frac{p+2}{p^{m+2}}+\frac{1}{p^{3m+4}}\right)
\left((p^{-1};p^{-1})_3\pm\frac{1-(p^{-1};p^{-1})_\infty}{p^{3m+2}}\right),
\end{multline}
and the limiting expectation $\lim_{n\rightarrow\infty}\E[Z_{(\pi^{3l+1}\mathcal{O}_K)\cap \mathcal{O}_K^{\new},n}]$ has lower and upper bound given by
\begin{multline}
\sum_{m\ge l}\frac{|\Disc_{K/\Q_p}|p^{-5m-l-1}}{(1-p^{-1})(1-p^{-2})}\left(\frac{p+1}{p}-\frac{2p+1}{p^{m+2}}+\frac{1}{p^{3m+3}}\right)\left((p^{-1};p^{-1})_3\pm\frac{1-(p^{-1};p^{-1})_\infty}{p^{3m+1}}\right)\\
+\sum_{m\ge l}\frac{|\Disc_{K/\Q_p}|p^{-5m-l-3}}{(1-p^{-1})(1-p^{-2})}\left(\frac{p+1}{p}-\frac{p+2}{p^{m+2}}+\frac{1}{p^{3m+4}}\right)
\left((p^{-1};p^{-1})_3\pm\frac{1-(p^{-1};p^{-1})_\infty}{p^{3m+2}}\right).
\end{multline}
In the above, replacing the sign $\pm$ by $-$ (resp. $+$) gives the lower (resp. upper) bound.
\end{thm}

\begin{proof}
For the unramified case, notice that for all $m\ge l$, the set
$$\mathcal{U}_{m,l}:=\{a+b\zeta+c\zeta^2:\val(a)\ge l,\min\{\val(b),\val(c)\}=m\}$$
has Haar probability measure $\mu_K(\mathcal{U}_{m,l})=(1-p^{-2})p^{-2m-l}$. Moreover, in $\mathcal{U}_{m,l}$, the correlation function $\rho_K^{(\infty)}$ has lower and upper bound provided by \Cref{prop: estimate of cubic correlation function}. Therefore, by \Cref{prop: limiting correlation function}, we have
\begin{align}
\begin{split}
\lim_{n\rightarrow\infty}\E[Z_{p^l\mathcal{O}_K^{\new},n}]&=\int_{p^l\mathcal{O}_K}\rho_K^{(\infty)}(x)dx\\
&\ge\sum_{m\ge l}\mu_K(\mathcal{U}_{m,l})\frac{p^{-3m}}{(p^{-1};p^{-1})_3}\left(\frac{(p+1)(p^2+p+1)}{p^3}-\frac{3(p^2+p+1)}{p^{m+3}}+\frac{3}{p^{3m+3}}\right)\\
&\times((p^{-1};p^{-1})_3-(1-(p^{-1};p^{-1})_\infty)p^{-3m})\\
&=\sum_{m\ge l}\frac{p^{-5m-l}}{(1-p^{-1})(1-p^{-3})}\left(\frac{(p+1)(p^2+p+1)}{p^3}-\frac{3(p^2+p+1)}{p^{m+3}}+\frac{3}{p^{3m+3}}\right)\\
&\times((p^{-1};p^{-1})_3-(1-(p^{-1};p^{-1})_\infty)p^{-3m}).
\end{split}
\end{align}
The upper bound follows analogously. 

For the ramified case, notice that for all $m\ge l$, the set
$$\mathcal{U}_{m,\ge m,l}:=\{a+b\pi+c\pi^2:\val(a)\ge l,\val(b)=m,\val(c)\ge m\}$$
has Haar probability measure $\mu_K(\mathcal{U}_{m,\ge m,l})=p^{-2m-l}(1-p^{-1})$, and the set
$$\mathcal{U}_{>m,m,l}:=\{a+b\pi+c\pi^2:\val(a)\ge l,\val(b)>m,\val(c)=m\}$$
has Haar probability measure $\mu_K(\mathcal{U}_{>m,m,l})=p^{-2m-l-1}(1-p^{-1})$. In $\mathcal{U}_{m,\ge m,l}$ and $\mathcal{U}_{>m,m,l}$, the correlation function $\rho_K^{(\infty)}$ has lower and upper bound provided by \Cref{prop: estimate of cubic correlation function}. By this result and \Cref{prop: limiting correlation function}, we have 
\begin{align}
\begin{split}
\lim_{n\rightarrow\infty}\E[Z_{p^l\mathcal{O}_K^{\new},n}]&=\int_{p^l\mathcal{O}_K}\rho_K^{(\infty)}(x)dx\\
&\ge\sum_{m\ge l}\frac{\mu_K(\mathcal{U}_{m,\ge m,l})|\Disc_{K/\Q_p}|p^{-3m}}{(1-p^{-1})^2(1-p^{-2})}\left(\frac{p+1}{p}-\frac{2p+1}{p^{m+2}}+\frac{1}{p^{3m+3}}\right)\\
&\times((p^{-1};p^{-1})_3-(1-(p^{-1};p^{-1})_\infty)p^{-3m-1})\\
&+\sum_{m\ge l}\frac{\mu_K(\mathcal{U}_{>m, m,l})|\Disc_{K/\Q_p}|p^{-3m-1}}{(1-p^{-1})^2(1-p^{-2})}\left(\frac{p+1}{p}-\frac{p+2}{p^{m+2}}+\frac{1}{p^{3m+4}}\right)\\
&\times((p^{-1};p^{-1})_3-(1-(p^{-1};p^{-1})_\infty)p^{-3m-2})\\
&=\sum_{m\ge l}\frac{|\Disc_{K/\Q_p}|p^{-5m-l}}{(1-p^{-1})(1-p^{-2})}\left(\frac{p+1}{p}-\frac{2p+1}{p^{m+2}}+\frac{1}{p^{3m+3}}\right)\\
&\times((p^{-1};p^{-1})_3-(1-(p^{-1};p^{-1})_\infty)p^{-3m-1})\\
&+\sum_{m\ge l}\frac{|\Disc_{K/\Q_p}|p^{-5m-l-2}}{(1-p^{-1})(1-p^{-2})}\left(\frac{p+1}{p}-\frac{p+2}{p^{m+2}}+\frac{1}{p^{3m+4}}\right)\\
&\times((p^{-1};p^{-1})_3-(1-(p^{-1};p^{-1})_\infty)p^{-3m-2}).\\
\end{split}
\end{align}
The upper bound follows analogously. Similarly, the cases  
$$\lim_{n\rightarrow\infty}\E[Z_{(\pi^{3l-1}\mathcal{O}_K)\cap \mathcal{O}_K^{\new},n}],\quad\lim_{n\rightarrow\infty}\E[Z_{(\pi^{3l+1}\mathcal{O}_K)\cap \mathcal{O}_K^{\new},n}]$$
follow from the same decomposition as above, with only the lower bounds on the valuations of the coefficients shifted by one step. This completes the proof.
\end{proof}

\subsection{Eigenvalues in totally ramified extensions of arbitrary degree}

The preceding cases concern extensions of degree at most three. In arbitrary degree, a simple expression is still available for the one-point correlation function at a uniformizer of a totally ramified extension. Integrating this formula over the first annulus and summing over extensions gives the expectation needed for the corresponding numerical test.

\begin{thm}\label{thm: correlation function at generator}
Let $K$ be a totally ramified extension over $\Q_p$ of degree $m$, and $x_1\in\mathcal{O}_K$ such that $\val(x_1)=1/m$. Then, we have
$$\rho_{K}^{(\infty)}(x_1)=|\Disc_{K/\Q_p}|.$$
\end{thm}

\begin{proof}
Observe that for all $a_0,\ldots,a_{m-1}\in\Z_p$,
$$\val(a_0+a_1x_1+\cdots+a_{m-1}x_1^{m-1})=\min\{\val(a_0),\val(a_1)+1/m,\ldots,\val(a_{m-1})+(m-1)/m\}.$$
Therefore, we have
\begin{align}
\begin{split}
\mathcal{O}_K&=\{a_0+a_1x_1+\cdots+a_{m-1}x_1^{m-1}:\val(a_0+a_1x_1+\cdots+a_{m-1}x_1^{m-1})\ge 0\}\\
&=\{a_0+a_1x_1+\cdots+a_{m-1}x_1^{m-1}:a_0,\ldots,a_{m-1}\in\Z_p\}\\
&=\Z_p[x_1].
\end{split}
\end{align}
Therefore, we can apply \Cref{thm: correlation for generator}, which completes the proof.
\end{proof}

\begin{thm}\label{thm: totally ramified 1/m}
Let $K$ be a totally ramified extension over $\Q_p$ of degree $m$, and $\pi$ be a uniformizer of $K$. Then, we have
$$\lim_{n\rightarrow\infty}\E[Z_{\pi\mathcal{O}_K\backslash\pi^2\mathcal{O}_K,n}]=\frac{p-1}{p^2}|\Disc_{K/\Q_p}|.$$
\end{thm}

\begin{proof}
Notice that the subset $\pi\mathcal{O}_K\backslash\pi^2\mathcal{O}_K\subset \mathcal{O}_K$ has Haar probability measure 
$$\mu_K(\pi\mathcal{O}_K\backslash\pi^2\mathcal{O}_K)=\frac{p-1}{p^2}.$$
Moreover, by \Cref{thm: correlation function at generator}, the limiting correlation function $\rho_K^{(\infty)}$ takes value $|\Disc_{K/\Q_p}|$ in $\pi\mathcal{O}_K\backslash\pi^2\mathcal{O}_K$. Therefore, by \Cref{prop: limiting correlation function}, we have
$$\lim_{n\rightarrow\infty}\E[Z_{\pi\mathcal{O}_K\backslash\pi^2\mathcal{O}_K,n}]=\int_{\pi\mathcal{O}_K\backslash\pi^2\mathcal{O}_K}\rho_K^{(\infty)}(x)dx=\frac{p-1}{p^2}|\Disc_{K/\Q_p}|.$$
\end{proof}

\begin{thm}\label{thm: summing slope 1/m}
Let $m\ge 1$. Then, we have
$$\lim_{n\rightarrow\infty}\sum_K\E[Z_{\pi\mathcal{O}_K\backslash\pi^2\mathcal{O}_K,n}]=\frac{m(p-1)}{p^{m+1}}.$$
Here, the sum ranges over all totally ramified extensions $K/\Q_p$ of degree $m$ with conjugate subfields counted separately, and $\pi$ is a corresponding uniformizer.
\end{thm}

\begin{proof}
Following Serre's mass formula (see \cite[(4.4)]{caruso2022zeroes}) by taking $m=r,p=q$, we have
$$\sum_K|\Disc_{K/\Q_p}|=\frac{m}{p^{m-1}}.$$
Therefore, by \Cref{thm: totally ramified 1/m}, we have
\begin{align}
\begin{split}\lim_{n\rightarrow\infty}\sum_K\E[Z_{\pi\mathcal{O}_K\backslash\pi^2\mathcal{O}_K,n}]&=\sum_K\lim_{n\rightarrow\infty}\E[Z_{\pi\mathcal{O}_K\backslash\pi^2\mathcal{O}_K,n}]\\
&=\sum_K\frac{p-1}{p^2}|\Disc_{K/\Q_p}|\\
&=\frac{m(p-1)}{p^{m+1}}.
\end{split}
\end{align}
Here, we can change the order of taking the limit and summation over all extensions $K/\Q_p$, because there are only finitely many totally ramified extensions $K/\Q_p$ of degree $m$.
\end{proof}

\section{Distributional consequences of the distinguished polynomial conjecture}
\label{sec:distributional-consequences}

The purpose of this section is to prove various explicit formulas for asymptotic statistics of zeros of $p$-adic $L$-functions which are implied by \Cref{conj:distinguished_polynomial_intro}, and in most cases require the cruder robustness statement of \Cref{conj:bad_sets_don't_effect_averages_intro} as well. Most of these consequences, though not all, were already stated in the Introduction. We also give formulas for the third moment of the number of eigenvalues in $p^r \Z_p$ (\Cref{cor:Zp_third_moment}) and zeros in cubic extensions (\Cref{cor:cubic_extensions}) which we did not state earlier because the statements are a bit more technical.

Recall that in \Cref{sec: Random matrix computations}, we showed that the distinguished polynomials $P_{A_n}$ of additive Haar-random matrices in $\Mat_n(\Z_p)$ converge weakly to a random polynomial $P_{A_\infty}\in\mathscr P$, which is almost surely squarefree (\Cref{thm:existence of limiting distinguished polynomial}). We also showed that, after shifting the eigenvalues by $1$, the original multiplicative Haar model of Ellenberg-Jain-Venkatesh \cite{ellenberg2011modeling} has the same limiting distinguished polynomial (\Cref{thm:existence of limiting distinguished polynomial GL}). Finally, we explicitly computed formulas for various limiting quantities. The work in this section consists in combining those results with Conjectures~\ref{conj:distinguished_polynomial_intro} and~\ref{conj:bad_sets_don't_effect_averages_intro} to yield the aforementioned predictions.

We will adapt the notation in \Cref{sec: Introduction} and \Cref{sec:preliminaries}. The space $\mathscr P$ is the disjoint union of the spaces $\mathscr P_d$ of distinguished polynomials of degree $d$, with the coefficient topology on each $\mathscr P_d$. For each twist $i\in\{0,\ldots,p-2\}$, let $X_i$ denote the family of even characters $\chi=\theta\omega^i$, where $\theta$ is quadratic with primitive conductor coprime to $p$, and, when $i=1$, we additionally require $\theta(p)\ne1$. For $D>0$, let $X_i(D) \subset X_i$ be the finite subset of such characters where in addition the primitive conductor of $\theta$ is at most $D$. We use notation $\chi_D$ to denote a character chosen uniformly at random from $X_i(D)$, where the twist $i$ is implicit in the notation.

For $\chi\in X_i$, the distinguished polynomial $P_\chi\in\mathscr P$ records the zeros in the open unit disk of the associated Iwasawa power series, as explained in \Cref{sec:preliminaries}, and we write $\lambda(\chi):=\deg(P_\chi)$. On the random matrix side, we retain the notation $A_n$, $P_{A_n}$, and $P_{A_\infty}$ from \Cref{subsec:limiting distinguished polynomial}. We write $\val$ for the $p$-adic valuation, with $\val(0)=\infty$, and adopt the convention $\Disc(P)=1$ when $\deg(P)\le1$.

\subsection{Generalities} We first consider the degree and discriminant, which are determined by the coefficients of the distinguished polynomial. Their joint limiting distribution is as follows.

\begin{thm}
\label{thm:conj-implies-disc-lambda}
Fix $i\in\{0,\ldots,p-2\}$ and assume \Cref{conj:distinguished_polynomial_intro}. As $D\to\infty$, the joint random variable
$$
\left(\val(\Disc(P_{\chi_D})),\lambda(\chi_D)\right)
$$
converges in distribution to
$$
\left(\val(\Disc(P_{A_\infty})),\deg(P_{A_\infty})\right).
$$
\end{thm}

We next consider the locations of individual zeros and of tuples of zeros. For a finite extension $K/\Q_p$, let $\mathcal O_K^{\new}=\{x\in\mathcal O_K:\Q_p(x)=K\}$. Fix finite extensions $K_1,\ldots,K_m$ of $\Q_p$. For $U\subset\mathcal O_{K_1}^{\new}\times\cdots\times\mathcal O_{K_m}^{\new}$ and $P\in\mathscr P$, define
$$
\begin{aligned}
Z_U(P)=\#\{(x_1,\ldots,x_m)\in U:\;&P(x_1)=\cdots=P(x_m)=0,\\
&x_1,\ldots,x_m\text{ are pairwise non-conjugate over }\Q_p\}.
\end{aligned}
$$
Thus $Z_U(P)$ counts ordered tuples of zeros of $P$ whose entries belong to distinct Galois-conjugacy classes. When $m=1$, it is the number of distinct zeros in $U$, as in \Cref{sec: Introduction}; in this case we also write $Z_{U,\chi_D}=Z_U(P_{\chi_D})$. If $U$ lies in the open unit disk in every coordinate, then $Z_U(P_{A_n})$ agrees with the eigenvalue-counting variable $Z_{U,n}$ of \Cref{sec: Random matrix computations}.

\begin{thm}
\label{thm:conj-implies-local-zero-distribution}
Fix $i\in\{0,\ldots,p-2\}$ and assume \Cref{conj:distinguished_polynomial_intro}. Let $K_1,\ldots,K_m$ be finite extensions of $\Q_p$, and suppose that $U\subset\mathcal O_{K_1}^{\new}\times\cdots\times\mathcal O_{K_m}^{\new}$ is clopen in the relative topology. Then $Z_U(P_{\chi_D})$ converges in distribution, as $D\to\infty$, to $Z_U(P_{A_\infty})$.
\end{thm}

The proofs of these two theorems use the continuous mapping theorem. By \Cref{thm:existence of limiting distinguished polynomial}, $P_{A_\infty}$ is almost surely squarefree, so it remains to verify continuity of the relevant statistics at squarefree polynomials.

The next lemma verifies the required continuity statements. 

\begin{lemma}
\label{lem:disc-degree-continuity}
The map $P\mapsto(\val(\Disc(P)),\deg(P))$ from $\mathscr P$ to $(\Z_{\ge0}\cup\{\infty\})\times\Z_{\ge0}$ is locally constant, in particular continuous, at every squarefree polynomial.
\end{lemma}

\begin{proof}
Let $P\in\mathscr P$ be squarefree. Since $\mathscr P$ is a disjoint union of the spaces $\mathscr P_d$, there is an open neighborhood of $P$ on which the degree is constant, equal to $\deg(P)$. On this fixed-degree component, the discriminant is a polynomial in the coefficients, hence is continuous as a $\Z_p$-valued function of the coefficients. Since $\Disc(P)\ne0$, the valuation map is locally constant at $\Disc(P)$, so $\val(\Disc(Q))=\val(\Disc(P))$ for all $Q$ in some neighborhood of $P$ inside that component. Therefore both coordinates of $(\val(\Disc(Q)),\deg(Q))$ are constant for all $Q$ sufficiently close to $P$, proving the claim.
\end{proof}

For local zero statistics, we use the root-counting continuity result of \cite[Lemma 3.7]{shen2026resultant}. It applies to the clopen sets appearing in \Cref{thm:conj-implies-local-zero-distribution}.

\begin{lemma}
\label{lem:root-count-continuity}
Let $K_1,\ldots,K_m$ be finite extensions of $\Q_p$, and let $U\subset\mathcal O_{K_1}^{\new}\times\cdots\times\mathcal O_{K_m}^{\new}$ be clopen in the relative topology. Then the map $P\mapsto Z_U(P)$ from $\mathscr P$ to $\Z_{\ge0}$ is continuous at every squarefree polynomial.
\end{lemma}

\begin{proof}
This is the specialization of the root-counting continuity result \cite[Lemma 3.7]{shen2026resultant} to the distinguished-polynomial space $\mathscr P$. In the notation of that lemma, take $S=\{1\}$, so that $\mathscr P^{\{1\}}=\mathscr P$ and the corresponding lifted domain is the open unit disk.
\end{proof}

We can now prove the two distributional convergence theorems stated above. 

\begin{proof}[Proof of \Cref{thm:conj-implies-disc-lambda}]
\Cref{conj:distinguished_polynomial_intro} gives convergence in distribution of $P_{\chi_D}$ to $P_{A_\infty}$. By \Cref{lem:disc-degree-continuity}, the map $P\mapsto(\val(\Disc(P)),\deg(P))$ is continuous at every squarefree polynomial, and \Cref{thm:existence of limiting distinguished polynomial} shows that $P_{A_\infty}$ is squarefree with probability one. The continuous mapping theorem gives the asserted convergence.
\end{proof}

\begin{proof}[Proof of \Cref{thm:conj-implies-local-zero-distribution}]
By \Cref{lem:root-count-continuity}, the map $P\mapsto Z_U(P)$ is continuous at every squarefree polynomial. Thus the continuous mapping theorem, \Cref{conj:distinguished_polynomial_intro}, and \Cref{thm:existence of limiting distinguished polynomial} give convergence in distribution of $Z_U(P_{\chi_D})$ to $Z_U(P_{A_\infty})$.
\end{proof}

\Cref{thm:conj-implies-disc-lambda} and \Cref{thm:conj-implies-local-zero-distribution} give convergence in distribution, but not automatically convergence of moments. In order to obtain this, we first record a uniform degree bound for the random matrix model. 

\begin{lemma}[Exponential tail bound for the distinguished degree]
\label{lem:random-matrix-degree-tail}
For every $n\ge1$ and $0\le d\le n$, we have
$$
\Pr(\deg(P_{A_n})=d)=p^{-d}\prod_{j=d+1}^n(1-p^{-j}).
$$
In particular, for every $n\ge1$ and $d\ge0$,
$$
\Pr(\deg(P_{A_n})\ge d)\le\frac{p^{-d}}{1-p^{-1}}.
$$
The same bound holds for $P_{A_\infty}$.
\end{lemma}

\begin{proof}
Let $\overline A_n$ denote the reduction of $A_n$ modulo $p$, which is uniformly distributed on $\Mat_n(\mathbb F_p)$. By the definition of the distinguished factor, $\deg(P_{A_n})$ equals the dimension of the generalized zero-eigenspace of $\overline A_n$.

Fix $0\le d\le n$. By Jordan normal form (or the $\Z/p\Z$ case of the decomposition argument in the proof of \Cref{thm:existence of limiting distinguished polynomial GL}) the generalized zero-eigenspace of $\overline A_n$ has dimension $d$ if and only if there is a unique decomposition $\mathbb F_p^n=X\oplus Y$ into $\overline A_n$-invariant subspaces, with $\dim X=d$ and $\dim Y=n-d$, such that $\overline A_n|_X$ is nilpotent and $\overline A_n|_Y$ is invertible. There are
$$
\frac{|\GL_n(\mathbb F_p)|}{|\GL_d(\mathbb F_p)|\,|\GL_{n-d}(\mathbb F_p)|}
$$
choices of such an ordered pair of complementary subspaces. For each choice, there are $p^{d^2-d}$ nilpotent endomorphisms of $X$ by the theorem of Fine-Herstein \cite[Theorem 1]{fine1958probability}, and $|\GL_{n-d}(\mathbb F_p)|$ invertible endomorphisms of $Y$. Since $|\GL_j(\mathbb F_p)|=p^{j^2}\prod_{k=1}^j(1-p^{-k})$, we obtain
$$
\begin{aligned}
\Pr(\deg(P_{A_n})=d)
&=\frac{1}{p^{n^2}}\frac{|\GL_n(\mathbb F_p)|}{|\GL_d(\mathbb F_p)|\,|\GL_{n-d}(\mathbb F_p)|}
p^{d^2-d}|\GL_{n-d}(\mathbb F_p)|\\
&=p^{-d}\prod_{j=d+1}^n(1-p^{-j})\\
&\le p^{-d}.
\end{aligned}
$$
It then follows that
$$
\Pr(\deg(P_{A_n})\ge d)
=\sum_{j=d}^n\Pr(\deg(P_{A_n})=j)
\le\sum_{j\ge d}p^{-j}
=\frac{p^{-d}}{1-p^{-1}}.
$$
Finally, since $P_{A_n}$ converges in distribution to $P_{A_\infty}$, the same bound holds for $P_{A_\infty}$.
\end{proof}

On the $p$-adic $L$-function side, the analogous moment control does not follow from weak convergence alone. This is where \Cref{conj:bad_sets_don't_effect_averages_intro} is required. Together with \Cref{conj:distinguished_polynomial_intro}, it gives the following result, including the moment identity \eqref{eq:moments_match_chi_tilde} stated in \Cref{sec: Introduction}.

\begin{thm}[Convergence of degree moments]
\label{thm:moment-convergence-from-conjectures}
Assume \Cref{conj:distinguished_polynomial_intro} and \Cref{conj:bad_sets_don't_effect_averages_intro}, and fix a twist $i\in\{0,\ldots,p-2\}$. Let $B\subset X_i$ have asymptotically negligible density in the sense of \Cref{conj:bad_sets_don't_effect_averages_intro}, and, for sufficiently large $D$, let $\widetilde\chi_D$ be chosen uniformly from $X_i(D)\setminus B$. Then, for every integer $k\ge0$,
$$
\lim_{D\to\infty}\E[\deg(P_{\chi_D})^k]=\lim_{D\to\infty}\E[\deg(P_{\widetilde\chi_D})^k]=\E[\deg(P_{A_\infty})^k],
$$
and in particular both limits exist.
\end{thm}

\begin{proof}
The case $k=0$ is immediate. Fix $k\ge1$. By \Cref{lem:random-matrix-degree-tail}, the degree of $P_{A_\infty}$ has an exponential tail. In particular, $\E[\deg(P_{A_\infty})^k]$ is finite. We first prove
$$
\liminf_{D\to\infty}\E[\deg(P_{\chi_D})^k]\ge\E[\deg(P_{A_\infty})^k],
$$
and then prove the reverse inequality for the limit superior.

For every fixed integer $M\ge0$, the function $P\mapsto\deg(P)^k\mathbf{1}_{\{\deg(P)\le M\}}$ is bounded and continuous on $\mathscr P$. By \Cref{conj:distinguished_polynomial_intro},
$$
\lim_{D\to\infty}\E[\deg(P_{\chi_D})^k\mathbf{1}_{\{\deg(P_{\chi_D})\le M\}}]
=\E[\deg(P_{A_\infty})^k\mathbf{1}_{\{\deg(P_{A_\infty})\le M\}}].
$$
Since all terms are nonnegative,
$$
\liminf_{D\to\infty}\E[\deg(P_{\chi_D})^k]
\ge\E[\deg(P_{A_\infty})^k\mathbf{1}_{\{\deg(P_{A_\infty})\le M\}}].
$$
Letting $M\to\infty$ proves the lower bound.

We now prove the upper bound by contradiction. Suppose that there exists $\varepsilon>0$ such that
$$
\limsup_{D\to\infty}\E[\deg(P_{\chi_D})^k]
>\E[\deg(P_{A_\infty})^k]+3\varepsilon.
$$
We will construct a subset $B\subset X_i$ of asymptotically negligible density whose contribution to the $k$th moment does not tend to zero, contradicting \Cref{conj:bad_sets_don't_effect_averages_intro}.

Enumerate the characters in $X_i$ as $\chi_1,\chi_2,\ldots$ in order of increasing conductor of the associated quadratic character $\theta$, breaking ties arbitrarily. Write $N(D)=|X_i(D)|$.

We first choose $D_1$ sufficiently large and set $N_1=N(D_1)$ so that
$$
\frac1{N_1}\sum_{n=1}^{N_1}\deg(P_{\chi_n})^k
\ge\E[\deg(P_{A_\infty})^k]+2\varepsilon
$$
and
$$
\frac1{N_1}\sum_{\substack{1\le n\le N_1\\\deg(P_{\chi_n})\le1}}\deg(P_{\chi_n})^k
\le\E[\deg(P_{A_\infty})^k]+\varepsilon.
$$
The first inequality follows from our assumption, while the second follows from \Cref{conj:distinguished_polynomial_intro}.

Next, choose $D_2>D_1$ sufficiently large and set $N_2=N(D_2)$ so that
$$
\frac1{N_2-N_1}\sum_{n=N_1+1}^{N_2}\deg(P_{\chi_n})^k
\ge\E[\deg(P_{A_\infty})^k]+2\varepsilon
$$
and
$$
\frac1{N_2-N_1}\sum_{\substack{N_1<n\le N_2\\\deg(P_{\chi_n})\le2}}\deg(P_{\chi_n})^k
\le\E[\deg(P_{A_\infty})^k]+\varepsilon.
$$
To see why such a choice of $D_2$ is possible, recall that our assumption provides  arbitrarily large $D$ for which the full average over $X_i(D)$ is at least $\E[\deg(P_{A_\infty})^k]+3\varepsilon$. Since the first $N_1$ characters contribute a fixed finite sum, their contribution becomes negligible as $D\to\infty$. Thus, for sufficiently large such $D$, the average over the interval $(N_1,N(D)]$ is at least $\E[\deg(P_{A_\infty})^k]+2\varepsilon$. Similarly, by \Cref{conj:distinguished_polynomial_intro}, the average restricted to degree at most $2$ over this interval is at most $\E[\deg(P_{A_\infty})^k]+\varepsilon$ for all sufficiently large $D$.

Continuing inductively, we obtain an increasing sequence $N_1<N_2<\cdots$, with $N_j=N(D_j)$ and $N_0=0$, such that for every $j\ge1$,
$$
\begin{aligned}
\frac1{N_j-N_{j-1}}\sum_{n=N_{j-1}+1}^{N_j}\deg(P_{\chi_n})^k
&\ge\E[\deg(P_{A_\infty})^k]+2\varepsilon,\\
\frac1{N_j-N_{j-1}}\sum_{\substack{N_{j-1}<n\le N_j\\\deg(P_{\chi_n})\le j}}\deg(P_{\chi_n})^k
&\le\E[\deg(P_{A_\infty})^k]+\varepsilon.
\end{aligned}
$$

Now define
$$
B=\bigcup_{j\ge1}\{\chi_n:N_{j-1}<n\le N_j,\ \deg(P_{\chi_n})>j\}.
$$
We first show that $B$ has asymptotically negligible density. By construction, for every fixed nonnegative integer $M$, only finitely many characters in $B$ have degree at most $M$. Consequently, for every fixed $M\ge1$,
$$
\limsup_{D\to\infty}\frac{|B\cap X_i(D)|}{N(D)}
\le\lim_{D\to\infty}\Pr(\deg(P_{\chi_D})>M)
=\Pr(\deg(P_{A_\infty})>M),
$$
where the equality follows from \Cref{conj:distinguished_polynomial_intro}. Letting $M\to\infty$, we conclude that $B$ has asymptotically negligible density.

On the other hand, for every $j\ge1$, our choice of $N_j$ gives
$$
\begin{aligned}
\frac1{N_j-N_{j-1}}\sum_{\substack{N_{j-1}<n\le N_j\\\chi_n\in B}}\deg(P_{\chi_n})^k
&=\frac1{N_j-N_{j-1}}\left(\sum_{n=N_{j-1}+1}^{N_j}\deg(P_{\chi_n})^k
-\sum_{\substack{N_{j-1}<n\le N_j\\\deg(P_{\chi_n})\le j}}\deg(P_{\chi_n})^k\right)\\
&\ge (\E[\deg(P_{A_\infty})^k]+2\varepsilon)-(\E[\deg(P_{A_\infty})^k]+\varepsilon)=\varepsilon.
\end{aligned}
$$
Thus, it follows that
$$
\frac1{N_j}\sum_{\chi\in B\cap X_i(D_j)}\deg(P_\chi)^k
\ge\frac{1}{N_j}(N_1\varepsilon+(N_2-N_1)\varepsilon+\cdots+(N_j-N_{j-1})\varepsilon)
\ge\varepsilon.
$$
This contradicts \Cref{conj:bad_sets_don't_effect_averages_intro}, since $B$ has asymptotically negligible density. Therefore,
$$
\limsup_{D\to\infty}\E[\deg(P_{\chi_D})^k]\le\E[\deg(P_{A_\infty})^k].
$$
Combining this with the lower bound proves
$$
\lim_{D\to\infty}\E[\deg(P_{\chi_D})^k]=\E[\deg(P_{A_\infty})^k].
$$

Finally, let $B\subset X_i$ be any subset of asymptotically negligible density. For all sufficiently large $D$,
$$
\E[\deg(P_{\widetilde\chi_D})^k]
=\frac{\E[\deg(P_{\chi_D})^k]-\frac1{N(D)}\sum_{\chi\in B\cap X_i(D)}\deg(P_\chi)^k}
{1-\frac{|B\cap X_i(D)|}{N(D)}}.
$$ 
By the assumption in \Cref{conj:bad_sets_don't_effect_averages_intro}, the sum in the numerator tends to zero, while the denominator tends to $1$. Hence
$$
\lim_{D\to\infty}\E[\deg(P_{\widetilde\chi_D})^k]
=\E[\deg(P_{A_\infty})^k],
$$
which completes the proof.
\end{proof}

\begin{rmk}[Moment tightness and negligible-density sets]\label{rmk: moment tightness and negligible density sets}
Under \Cref{conj:distinguished_polynomial_intro}, the negligible-density condition in \Cref{conj:bad_sets_don't_effect_averages_intro} is equivalent to the moment-tightness condition used in \cite[Theorem 3.3]{shen2026resultant}, imposed for every fixed integer $k\ge1$. In our notation, this condition reads
$$
\lim_{M\to\infty}\limsup_{D\to\infty}\E[\deg(P_{\chi_D})^k\mathbf{1}_{\{\deg(P_{\chi_D})>M\}}]=0.
$$
To see this, \Cref{thm:moment-convergence-from-conjectures}, together with the weak convergence in \Cref{thm:conj-implies-disc-lambda}, implies that for every fixed $M$ the expectation inside the limit converges to $\E[\deg(P_{A_\infty})^k\mathbf{1}_{\{\deg(P_{A_\infty})>M\}}]$, which tends to zero as $M\to\infty$. Conversely, if moment tightness holds and $B\subset X_i$ has asymptotically negligible density, then for every $M\ge1$,
$$
\frac{1}{|X_i(D)|}\sum_{\chi\in B\cap X_i(D)}\deg(P_\chi)^k
\le M^k\frac{|B\cap X_i(D)|}{|X_i(D)|}+\E[\deg(P_{\chi_D})^k\mathbf{1}_{\{\deg(P_{\chi_D})>M\}}].
$$
Taking the limit superior as $D\to\infty$ and then letting $M\to\infty$ gives \Cref{conj:bad_sets_don't_effect_averages_intro}. 
\end{rmk}

We finally return to local zero statistics. \Cref{thm:local-zero-moment-transfer} below shows a general version of the convergence statement needed, which can be specialized to each case. 

\begin{prop}[Convergence of local zero moments]
\label{thm:local-zero-moment-transfer}
Assume \Cref{conj:distinguished_polynomial_intro} and \Cref{conj:bad_sets_don't_effect_averages_intro}, and fix a twist $i\in\{0,\ldots,p-2\}$. For $K_1,\ldots,K_m$ and $U$ as in \Cref{thm:conj-implies-local-zero-distribution}, and every integer $k\ge1$,
$$
\lim_{D\to\infty}\E[Z_U(P_{\chi_D})^k]=\E[Z_U(P_{A_\infty})^k]=\lim_{n\to\infty}\E[Z_U(P_{A_n})^k].
$$
\end{prop}

\begin{proof}
The idea is that since $0\le Z_U(P)\le\deg(P)^m$, the convergence of degree moments in \Cref{thm:moment-convergence-from-conjectures} controls the contribution from polynomials of large degree. 

Fix an integer $M\ge 1$. Since $0\le Z_U(P)\le\deg(P)^m$, the function
$$
P\mapsto Z_U(P)^k\mathbf{1}_{\{\deg(P)\le M\}}
$$
is bounded by $M^{mk}$. Moreover, by \Cref{lem:root-count-continuity}, it is continuous at every squarefree polynomial, since the degree is locally constant on $\mathscr P$. By \Cref{thm:existence of limiting distinguished polynomial}, the limiting polynomial $P_{A_\infty}$ is almost surely squarefree. Thus \Cref{conj:distinguished_polynomial_intro} implies that
$$
\lim_{D\to\infty}\E[Z_U(P_{\chi_D})^k\mathbf{1}_{\{\deg(P_{\chi_D})\le M\}}]
=\E[Z_U(P_{A_\infty})^k\mathbf{1}_{\{\deg(P_{A_\infty})\le M\}}].
$$

We next control the contribution from polynomials of degree greater than $M$. For any $P\in\mathscr P$ with $\deg(P)>M$, we have
$$
Z_U(P)^k\le\deg(P)^{mk}\le\frac{\deg(P)^{mk+1}}{M},
$$
where the second inequality follows from $\deg(P)/M>1$. Consequently,
$$
\E[Z_U(P_{\chi_D})^k\mathbf{1}_{\{\deg(P_{\chi_D})>M\}}]
\le\frac{1}{M}\E[\deg(P_{\chi_D})^{mk+1}].
$$
By \Cref{thm:moment-convergence-from-conjectures}, the expectation $\E[\deg(P_{\chi_D})^{mk+1}]$ converges to the finite quantity $\E[\deg(P_{A_\infty})^{mk+1}]$. The same tail estimate applies to $P_{A_\infty}$. Thus the contributions from degrees greater than $M$ are bounded by a constant times $M^{-1}$, uniformly for all sufficiently large $D$ and for the limiting polynomial. Combining this estimate with the convergence of truncated expectations and letting $M\to\infty$ proves the first equality.

For the random matrix model, the same argument gives convergence of the truncated expectations, since $P_{A_n}$ converges weakly to $P_{A_\infty}$. To control the remaining contribution, \Cref{lem:random-matrix-degree-tail} gives
$$
\mathbb P(\deg(P_{A_n})=r)
=p^{-r}\prod_{j=r+1}^n(1-p^{-j})\le p^{-r}
$$
for every $0\le r\le n$. It follows that
$$
\E[\deg(P_{A_n})^{mk+1}]
\le\sum_{r=1}^\infty r^{mk+1}p^{-r}<\infty,
$$
uniformly in $n$, since the series converges by exponential decay. Therefore,
$$
\E[Z_U(P_{A_n})^k\mathbf{1}_{\{\deg(P_{A_n})>M\}}]
\le\frac{1}{M}\sum_{r=1}^\infty r^{mk+1}p^{-r}.
$$
The same bound holds for $P_{A_\infty}$ by \Cref{lem:random-matrix-degree-tail}. Hence the contributions from degrees greater than $M$ tend to zero uniformly in $n$ as $M\to\infty$. Together with the convergence of truncated expectations, this proves the second equality.
\end{proof}

\subsection{Proofs of Consequences~\ref{thm:lambda_intro} through~\ref{cor:cubic_extensions}}

We now apply the preceding convergence results to the explicit random matrix computations of \Cref{sec: Random matrix computations}. 

\begin{proof}[Proof of \Cref{thm:lambda_intro}]
By \Cref{thm:conj-implies-disc-lambda}, the distribution of $\deg(P_{\chi_D})$ converges to that of $\deg(P_{A_\infty})$. Hence, for every $d\ge0$,
$$
\lim_{D\to\infty}\Pr\bigl(\deg(P_{\chi_D})=d\bigr)=\Pr\bigl(\deg(P_{A_\infty})=d\bigr).
$$
By \Cref{lem:random-matrix-degree-tail},
$$
\Pr\bigl(\deg(P_{A_\infty})=d\bigr)
=p^{-d}\prod_{t>d}(1-p^{-t}).
$$
This proves the claim.
\end{proof}

\begin{proof}[Proof of \Cref{thm:disc_ejv_intro}]
By \Cref{thm:conj-implies-disc-lambda}, for every $s\geq1$,
$$
\lim_{D\to\infty}\Pr\bigl(\val(\Disc(P_{\chi_D}))=s,\deg(P_{\chi_D})=2\bigr)=\Pr\bigl(\val(\Disc(P_{A_\infty}))=s,\deg(P_{A_\infty})=2\bigr).
$$
By \Cref{thm:existence of limiting distinguished polynomial GL}, the distribution computed by Ellenberg-Jain-Venkatesh \cite[Section 4.5]{ellenberg2011modeling} for the multiplicative Haar model also applies to $P_{A_\infty}$. Write $\eta=\prod_{j\ge1}(1-p^{-j})$. Their computation gives
$$
\Pr\bigl(\val(\Disc(P_{A_\infty}))=s,\deg(P_{A_\infty})=2\bigr)=\eta\frac{p^{-s-1}}{1-p^{-1}}
\begin{cases}
1-p^{-(s+1)/2}, & s\text{ odd},\\
1-\dfrac{p^{-s/2-1}}{1+p^{-1}}, & s\text{ even}.
\end{cases}
$$
On the other hand, \Cref{lem:random-matrix-degree-tail} gives
$$
\Pr\bigl(\deg(P_{A_\infty})=2\bigr)
=p^{-2}\prod_{t>2}(1-p^{-t})
=\frac{\eta p^{-2}}{(1-p^{-1})(1-p^{-2})}.
$$
Dividing the two expressions yields
$$
\lim_{D\to\infty}\Pr\bigl(\val(\Disc(P_{\chi_D}))=s
\,\big|\,\deg(P_{\chi_D})=2\bigr)=
\begin{cases}
p^{1-s}(1-p^{-2})(1-p^{-(s+1)/2}),
& s\text{ odd},\\
p^{1-s}(1-p^{-2})
\left(1-\dfrac{p^{-s/2-1}}{1+p^{-1}}\right),
& s\text{ even},
\end{cases}
$$
as claimed.
\end{proof}

We turn to the local statistics of zeros. 

\begin{proof}[Proof of \Cref{cor:Zp_zeros_intro}]
Let $B=a+p^r\Z_p$. By \Cref{thm:local-zero-moment-transfer}, the expectations of the corresponding root-counting statistics converge to their random matrix counterparts. By \Cref{thm:1-pt_cor_fn_Zp}, the limiting one-point correlation function on $\Z_p$ is identically $1$. Hence
$$
\lim_{D\to\infty}\E[Z_{B,\chi_D}]=\int_B 1\,dx
=p^{-r},
$$
which proves \eqref{eq:zp_density_intro}.

Now let $B_1=a_1+p^r\Z_p$ and $B_2=a_2+p^r\Z_p$ be disjoint. Since the disks are disjoint, for every polynomial $P$ we have $Z_{B_1}(P)Z_{B_2}(P)=Z_{B_1\times B_2}(P)$. Thus \Cref{thm:local-zero-moment-transfer} gives
$$
\lim_{D\to\infty}\E[Z_{B_1,\chi_D}Z_{B_2,\chi_D}]=\lim_{n\to\infty}\E[Z_{B_1\times B_2,n}].
$$
Write $|a_1-a_2|=p^{-m}$. Since $B_1$ and $B_2$ are disjoint, we have $m<r$, and hence
$|x_1-x_2|=p^{-m}$ for all $(x_1,x_2)\in B_1\times B_2$. Therefore, by \Cref{prop: limiting correlation function} and \Cref{thm: 2pt in Zp},
$$
\lim_{n\to\infty}\E[Z_{B_1\times B_2,n}]=p^{-2r}\rho^{(\infty)}_{\Q_p\times\Q_p}(a_1,a_2).
$$
The series expression in \Cref{thm: 2pt in Zp} may now be rewritten using the Jacobi theta function. Using the series definition of $\theta_3$ given in \Cref{sec: Introduction}, and separating the positive and negative indices, we obtain

$$
\begin{aligned}
\theta_3(-\sqrt p;p^{-2m-1})
&=1+\sum_{k\geq0}(-1)^{k+1}
p^{-m-\frac{(2m+1)k^2+(4m+1)k}{2}}+\sum_{k\geq0}(-1)^{k+1}
p^{-m-1-k-\frac{(2m+1)k^2+(4m+1)k}{2}}\\
&=1-\rho^{(\infty)}_{\Q_p\times\Q_p}(a_1,a_2).
\end{aligned}
$$
Combining this with the first moment given in \eqref{eq:zp_density_intro}, we obtain
$$
\lim_{D\to\infty}\Cov(Z_{B_1,\chi_D},Z_{B_2,\chi_D})=p^{-2r}\left(\rho^{(\infty)}_{\Q_p\times\Q_p}(a_1,a_2)-1\right)
=-p^{-2r}\theta_3\left(-\sqrt p;\frac{|a_1-a_2|^2}{p}\right),
$$
which proves \eqref{eq:covariances_intro}.
\end{proof}

\begin{proof}[Proof of \Cref{cor:Zp_higher_moments_intro}]
Let $B=a+p^r\Z_p$, where $a\in p\Z_p$ and $r\ge 1$. By \Cref{thm:local-zero-moment-transfer}, the first two moments of $Z_{B,\chi_D}$ converge to the corresponding random matrix moments. More precisely, for $k=1,2$, we have
$$
\lim_{D\to\infty}\E[Z_{B,\chi_D}^k]
=\lim_{n\to\infty}\E[Z_{B,n}^k].
$$
For a Haar-distributed matrix $A_n\in\Mat_n(\Z_p)$, the shifted matrix $A_n-aI_n$ has the same distribution as $A_n$, since additive Haar measure on $\Mat_n(\Z_p)$ is translation invariant. Moreover, the eigenvalues of $A_n-aI_n$ are obtained by subtracting $a$ from those of $A_n$. Consequently, $Z_{B,n}$ and $Z_{p^r\Z_p,n}$ have the same distribution. It therefore suffices to compute the limiting first and second moments of $Z_{p^r\Z_p,n}$.

Recall from \Cref{thm: 1pt expectation in plZp} that
$\E[Z_{p^r\Z_p,n}]=p^{-r}$ for every $n\ge 1$. This implies that
$$
\lim_{D\to\infty}\E[Z_{B,\chi_D}]=p^{-r}.
$$
To compute the second moment, recall that $Z_{p^r\Z_p\times p^r\Z_p,n}$ counts ordered pairs of distinct eigenvalues lying in $p^r\Z_p$. Separating the diagonal and off-diagonal contributions gives
$$
Z_{p^r\Z_p,n}^2
=Z_{p^r\Z_p,n}+Z_{p^r\Z_p\times p^r\Z_p,n}.
$$
Taking expectations and applying \Cref{thm: 2pt expectation in plZp}, we obtain
$$
\begin{aligned}
\lim_{n\to\infty}\E[Z_{p^r\Z_p,n}^2]
&=p^{-r}+\lim_{n\to\infty}\E[Z_{p^r\Z_p\times p^r\Z_p,n}]\\
&=p^{-r}+(1-p^{-1})\sum_{k\ge 0}(-1)^k(1+p^{-k-1})
\frac{p^{-\frac{k^2+k}{2}}p^{-r(k^2+2k+3)}}{1-p^{-(k^2+2k+2)}}.
\end{aligned}
$$
Therefore, we conclude that
$$
\begin{aligned}
\lim_{D\to\infty}\E\left[\left(Z_{B,\chi_D}-\E[Z_{B,\chi_D}]\right)^2\right]&=\lim_{n\to\infty}\E[Z_{p^r\Z_p,n}^2]-p^{-2r}\\
&=p^{-r}-p^{-2r}+(1-p^{-1})\sum_{k\ge 0}(-1)^k(1+p^{-k-1})
\frac{p^{-\frac{k^2+k}{2}}p^{-r(k^2+2k+3)}}{1-p^{-(k^2+2k+2)}},
\end{aligned}
$$
which is precisely the expression in \eqref{eq:zp_variance_intro}.
\end{proof}

We next treat zeros in quadratic extensions. 

\begin{proof}[Proof of \Cref{cor:quadratic_intro}]
By \Cref{thm:local-zero-moment-transfer}, the expected number of zeros in $B_K$ converges to the corresponding random matrix expectation:
$$
\lim_{D\to\infty}\E[Z_{B_K,\chi_D}]=\lim_{n\to\infty}\E[Z_{B_K,n}].
$$

Suppose first that $K/\Q_p$ is unramified. Since the valuation group of $K$ is $\Z$, we have $B_K=p\mathcal O_K^{\new}$. Hence, applying the unramified formula in \Cref{thm: expected number in quadratic extensions} by taking $l=1$, we have
$$
\begin{aligned}
\lim_{D\to\infty}\E[Z_{B_K,\chi_D}]
&=\sum_{k\geq0}(1-p^{-k-1})p^{-\frac{k^2+k}{2}}\\
&\quad\times
\frac{(1+p^{-1})p^{-(k^2+2k+3)}
(1-p^{-(k^2+2k+3)})-2p^{-1}p^{-(k^2+2k+4)}
(1-p^{-(k^2+2k+2)})}{(1-p^{-(k^2+2k+2)})
(1-p^{-(k^2+2k+3)})}\\
&=p^{-3}-2p^{-5}+2p^{-6}+p^{-7}+O(p^{-8}),
\end{aligned}
$$
which is the first assertion. Now suppose that $K/\Q_p$ is ramified. Let $x_1$ be as in \Cref{thm: quadratic correlation function}, so that $|x_1|=p^{-1/2}$. Then
$B_K=x_1\mathcal O_K\cap\mathcal O_K^{\new}$. Applying the second ramified formula in \Cref{thm: expected number in quadratic extensions} by taking $l=0$, and using $|\Disc_{K/\Q_p}|=p^{-1}$ gives
$$
\begin{aligned}
\lim_{D\to\infty}\E[Z_{B_K,\chi_D}]
&=p^{-2}\sum_{k\geq0}\frac{(1-p^{-2k-2})p^{-(k^2+k)}}{(1-p^{-(k^2+2k+2)})(1-p^{-(k^2+2k+3)})}\\
&\qquad\times\left[(1-p^{-(k^2+2k+3)})-p^{-1}(1-p^{-(k^2+2k+2)})\right]\\
&=p^{-2}(1-p^{-1})\sum_{k\geq0}\frac{(1-p^{-2k-2})p^{-(k^2+k)}}{(1-p^{-(k^2+2k+2)})(1-p^{-(k^2+2k+3)})}\\
&=p^{-2}-p^{-3}+p^{-4}-p^{-6}+O(p^{-8}),
\end{aligned}
$$
as claimed.
\end{proof}

For totally ramified extensions of general degree, we consider the first annulus, where the sum over all extensions has a particularly simple expression.

\begin{proof}[Proof of \Cref{cor:ramified_shell_intro}]
Fix a totally ramified extension $K/\Q_p$ of degree $m$ and a uniformizer $\pi\in\mathcal O_K$. Every element of $\pi\mathcal O_K\setminus\pi^2\mathcal O_K$ is again a uniformizer of $K$, and hence generates $K/\Q_p$. 
By \Cref{thm:local-zero-moment-transfer},
$$
\lim_{D\to\infty}
\E[Z_{\pi\mathcal O_K\setminus\pi^2\mathcal O_K,\chi_D}]
=\lim_{n\to\infty}\E[Z_{\pi\mathcal O_K\setminus\pi^2\mathcal O_K,n}].
$$
By the finiteness theorem for extensions of local fields \cite[Proposition 7.64]{milne2020algebraic}, there are only finitely many totally ramified extensions of $\Q_p$ of degree $m$. Summing over $K$ and applying \Cref{thm: summing slope 1/m} gives
$$
\lim_{D\to\infty}\sum_K\E[Z_{\pi\mathcal O_K\setminus\pi^2\mathcal O_K,\chi_D}]
=\lim_{n\to\infty}\sum_K\E[Z_{\pi\mathcal O_K\setminus\pi^2\mathcal O_K,n}]
=\frac{m(p-1)}{p^{m+1}},
$$
as claimed.
\end{proof}

We finish with two further consequences of the random matrix computations. The first concerns the centered third moment of the number of zeros in a disk. 

\begin{consequence}[The centered third moment in $p^r\Z_p$]
\label{cor:Zp_third_moment}
With the same assumptions and setting as \Cref{cor:Zp_zeros_intro}, let $B=a+p^r\Z_p\subset p\Z_p$ with $r\ge1$, and let $V_r$ denote the right-hand side of \eqref{eq:zp_variance_intro}. Then
\begin{multline*}
\lim_{D\to\infty}\E\left[\left(Z_{B,\chi_D}-\E[Z_{B,\chi_D}]\right)^3\right]-(3(1-p^{-r})V_r-2p^{-r}+3p^{-2r}-p^{-3r})\\
\in\left[
\frac{p^{-6r}(1-p^{-1}+p^{-2})(p^{-1};p^{-1})_\infty}{(1-p^{-1})(1-p^{-2})(1-p^{-5})},
\frac{p^{-6r}(1-p^{-1}+p^{-2})(p^{-1};p^{-1})_3}{(1-p^{-1})(1-p^{-2})(1-p^{-5})}
\right].
\end{multline*}
\end{consequence}

\begin{proof}
By \Cref{thm:local-zero-moment-transfer}, the first three moments of $Z_{B,\chi_D}$ converge to the corresponding random matrix moments. By translation invariance of additive Haar measure on $\Mat_n(\Z_p)$, $Z_{B,n}$ and $Z_{p^r\Z_p,n}$ have the same distribution. Hence it suffices to compute the limiting centered third moment of $Z_{p^r\Z_p,n}$.

Since $Z_{p^r\Z_p\times p^r\Z_p,n}$ and $Z_{p^r\Z_p\times p^r\Z_p\times p^r\Z_p,n}$ count ordered tuples of distinct eigenvalues, we have
$$
Z_{p^r\Z_p,n}^2=Z_{p^r\Z_p\times p^r\Z_p,n}+Z_{p^r\Z_p,n},
\qquad
Z_{p^r\Z_p,n}^3=Z_{p^r\Z_p\times p^r\Z_p\times p^r\Z_p,n}+3Z_{p^r\Z_p\times p^r\Z_p,n}+Z_{p^r\Z_p,n}.
$$
By \Cref{thm: 1pt expectation in plZp}, $\E[Z_{p^r\Z_p,n}]=p^{-r}$, while \Cref{cor:Zp_higher_moments_intro} gives
$$
\lim_{n\to\infty}\E[Z_{p^r\Z_p\times p^r\Z_p,n}]=V_r-p^{-r}+p^{-2r}.
$$
Therefore,
$$
\lim_{n\to\infty}\E\left[\left(Z_{p^r\Z_p,n}-p^{-r}\right)^3\right]=\lim_{n\to\infty}\E[Z_{p^r\Z_p\times p^r\Z_p\times p^r\Z_p,n}]+3(1-p^{-r})V_r-2p^{-r}+3p^{-2r}-p^{-3r}.
$$
Finally, applying \Cref{thm: 3pt zeros in pZp}, we have
$$
\lim_{n\to\infty}\E[Z_{p^r\Z_p\times p^r\Z_p\times p^r\Z_p,n}]
\in\left[\frac{p^{-6r}(1-p^{-1}+p^{-2})(p^{-1};p^{-1})_\infty}{(1-p^{-1})(1-p^{-2})(1-p^{-5})},
\frac{p^{-6r}(1-p^{-1}+p^{-2})(p^{-1};p^{-1})_3}{(1-p^{-1})(1-p^{-2})(1-p^{-5})}\right],
$$
which completes the proof.
\end{proof}

The second consequence concerns zeros in cubic extensions. 

\begin{consequence}[Zeros in cubic extensions]
\label{cor:cubic_extensions}
With the same assumptions and setting as \Cref{cor:Zp_zeros_intro}, let $K/\Q_p$ first be the unique unramified cubic extension. Then $\lim_{D\to\infty}\E[Z_{p\mathcal O_K^{\new},\chi_D}]$ has lower and upper bounds given by
\begin{multline}
\sum_{m\ge1}\frac{p^{-5m-1}}{(1-p^{-1})(1-p^{-3})}
\left(\frac{(p+1)(p^2+p+1)}{p^3}-\frac{3(p^2+p+1)}{p^{m+3}}+\frac{3}{p^{3m+3}}\right)\\
\times\left((p^{-1};p^{-1})_3\pm(1-(p^{-1};p^{-1})_\infty)p^{-3m}\right).
\end{multline}
For the ramified cubic extensions $K/\Q_p$, let $\pi$ be a uniformizer of $K$. Then $\lim_{D\to\infty}\sum_K\E[Z_{\pi\mathcal O_K\cap\mathcal O_K^{\new},\chi_D}]$ has lower and upper bounds given by
$$
\begin{aligned}
\frac{3(p-1)}{p^4}
&+\sum_{m\ge1}\frac{3p^{-5m-3}}{(1-p^{-1})(1-p^{-2})}
\left(\frac{p+1}{p}-\frac{2p+1}{p^{m+2}}+\frac{1}{p^{3m+3}}\right)
\left((p^{-1};p^{-1})_3\pm\frac{1-(p^{-1};p^{-1})_\infty}{p^{3m+1}}\right)\\
&+\sum_{m\ge0}\frac{3p^{-5m-5}}{(1-p^{-1})(1-p^{-2})}
\left(\frac{p+1}{p}-\frac{p+2}{p^{m+2}}+\frac{1}{p^{3m+4}}\right)
\left((p^{-1};p^{-1})_3\pm\frac{1-(p^{-1};p^{-1})_\infty}{p^{3m+2}}\right),
\end{aligned}
$$
where the sum is over all ramified cubic extensions $K/\Q_p$, and taking the minus sign gives the lower bound while taking the plus sign gives the upper bound.
\end{consequence}

\begin{proof}
By \Cref{thm:local-zero-moment-transfer}, the expected number of zeros in $p\mathcal O_K^{\new}$ or $\pi\mathcal O_K\cap\mathcal O_K^{\new}$ converges to the corresponding random matrix expectation.

Suppose first that $K/\Q_p$ is unramified. The first claim follows immediately from \Cref{thm:cubic_small_ball_expectation} by taking $l=1$ in \eqref{eq:cubic_unramified_small_ball}.

Now let $K/\Q_p$ be ramified cubic and let $\pi$ be a uniformizer. Since a cubic extension has prime degree, every element of $K\setminus\Q_p$ generates $K$. In particular, every element of $\pi\mathcal O_K\setminus\pi^2\mathcal O_K$ generates $K/\Q_p$, and hence
$$
\pi\mathcal O_K\cap\mathcal O_K^{\new}
=\left(\pi\mathcal O_K\setminus\pi^2\mathcal O_K\right)\bigsqcup\left(\pi^2\mathcal O_K\cap\mathcal O_K^{\new}\right).
$$
Since there are only finitely many ramified cubic extensions of $\Q_p$, by \Cref{thm:local-zero-moment-transfer} we have
$$
\lim_{D\to\infty}\sum_K\E[Z_{\pi\mathcal O_K\cap\mathcal O_K^{\new},\chi_D}]
=\lim_{n\to\infty}\sum_K\E[Z_{\pi\mathcal O_K\setminus\pi^2\mathcal O_K,n}]
+\lim_{n\to\infty}\sum_K\E[Z_{\pi^2\mathcal O_K\cap\mathcal O_K^{\new},n}].
$$
For the first contribution, \Cref{thm: summing slope 1/m} with $m=3$ gives
$$
\lim_{n\to\infty}\sum_K\E[Z_{\pi\mathcal O_K\setminus\pi^2\mathcal O_K,n}]
=\frac{3(p-1)}{p^4}.
$$
For the second contribution, \Cref{thm:cubic_small_ball_expectation} with $l=1$ in \eqref{eq:cubic_ramified_3lminus1_small_ball} gives the two bounds for $\pi^2\mathcal O_K\cap\mathcal O_K^{\new}$. Every term is proportional to $|\Disc_{K/\Q_p}|$. Comparing the individual and summed formulas in \Cref{thm: totally ramified 1/m} and \Cref{thm: summing slope 1/m} gives
$$
\sum_K|\Disc_{K/\Q_p}|=\frac{3}{p^2}.
$$
Summing the bounds in \eqref{eq:cubic_ramified_3lminus1_small_ball} over $K$ and using this identity gives the two series displayed above. Adding the contribution from $\pi\mathcal O_K\setminus\pi^2\mathcal O_K$ proves the ramified claim.
\end{proof}

\section{Computational methods}\label{sec:computational_methods}

Our computation proceeds as follows:
\begin{enumerate}
    \item Compute special values of the classical Dirichlet $L$-function $L(s,\chi)$ at nonpositive integers.
    \item Use these to obtain an approximation of the Iwasawa power series via Lagrange interpolation.
    \item Compute the Newton polygon of the approximate Iwasawa power series. For sufficiently many interpolation points, this is provably the same as the Newton polygon of the Iwasawa power series and its distinguished polynomial.
    \item Compute additional information on zeros not captured by the Newton polygon (which only gives their $p$-adic valuations), namely which subextensions $K \subset \bar{\Q}_p$ and which sub-disks of the unit disk in $K$ they lie in.
\end{enumerate}

In this section, we discuss these steps in finer detail in order. We remind the reader that the code itself may be accessed at \url{https://github.com/knospe/iwasawa-power-series}.

\subsection{Special values of $L$-functions}
\label{sec:special-values}
Let $\theta$ be a quadratic character with fundamental discriminant $\mf{D}$ and conductor $\mathcal{D}=|\mf{D}|$, and let $\omega^i$ be a Teichm\"uller twist where $i=0,\dots,p-2$.  We often denote $\chi = \theta \omega^i$ as before. We approximate the Iwasawa power series $F_{\chi}$ 
via interpolation at special values. To this end, we use formulas \eqref{eq:powerseries} and \eqref{eq:interpol_general} and compute
special values of the $\chi$-twisted $L$-function at negative integers.  The straightforward way to obtain the special values is via twisted Bernoulli numbers, which requires time $O(\mathcal{D})$, done for instance in \cite{delbourgo2015invariants,delbourgo2023iwasawa,knospe2024special}. This can conveniently be employed for conductors up to around $10^6$, but we observed that the empirical averages were still not matching random matrix predictions in this range, and hence used two other methods for the main data presented here. 

One method, essentially that of \cite{ellenberg2011modeling}, is based on computations with analytic methods running in time $O(\mathcal{D}^{1/2 + \epsilon})$. We used the {\ttfamily PARI} software library for these computations. For the large conductors $\sim 10^{10}$ which we tested, these methods were the most efficient, but additional speedups optimized to a few specific statistics, detailed in \Cref{subsec:high_conductor_speedups}, were still required to compute sufficient data.

For our main dataset, consisting of characters with conductors up to approximately $10^7$, we found that the most efficient method was to compute special values using Cohen-Eisenstein series, defined e.g. in \cite[Definition 10.3.11]{cohen2007number}. These are certain series $\sum_{n\geq 0} H_k(n)q^n$ defined for integers $k \geq 2$, where for those $n$ such that $n=|\mf{D}|$ for a fundamental discriminant $\mf{D}$ of sign $(-1)^k$, the coefficient is given by 
\begin{equation}
    H_k(n) = L(1-k,\theta_{\mf{D}})
\end{equation}
where $\theta_{\mf{D}}$ is the unique quadratic character with discriminant $\mf{D}$ (the definition of the other coefficients is more involved, but is explained in the same reference). This series is a modular form of weight $k + \frac{1}{2}$ on the congruence subgroup $\Gamma_0(4)$ \cite[Theorem 10.3.12]{cohen2007number}, which allows {\ttfamily PARI} to compute the coefficients efficiently. 

We used {\ttfamily PARI} to compute the first $10^7$ Fourier coefficients of the Cohen-Eisenstein series for weights $k+\frac{1}{2}$, where $2 \leq k \leq 105$. Each series simultaneously provides the special values $L(1-k,\theta)$ for all quadratic Dirichlet characters $\theta$ of conductor $<10^7$.

In a second step, we extended this to $2 \times 10^7$ Fourier coefficients, but this required massive computing resources, including several hundred gigabytes of memory and several days of running time on a large server. The main advantage of this approach is that computations of the power series $F_{\chi}$ are rather fast once the required Cohen-Eisenstein series are computed. Using the Fourier coefficients of $\sim 100$ series, we could calculate the Iwasawa power series for all characters $\chi=\theta \omega^i$ in the families $X_i$, for $p=3,5,7$, with primitive conductor of $\theta$ less than $10^7$ (and even $<2 \times 10^7$), which had not been done before. 

\subsection{Approximating Iwasawa power series via interpolation} We leverage the following observation (which follows immediately from equation \eqref{eq:interpol_general}) to handle the Teichm\"uller twists $\theta \omega^i$ of a quadratic character $\theta$. We avoid non-trivial Teichm\"uller twists of the complex $L$-functions since their values are harder to compute. Note that $\theta \omega^i$ is no longer a quadratic character, unless $i=0$ or $i=\frac{p-1}{2}$.

\begin{prop}
Let $i \in \Z$ and let $\theta \omega^i$ be a twisted Dirichlet character. Suppose that $k = i + m(p-1) \geq 1$ with $m \in \Z$. Then
$$ L_p(1-k,\theta \omega^i) = (1- \theta(p) p^{k-1}) L(1-k,\theta) . $$
\label{prop:special-values}
\end{prop}
In other words, by choosing $k$ appropriately, we make it so that the special value on the right-hand side involves only a quadratic character rather than a twist of one. For a fixed Teichm\"uller twist $\omega^i$, we use Proposition \ref{prop:special-values} to compute around $20$ special values $L_p(1-k,\theta \omega^i)$, where $k=i+(p-1)$, $i+2(p-1),\dots$. This gives special values of the power series $F_{\theta \omega^i}(T)$ at $T=(1+dp)^{1-k}-1$, where $d$ is the conductor of $\theta$. Note that $v_p(T)\geq 1$ for $T=(1+dp)^{1-k}-1$. Now we use Lagrange interpolation to obtain a polynomial $L(T)$ with  the same function values at the chosen set of points. What is the relation between $F_{\theta \omega^i}(T)$ and $L(T)$? The following proposition controls the error between their coefficients. As expected, the number of points used for interpolation increases the accuracy.  

\begin{prop} Let $f(T)=\sum a_j T^j \in \Z_p[[T]]$ and let $x_1, \dots, x_n \in p\Z_p$ be distinct points. Let 
$$W(T)=\prod_{j=1}^n (T-x_j),$$
and let
$$L(T)= \sum b_j T^j \in \Z_p[T]$$ 
be the polynomial of degree $<n$ that interpolates $f(T)$ at the points $(x_j,f(x_j))$ for $j=1,\dots,n$. Then 
$$ f(T)=q(T) W(T) + L(T),$$
for some $q(T) \in \Z_p[[T]]$. Furthermore, we have
$$a_j \equiv b_j \mod p^{n-j}$$ 
for $j=0,1,\dots, n-1$. In particular,
$$ v_p(a_j) = v_p(b_j) $$
if $v_p(b_j) < n-j  $.

\label{prop:interpolation-polynomial}
\end{prop}
\begin{proof}
Since we assumed that the points $x_i$ are distinct, the interpolation polynomial $L(T)$ of degree $<n$ is unique.
Since all $x_i$ have positive valuation, the polynomial $W(T)$ is distinguished. Hence there exists a unique polynomial $\tilde{L}(T)$ of degree $<n$ such that%
$$ f(T)=q(T) W(T) + \tilde{L}(T)$$
(this is standard and typically appears as an intermediate step in the $p$-adic Weierstrass preparation theorem). By construction, $W$ vanishes at the points $x_j$ and hence $f(x_j)=\tilde{L}(x_j)$ for all $j=1,\dots,n$. The uniqueness of interpolation then implies $L=\tilde{L}$.
Expanding the polynomial $W(T)$ gives
$$ W(T) = \sum_{j=0}^n c_j T^j ,$$
where the coefficients $c_j$ are given by evaluating the elementary symmetric polynomial of degree $n-j$ in $n$ variables  at our points $x_1, \dots, x_n$ (up to sign). Since $v_p(x_j)\geq 1$ we have $v_p(c_j) \geq n-j $. Let $q(T)=\sum_{j=0}^{\infty} q_j T^j \in \Z_p[[T]]$. Matching the coefficients of $f(T)-L(T)$ and $q(T)W(T)$ then yields
$$ a_j - b_j = q_0 c_j + q_1 c_{j-1} + \dots + q_j c_0 .$$
This implies $v_p(a_j - b_j) \geq n-j$, which completes the proof.
\end{proof}

The Newton polygon is the lower convex hull of the points $(j,v_p(a_j))$. We list its slopes in nondecreasing order, each repeated according to the horizontal length of the corresponding segment. The basic theorem of Newton polygons asserts that the roots of such a series in the open unit ball in $\bar{\Q}_p$ are in bijection with the \emph{negative} slopes $s_j$, and the valuation of the root is the (negative of) the corresponding slope. See, for instance, \cite[Chapter 6]{cassels1986local}. We refer to the negatives of these $s_j$, which are positive, as the \emph{principal slopes} of the power series or the Newton polygon. Let $m_j \geq 1$ be the multiplicity of the slope $s_j$. Then the polynomial (with coefficients in $\Z_p$) factorizes over $\Z_p$ into polynomials of degree $m_j$. This follows from the fact that conjugate roots have the same valuation. 

Since $F_{\theta \omega^i}$ and $P_{\theta \omega^i}$ have the same roots in the unit disk, they have the same principal slopes. Therefore, it suffices to compute the Newton polygon of $F_{\theta \omega^i}$. Since the Iwasawa $\mu$-invariant is zero by the Ferrero-Washington theorem, there exists a smallest index $\lambda \geq 0$ such that the coefficient $a_{\lambda}$ of $F_{\theta \omega^i}$ is a $p$-adic unit. This $\lambda$-invariant is clearly also the degree of the associated distinguished polynomial $P_{\theta \omega^i}$, and the slopes $s_j$ are negative exactly when $j=1,\ldots,\lambda$. 

Hence to get the valuations of the roots of $P_{\theta \omega^i}$, we need only compute the Newton polygon of $F_{\theta \omega^i}$, which may be read off from its coefficients. Of course, we only know an approximation to $F_{\theta \omega^i}$ gotten by interpolation of a few special values, but the following Corollary to 
Proposition \ref{prop:interpolation-polynomial} shows that we recover the correct Newton polygon provided that the number of points is greater than $\lambda$ and the valuation of the coefficients is not too large. 

\begin{cor} Let $f(T)=\sum a_j T^j \in \Z_p[[T]]$, where $v_p(a_{\lambda})=0$ and $v_p(a_j)>0$ for $j=0,\dots,\lambda-1$, and let  $L(T)= \sum b_j T^j \in \Z_p[T]$ be defined through interpolation at $n$ points as in Proposition \ref{prop:interpolation-polynomial}. Assume that $n>\lambda$ and $v_p(b_j)<n-j$ for every $0\leq j\leq\lambda$. Then the Newton polygons of $f(T)$ and $L(T)$ coincide on the segment $[0,\lambda]$.
\end{cor}

Since Iwasawa power series with $\lambda>10$ are very scarce (and we have not found any quadratic Dirichlet character of conductor $<2 \times 10^7$ with $\lambda>14$), we could limit the number of points used for interpolation.
We interpolated $F_{\theta \omega^i}$ at  $k=i+m(p-1)$ where $m\geq 1$ and $k<15p$, e.g., $21$ points for $p=3$ and twist $i=1$, and $18$ points for $p=5$ and twist $i=1$. With these choices, we have not found a single example where the above Newton-polygon condition $v_p(b_j)<n-j$ for $0\leq j\leq \lambda$ does not hold.

Another potential source of error which must be accounted for is that $p$-adic numbers are represented by finite expansions in powers of $p$. We used {\ttfamily SageMath}'s representation of $p$-adic numbers with a precision of $O(p^{30})$. However, computations with $p$-adic numbers may reduce their accuracy. This applies in particular to the Lagrange interpolation, which degrades the precision by several digits. The loss of precision of the coefficient $b_j$ of the interpolation polynomial $L(T)=\sum b_j T^j$ increases with $j$. We observed that the loss of precision is only relevant if $\lambda$ is large ($\geq 10$) or the $p$-adic valuation of the discriminant associated with a slope of the Newton polygon is large ($\geq 8$).
In this case, we increased the precision to $O(p^{50})$. We are therefore fairly confident that the Newton polygon and the associated slopes we compute have not been affected by such errors.

From the slopes of the Newton polygon, we immediately obtain the $\lambda$-invariant and the valuations of all roots of $P_{\theta \omega^i}$. However, we do \emph{not} obtain some natural information we would like, such as which extensions the roots lie in, or where within the unit disk they lie.

\begin{example} Let $p=3$. Consider the following examples of the Iwasawa power series $F_{\theta}$ associated with an even quadratic character $\theta$ of fundamental discriminant $\mf{D}$ and twist $i=0$. 
    \begin{enumerate}
     
        \item Let $\mf{D}=3989$. Then interpolation yields
        $$F_{\theta}(T) =  (2\cdot 3^3 + O(3^4))+ (3 + O(3^4))T + (1 + 3 + 2\cdot 3^2 + 2\cdot 3^3 + O(3^4))T^2 +O(T^3) .$$
        Because the first coefficient with valuation $0$ is in degree $2$, we have $\lambda=2$, and a quick computation gives principal slopes are $(2,1)$. Because there are two principal slopes, the roots in the unit disk must lie in an extension of degree at most $2$. Because they have different valuations, they cannot be conjugate roots in a quadratic extension, so they must both lie in $\Z_p$.
        \item Let $\mf{D}=77$. Then
        $$F_{\theta}(T) =  (3 + 3^2  + O(3^4))+ (3 + O(3^4))T+ (2 + 3 + 3^2 + 2\cdot 3^3 + O(3^4))T^2 + O(T^3)  .$$
        Again $\lambda=2$, but now the principal slopes are $(\tfrac{1}{2},\tfrac{1}{2})$. Then both roots must lie in a ramified quadratic extension.
        \item Consider $\mf{D}=2908$, which also appears in the computations of \cite[Supplement, Table II]{ervanmet1992}. Then 
$$ F_{\theta}(T) =  (3^2 + 3^3 + O(3^4))+ (3^4 + O(3^5) )T +(1 + 3^3 + O(3^4))T^2 + O(T^3). $$ 
The principal slopes are $(1,1)$. Then the two roots of $P_{\theta}$ may either both lie in $\Z_p$ or be conjugate elements of a quadratic extension (either ramified or unramified), and there is no way to tell from the Newton polygon alone which is the case. We extract the associated distinguished polynomial of slope $1$, which also requires the coefficients of $ F_{\theta}(T)$ of degree $>2$ for better precision. This gives
$$ P_{\theta}(T) = T^2 + (2 \cdot 3^2 + 3^4 + O(3^5) ) T + (3^2 + 3^3 + 3^4 + O(3^5)) .$$
The discriminant of $ P_{\theta}(T)$ is $2 \cdot 3^2 + 3^4 + O(3^5)$. Since the valuation of the discriminant is even, the extension field must be either trivial or unramified. Finally, we conclude that the extension is unramified as the leading coefficient of the $3$-adic expansion is not a square modulo $3$. 
    \end{enumerate}
    \label{example:ramification}
\end{example}

For many Iwasawa series, such as the first two examples, the Newton polygon alone is enough to tell the extension type of the roots. For others, such as the last case, other techniques are needed, as we describe next. 

\subsection{Additional information on zeros} 

If integer slopes of multiplicity $m \geq 2$ are present, then the slopes do not provide sufficient information to determine the type of the extension field (trivial, unramified or ramified). For these slopes, we computed the corresponding distinguished polynomial factor of degree $m$ from the interpolation polynomial $L(T)$. For $m=2$, the discriminant of that quadratic polynomial provides the required information. Indeed, when $p$ is odd, the extension is ramified if and only if the valuation of the discriminant is odd. Otherwise, the roots are in $\Z_p$ if and only if the leading term of the $p$-adic expansion of the discriminant is a square modulo $p$. 

The cases with integer slopes of multiplicity $m \geq 3$ are more subtle, and the computational errors are more difficult to control.  We first compute the corresponding distinguished polynomial of degree $m$ and try to factorize the polynomial over $\Z_p$. Then we analyze the field extension associated with each irreducible factor of degree $\geq 2$. The ramification of quadratic factors depends only on whether the valuation of the discriminant is even or odd. For an irreducible factor $f$ of degree $m>2$, we try using {\ttfamily SageMath} to construct the field extension of $\Q_p$ and to determine the ramification index. However, this only works for a subset of unramified extensions (and for Eisenstein extensions, which does not apply to our case with integer slopes). We rescale $f$ from slope $s$ to slope $0$; 
the resulting polynomial of slope $0$ is given by
$$\tilde{f} (T) = p^{-ms} f(p^s T) .$$
Since $s \in \Z$, the polynomials $f$ and $\tilde{f}$ generate the same field extension $K/\Q_p$. 

When {\ttfamily SageMath} does not directly determine the ramification (which is in fact the case for some unramified and all ramified extensions), we instead use the following arithmetic tests. 
\begin{enumerate}
\item If the discriminant of $\tilde{f}$ has odd valuation, then the extension must be ramified. 
\item If the  discriminant of $\tilde{f}$ is a unit, then the extension must be unramified. This follows from the discriminant-index formula:
$$ v_p(\text{disc}(\tilde{f})) = v_p(\text{disc}(K)) + 2 v_p(\mathcal{O}_K : (\Z_p[T]/(\tilde{f}))).$$
When the left-hand side of the above equation is zero, the discriminant of $K$ is a unit and $K$ is unramified. In this case, $\tilde{f}$ defines the maximal order $\mathcal{O}_K$.
\item If the reduction of $\tilde{f}$ modulo $p$ is irreducible, then the extension is unramified.
\end{enumerate}

There are a few remaining ambiguous cases (mostly when $p=3$), where  $v_p(\text{disc}(\tilde{f}))$ is positive and even. In most of the remaining cases, we computed the index using the Montes algorithm and then determined the ramification from the discriminant-index formula. A few were not resolved, leading to the small error bars mentioned at the end of \Cref{subsec:numerics_intro}.

The result can be numerically unstable if the discriminant of $f$ has a large $p$-adic valuation, i.e., if the zeros are close to each other. In this case, the precision should be increased. 

\begin{example} Let $p=3$. Let $\theta$ be the quadratic character of fundamental discriminant $\mf{D}=-709255$, which is odd, and consider twist $i=1$. We obtain
\begin{multline*}
     F_{\theta \omega}(T)= 2\cdot3^4 + 2\cdot3^5 + O(3^6) + (2\cdot 3^3 + 2 \cdot 3^4 + 3^5 + O(3^6))T +  O(3^6) T^2  \\ + (3 + 3^2 + 3^4  + O(3^6))T^3 + (2 + 2\cdot 3 + 3^3 + 2 \cdot 3^4 + O(3^6))T^4 + O(T^5)
\end{multline*}
We have $\lambda=4$ and the principal slopes are $(1,1,1,1)$. The associated distinguished polynomial $f(T)=P_{\theta \omega}(T)$ of degree $4$ is 
$$T^4 + (2\cdot 3 + 2 \cdot 3^2 + 2 \cdot 3^5 + O(3^6))T^3  +  (2 \cdot 3^3 + 2 \cdot 3^5 + O(3^6))T^2 + (3^3 + 3^4 + 3^5 + O(3^6))T + 3^4 + 2 \cdot 3^5 + O(3^6). $$
Rescaling gives the polynomial $\tilde{f}(T)=p^{-4}f(pT)$ having slopes $(0,0,0,0)$. The reduction of $\tilde{f}$ modulo $p$ is 
$$ T^4 + 2 T^3 + T + 1 \in \mathbb{F}_3[T],$$
which is irreducible. Then
 $f(T)$ is also irreducible, and the extension is unramified.

\end{example}

\subsection{Optimizations for high conductors} \label{subsec:high_conductor_speedups}

For conductors $>10^{10}$, we used analytic methods similar to those of \cite{ellenberg2011modeling} as mentioned above. However, the special values required for interpolation are still quite time-consuming to compute, so it is important to avoid computing too many special values for a given character unless it is really necessary. Hence we only computed data on two statistics, the discriminant conditional on $\lambda=2$ (\Cref{tab:lambda2-discriminant-low-high-comparison}), and the number of roots in a small ball $p^3\Z_p$ (Figure~\ref{fig:twist1_smallball_high_low}). In each case we used additional optimizations described here. These work for twist $1$, and also for the other `special' twist $1+\frac{p-1}{2}$ (though we did not implement the latter systematically), but they do not work for general twists. Fortuitously, these are exactly the ones where we needed more data, as they are the only ones which disagree with the random matrix predictions at lower conductor.

The basic idea is that one computes the special value $L_p(0,\theta\omega^i)$ at $0$ using the classical $L$-value $L(0,\theta \omega^{i-1})$. This can be done efficiently if $\theta \omega^{i-1}$ is a quadratic character, i.e. in the cases of the two special twists $i=1$ and $1+\frac{p-1}{2}$. Computing the special value at $0$ of a Dirichlet $L$-function corresponding to a quadratic character is much faster than computing other special values. Furthermore, the constant term of the associated Iwasawa series and of its distinguished polynomial are the same as this special value, up to a unit. 

For the average number of zeros in $p^3\Z_p$, we know that if there is a zero in this disk, then the constant term of the distinguished polynomial must have valuation at least $3$. Since the constant term of the distinguished polynomial has the same $p$-adic valuation as the classical $L$-value $L(0,\theta\omega^{i-1})$ for either of these two special twists, we could compute just this one special value, and only compute the rest of the special values needed for interpolation when it had valuation $\geq 3$, which is a very small number of cases. 

For the discriminant, we used only one or two special $L$-values to detect the cases $\lambda = 0$, $\lambda=1$ and $\lambda \geq 2$ (see \cite{knospe2024special}). The former two make up the great majority of cases, and are irrelevant for the discriminant. These preliminary tests eliminate the vast majority of characters before interpolation, making the computations in the conductor windows near
$10^{10}$ feasible.

\bibliographystyle{alpha_abbrvsort}
\bibliography{references.bib}

\end{document}